\documentclass[11pt]{article}
\renewcommand{\theequation}{\thesection.\arabic{equation}}
\date{}
\makeatletter
\@addtoreset{equation}{section}
\makeatother
\makeatother
\usepackage{amssymb,amsmath}
\usepackage{amsmath}
\usepackage{color}
\usepackage[colorlinks,linkcolor=blue,citecolor=blue]{hyperref}
\usepackage{cite}
\usepackage{graphicx}
\usepackage[misc]{ifsym}
\usepackage{enumitem}
\usepackage{caption}
\usepackage{times}

\renewcommand{\baselinestretch}{1}
 \numberwithin{equation}{section}
 \newtheorem{thm}{Theorem}[section]
\newtheorem{lem}[thm]{Lemma}
\newtheorem{rem}[thm]{Remark}
\newtheorem{pro}[thm]{Proposition}
\newtheorem{cor}[thm]{Corollary}
\newtheorem{myDef}[thm]{Definition}

\newcommand{\pf}{{\bf Proof.\ }}
\newcommand{\hbx}{\hfill$\Box$}
\newcommand{\Rn}{{\mathbb{R}^n}}

\newcommand{\W}{W^{1,2}(\Omega)}

\newcommand{\WW}{W^{1,2}(\Omega)\times W^{1,2}(\Omega)}
\newcommand{\WWZ}{W^{1,2}_0(\Omega)\times W^{1,2}_0(\Omega)}
\newcommand{\K}{\mathcal{K}}
\newcommand{\KW}{\widetilde{\mathcal{K}}}

\newcommand{\dist}{{\rm dist}}
\newcommand{\supp}{{\rm supp }}
\newcommand{\loc}{{\rm loc}}
\newcommand{\Int}{{\rm Int}}
\newcommand{\diam}{{\rm diam}}
\newcommand{\tr}{{\rm tr}}

\allowdisplaybreaks

\begin{document}

\title{\Large On a class of coupled obstacle systems}
\author{\small  Lili Du$^{1,2}$, Xu Tang$^2$, Cong Wang$^3$\footnote{Corresponding author. \  \hfill\break\indent  E-mail:  dulili@szu.edu.cn(L. Du), tangxu8988@163.com(X. Tang), wc252015@163.com(C. Wang).}\\
{\small $^{1}$School of Mathematical Sciences, Shenzhen University}\\ {\small Shenzhen, 518060, People's Republic of China}\\
{\small $^{2}$Department of Mathematics, Sichuan University}\\ {\small Chengdu 610064, People's Republic of China}\\
{\small  $^{3}$School of  Mathematics, Southwest Jiaotong University}\\ {\small Chengdu  610031, People's Republic of China}}

\maketitle
\begin{abstract}
In this paper, we explore cooperative and competitive coupled  obstacle systems, which, up to now, are new type obstacle systems and formed by coupling two equations belonging to classical obstacle problem. On one hand, applying the constrained minimizer in variational methods we establish the existence of solutions for the systems. Moreover, the optimal regularity of solutions is obtained, which is the cornerstone for further research on so-called free boundary.  Furthermore, as coefficient $\lambda\to0$, there exists a sequence of solutions converging to solutions of the single classical obstacle equation.  On the other hand, motivated by the heartstirring ideas of single classical obstacle problem, based on the corresponding blowup  methods, Weiss type monotonicity formula and  Monneau type monotonicity formula of systems to be studied, we investigate the regularity of free boundary, and on the regular and singular points in particular, as it should be, which is more challenging but exceedingly meaningful in solving free boundary problems.
 \end{abstract}

\renewcommand{\contentsname}{\centering Contents}
{\renewcommand{\baselinestretch}{0.8} \tableofcontents }

\section{Introduction}
In present paper, we consider cooperative($\lambda>0$) and competitive($\lambda<0$) obstacle system,
\begin{eqnarray}
\left\{
  \begin{array}{ll}
\Delta u=(1+\lambda v)\chi_{\{u>0\}},\ &{\rm in}\ \Omega,\\
\Delta v=(1+\lambda u)\chi_{\{v>0\}},\ &{\rm in}\ \Omega,\\
u, v\geq0, \ &{\rm in}\ \Omega,\\
u=g_1, v=g_2,\ &{\rm on}\ \partial\Omega,
  \end{array}
\right.\label{system}
\end{eqnarray}
where $\Omega$ is a smooth open domain in $\mathbb{R}^n$, $(u,v)\in \W\times\W, (g_1, g_2)\in \W\times \W, g_1,g_2\geq0$ on $\partial\Omega$ in the sense of $(g_1^-, g_2^-)\in W^{1,2}_0(\Omega)\times W^{1,2}_0(\Omega)$  and $g_1,g_2\not\equiv0$, where $g_i^{\pm}=\max\{\pm g_i,0\}, \lambda\in \mathbb{R}$. System \eqref{system} as $\lambda=0$ reduces to the classical obstacle problem. It is worth noting that the boundary conditions above exclude the possibility of semi-trivial solution $(u,0)$ or $(0,v)$ of problem \eqref{system}. The corresponding free boundary of system \eqref{system} is  $\Gamma(u,v)=\Gamma(u)\cup\Gamma(v)$, where
\begin{eqnarray*}
\Gamma(u)=\Omega\cap \partial\Omega(u)=\Omega\cap \partial\Lambda(u), \Gamma(v)=\Omega\cap \partial\Omega(v)=\Omega\cap\partial\Lambda(v)
\end{eqnarray*}
with $\Omega(u)=\{u>0\},\Lambda(u)=\{u=0\},\Omega(v)=\{v>0\},\Lambda(v)=\{v=0\}$.

Taking $\lambda=0$ in the system \eqref{system}, it reduces to the classical obstacle problem
\begin{eqnarray}\label{classicalP}
\Delta u=\chi_{\{u>0\}}\ {\rm and}\ u\geq0\  {\rm in}\ \Omega,
\end{eqnarray}
which has always been a hot research topic in the field of free boundary problems. Two natural and interesting questions arise at this point:
\begin{description}
  \item[${\rm (Q_1)}$] How regular is the solution of problem \eqref{classicalP}?
  \item[${\rm (Q_2)}$]  Is the free boundary a smooth surface of problem \eqref{classicalP}?
\end{description}

Concerning the  question ${\rm (Q_1)}$, it was shown in  \cite{1972Frehse} that the solution is in $C^{1,1}(\Omega)$ and which is optimal. On the
question ${\rm (Q_2)}$, the first fundamental accomplishment is due to  Kinderlehrer and  Nirenberg \cite{1977KinderlehrerNirenberg}, who found that  $C^1$ free boundaries are  analytic. The breakthrough that removing the assumption of smoothness $C^1$ came from the seminal paper \cite{1977Caffarelli}, where Caffarelli  introduced one of the most transformative ideas in the theory of free boundary problems, namely the use of blowup arguments. Moreover, the highly nontrivial and deep issues that regular and singular points should exhaust the whole free boundary have been answered by  Caffarelli \cite{1977Caffarelli}, this is the so-called Caffarelli's dichotomy theorem. Combining with \cite{1977Caffarelli,1998Caffarelli,1977KinderlehrerNirenberg}, using the directional monotonicity property, allows one to conclude that the free boundary is an analytic hypersurface in a neighborhood of a regular point. In \cite{1998Caffarelli}, the uniqueness about blowup limit is accomplished by applying the celebrated Alt-Caffarelli-Friedman monotonicity formula to the first derivatives of the solution. Caffarelli \cite{1998Caffarelli} established  the first decisive result that  concerning the structure of the singular set. Namely, he
proved that  the singular set is locally contained in a $C^1$ manifold.  Weiss \cite{1999Weiss} introduced an main tools in the study of free boundary problems of a variational nature with a homogeneous structure.  Monneau devised in \cite{2003Monneau} a monotonicity formula that nowadays bears his name based on the Weiss's monotonicity formula, applying Monneau's  monotonicity formula, one infers not only uniqueness of the blowup, but also the continuous dependence from the singular point. In the two-dimensional space, Weiss \cite{1999Weiss} showed  that $C^1$--regularity for the structure of the singular set can be improved to $C^{1,\alpha}$ by means of an epiperimetric-type approach. Very recently, Figalli and  Serra \cite{2019FigalliSerra} obtained an exhaustive characterization of the singular set:  singular set is locally contained in a $C^2$ curve in two dimensions, whereas  singular points(except anomalous ones) can be covered by $C^{1,1}$--and in some cases $C^2$--manifolds in three dimensional and higher dimensional space.  For more classical results, see \cite{1984AltCaffarelliFriedman,1980Caffarelli,2018Figalli,2012PetrosyanShahgholianUraltseva,
1984FriedmanPhillips,1985Kawohl,2003Monneau,1964Stampacchia}, there lists are, of course, far from being exhaustive. The motivation for studying obstacle problems in general, not only because of its profound physical backgrounds and extensive application in elasticity, biology, engineering and mathematical finance, but also thanks to the promotion and enlightening of mathematical theories generated by researching classical obstacle problem to other branch of mathematics(such as potential theory, complex analysis, geophysics) and other types of obstacle problems(such as obstacle problems involving $p$-Laplacian, higher order operator, fractional order operator, and thin obstacle problems, vector-valued types obstacle problems or obstacle systems).

As a matter of research on fundamental mathematical theory,  system \eqref{system} may be seem as one of the extensions of the classical obstacle problem to the vector-valued case by weakly coupling two classical scalar equations \eqref{classicalP}. This coupling is achieved through the presence of the terms $\lambda v\chi_{\{u>0\}}$ and $\lambda u\chi_{\{v>0\}}$.

The one of the important backgrounds of such systems, from a biological standpoint, is the corresponding reaction-diffusion system,
\begin{eqnarray*}
\left\{
  \begin{array}{ll}
u_t=-\Delta u+1+\lambda v,\\
v_t=-\Delta v+1+\lambda u,
  \end{array}
\right.
\end{eqnarray*}
from which it is means that, recording the concentrations $u$ and $v$ of two species/reactants, each species/reactant slows down or accelerates the extinction/reaction of the other species. The proper selection of reaction kinetics would assure a constantdecay/reaction rate in the case that $u$ and $v$ are of comparable size. From a physical perspective,  system \eqref{system} also arises from, in special circumstances, Bose-Einstein condensates(see \cite{2004DKNF,1996MCSS,2003RCFGKMWHV}).

In recent years, there has been increasing attention towards obstacle type systems. Andersson, Shahgholian, Uraltseva and  Weiss \cite{2015AnderssonShahgholianUraltsevaWeiss} considered the singular system
\begin{eqnarray*}
\Delta \textbf{U}=\frac{\textbf{U}}{|\textbf{U}|}\chi_{\{|\textbf{U}|>0\}},
\end{eqnarray*}
where $\textbf{U}=(u_1,u_2,\cdots, u_m)\in W^{1,2}(\Omega,\mathbb{R}^m)$ is a vector function: $\Omega\to \mathbb{R}^m$. The authors were interested in qualitative behaviour of the minimisers $\textbf{U}$ as well as the free boundary $\partial\{x:| \textbf{U}(x)|>0\}$. They noted that the part of the free boundary where the gradient $\nabla \textbf{U}\not=0$ is, by the implicit function theorem, locally a $C^{1,\beta}$-surface, so that they were more concerned with the part where the gradient vanishes. The main result of the paper stated that the set of ``regular'' free boundary points of the minimizers $\textbf{U}$ is locally a $C^{1,\beta}$-surface. In order to proving this result, the authors employed various  technical tools including monotonicity formulas, quadratic growth of solutions, and an epiperimetric inequality for the balanced energy functional.

Figalli and  Shahgholian \cite{2015FigalliShahgholian} provided a comprehensive introduction for the unconstrained free boundary problems, including the unconstrained system, and some very meaningful open questions in the field. El Hajj and Shahgholian in \cite{21ElHajjShahgholian} focused on, the interesting is whether the super-level sets of solutions inherit certain geometric properties of given domain, the convexity of free boundaries for the Bernoulli and obstacle type problem, both in scalar and system cases. They presented heuristic ideas for several free boundary problems, emphasizing the persisting allure of investigating the convexity of the obstacle system. They proposed some intriguing open problems and expressed their belief that  these problems will be of interest for the future study.  Additional information and open questions related to this topic were brought up in \cite{21ElHajjShahgholian}, what this to manifest is the significance of studying  obstacle problems \eqref{system}. The works presented in \cite{1992Adams,2000Adams,1985ChipotVergara-Caffarelli,1986DuzaarFuchs,1979EvansFriedman,1971Lions} have shed considerable light on various aspects related to optimal switching, multi-membranes, control of systems, constrained weakly systems, vector-valued types obstacle systems.

The Bernoulli type cooperative free boundary systems have been studied widely. In the work of Caffarelli, Shahgholian and Yeressian \cite{2018CaffarelliShahgholianYeressian}, they studied the Bernoulli cooperative systems, they proved that the regularity of the free boundary and almost minimizers are Lipschitz continuous. Mazzoleni, Terracini and Velichkov \cite{2020MazzoleniTerraciniVelichkov} considered the regularity of the free boundary for the Bernoulli vectorial minimization problem, they proved that the free boundary, the regular part is relatively open and locally the graph of $C^\infty$ function, the singular part is relative closed and has finite $\mathcal{H}^{n-1}$  measure. Mazzoleni, Terracini and Velichkov \cite{2017MazzoleniTerraciniVelichkov} extended the regularity of free boundary for the one-phase free boundary problems to the vector valued case by studying regularity properties of optimal sets for the shape optimization problem. Spolaor,  Velichkov in \cite{2019SpolaorVelichkov} established the $C^{1,\alpha}$ regularity of the free boundary  in the two-dimensional case for the classical one-phase problems in the  vector cases without any restriction on the sign of the component functions. Aleksanyan,  Fotouhi, Shahgholian and  Weiss \cite{2022AFSW} conducted a thorough investigation into the parabolic cooperative free boundary system. They demonstrated that an optimal growth rate for solutions near free boundary points and the regularity of the free boundary.

Our focus will be on the obstacle system \eqref{system}. And,  as we said, no results for system \eqref{system} have been obtained thus far.
\subsubsection*{The main result and structure of the paper}

Some notations and  preknowledges will be introduced in Section \ref{secpre}. In Section \ref{sec2}, based the constrained minimizer, we prove the existence of solutions for obstacle system \eqref{system}. The results are as follows.
\begin{thm}[Existence]\label{thm3.2}
		There exists $(u,v)\in \WW$ solving system \eqref{system} for any $\lambda>-\lambda_1$, where $\lambda_1$ is the first eigenvalue of operator $-\Delta$ in $W^{1,2}_0(\Omega)$.
\end{thm}
\begin{rem}\label{remeigen}
 Indeed, by the variational methods and compact embedding theorems, it is not difficult to know that there exists a first eigenvalue $\lambda_1>0$ of operator $-\Delta$ equipping with  the corresponding eigenfunction $\varphi_1$ in $W^{1,2}_0(\Omega)$ with nonnegative definiteness,  which satisfies
	\begin{eqnarray}
		\displaystyle \lambda_1=\frac{\int_\Omega|\nabla \varphi_1|^2dx}{\int_\Omega|\varphi_1|^2dx}=\inf_{\varphi\in W^{1,2}_0(\Omega)\setminus\{0\}}\frac{\int_\Omega|\nabla \varphi|^2dx}{\int_\Omega|\varphi|^2dx}\label{33.1}
	\end{eqnarray}
	and the eigenvalue problem
	\begin{eqnarray*}
		\left\{
		\begin{array}{ll}
			-\Delta \varphi_1=\lambda_1 \varphi_1\ {\rm in}\ \Omega,\\[0.5em]
			\varphi_1\in W^{1,2}_0(\Omega)\setminus\{0\}.
		\end{array}
		\right.
	\end{eqnarray*}
\end{rem}

\begin{rem}\label{rem1.2}
Theorem {\rm\ref{thm3.2}} is proved by using approximation method and which shows that the minimizer of $J$ over $\K$ is exactly solution of cooperative($\lambda>0$) and competitive($\lambda<0$) obstacle problems. The restriction $\lambda>-\lambda_1$ is to ensure that the boundedness from below  functional of
 \begin{eqnarray*}
		J(u,v)=\frac{1}{2}\int_\Omega |\nabla u|^2+|\nabla v|^2dx+\int_\Omega u+v+\lambda uv dx
	\end{eqnarray*}
and the corresponding functional of the regularity problem(see \eqref{system1}), as well as the uniform boundedness of approximation solutions sequence.
	
	If $\lambda<-\lambda_1$, as $\|\nabla g_1^+\|_{L^2(\Omega)},\|\nabla g_2^+\|_{L^2(\Omega)}, \|g_1^+\|_{L^2(\Omega)}, \|g_2^+\|_{L^2(\Omega)}$ small enough such that
	\begin{eqnarray*}
		\bar{C}:&=&(\lambda_1+\lambda)\|\varphi_1\|_{L^2(\Omega)}^2+\frac12\|\nabla g_1^+\|^2 +\frac12\|\nabla g_2^+\|^2   \\
		&&+\|\varphi_1\|_{L^2(\Omega)}\left(\lambda_1\|\nabla g_1^+\|_{L^2(\Omega)}+\lambda_1\|\nabla g_2^+\|_{L^2(\Omega)}+\|g_1^+\|_{L^2(\Omega)}+\|g_2^+\|_{L^2(\Omega)}\right)\\
		&<&0,
	\end{eqnarray*}
	taking $(u,v)=\mu(\varphi_1+g_1^+,\varphi_1+g_2^+), \mu>0$  such that $\mu(\varphi_1+g_1^+,\varphi_1+g_2^+)\in \mathcal{K}$. Then, using H\"{o}lder inequality,
	\begin{eqnarray*}
		&&J(\mu(\varphi_1+g_1^+),\mu(\varphi_1+g_2^+))\\
		&=&\frac{\mu^2}{2}\int_\Omega |\nabla (\varphi_1+g_1^+)|^2+ |\nabla (\varphi_1+g_2^+)|^2dx+\lambda \mu^2\int_\Omega (\varphi_1+g_1^+)(\varphi_1+g_2^+) dx\\
		&&\ \ +\mu\int_\Omega 2\varphi_1+g_1^++g_2^+dx\\
		&\leq&\bar{C}\mu^2+\mu\int_\Omega 2\varphi_1+g_1^++g_2^+dx,
	\end{eqnarray*}
	which implies that $J(\mu(\varphi_1+g_1^+),\mu(\varphi_1+g_2^+))\to-\infty$ as $\mu\to+\infty$. And the functional $J$ is unbounded from below for $\lambda<-\lambda_1$.
\end{rem}

By the standard bootstrap argument, we know that $(u,v)\in C^{1,\alpha}_{\rm loc}(\Omega)\times C^{1,\beta}_{\rm loc}(\Omega)$ with any  $\alpha,\beta\in(0,1)$. The optimal regularity can be found in Section \ref{sec4} using the bootstrap argument and  the proof of \cite[Theorem 2.14]{2012PetrosyanShahgholianUraltseva}.
\begin{thm}[Regularity]\label{thmopimalregularity}
		Let $(u,v)$ be the solution of system \eqref{system}. Then $(u,v)\in C_{\rm loc}^{1,1}(\Omega)\times  C_{\rm loc}^{1,1}(\Omega)$ and
		\begin{eqnarray*}
			\|u\|_{C^{1,1}(K)}\leq C\left(1+\|v\|_{L^\infty(K)}+\|\nabla v\|_{L^\infty(K_{2\delta})}+\|\nabla u\|_{L^\infty(K_\delta)}\right),\\
			\|v\|_{C^{1,1}(K)}\leq C\left(1+\|u\|_{L^\infty(K)}+\|\nabla u\|_{L^\infty(K_{2\delta})}+\|\nabla v\|_{L^\infty(K_\delta)}\right)
		\end{eqnarray*}
		for any compact set $K\subset \Omega$, where $C=C(n,\delta,\lambda)$ and $\delta:=\frac{1}{3}\dist\ (K,\partial\Omega)$, $K_\delta:=\{x:\dist(x,K)<\delta\}$.
\end{thm}

\begin{rem}\label{ren1.4}
It is easy to see that this is the best possible result we could hope for as $\lambda>0$ or under   the assumptions
\begin{align}
\lambda<0,\ 1+\lambda v(z_0)\not=0, \  \exists\ z_0\in\Gamma(u)\tag{$\Gamma_u$}.\label{gammau}\\
\lambda<0,\ 1+\lambda u(z_0)\not=0, \  \exists\ z_0\in\Gamma(v)\tag{$\Gamma_v$}.\label{gammav}
\end{align}
When across the boundary of the contact set $\{u=0\}$ the Laplacian of $u$  ``jump'' from  $0$ to $1+\lambda v$. It follows from $\lambda>0$  that $1+\lambda v>0$ on $\Gamma(u)$ and from \eqref{gammau} that $1+\lambda v\not=0$ on certain neighborhood of $z_0$.  So $u$  is not necessarily $C^2$ in $\Omega$. Using the same analysis the optimal of $v$ can be  established.
\end{rem}

In Section \ref{secbehavior}, we consider the asymptotic behavior of solutions of obstacle system \eqref{system} as $\lambda\to0$.  For $\lambda\to0$, we get that the asymptotic behavior of solutions of system \eqref{system}.
	\begin{thm}[Convergence]\label{thmasy}
		Assume that $(u_\lambda,v_\lambda)$ solves obstacle problem \eqref{system} for any $\lambda>-\lambda_1$. As $\lambda\to0$, then there exist a sequence of solutions $\{(u_{\lambda_n}, v_{\lambda_n})\}$ with $\lim_{n\to\infty}\lambda_n=0$ and $(u_0, v_0)\in \WW$ such that $u_{\lambda_n}\to u_0, v_{\lambda_n}\to v_0$  in $W^{1,2}(\Omega)$, where $(u_0, v_0)$ is the solution of classical obstacle problem
		\begin{eqnarray}
			\left\{
			\begin{array}{ll}
				\Delta u= \chi_{\{u>0\}},\ {\rm in}\ \Omega,\\
				\Delta v= \chi_{\{v>0\}},\ {\rm in}\ \Omega,\\
                u,v\geq0,\ {\rm in}\ \Omega,\\
				u=g_1, v=g_2,\ {\rm on}\ \partial\Omega.
			\end{array}
			\right.\label{classicalobstacle}
		\end{eqnarray}
	\end{thm}

In Section \ref{sec5}, we main focus is on the regularity and properties of the free boundaries which may encounter multiple cases. It is possible that free boundaries $\Gamma(u), \Gamma(v)$ of obstacle system \eqref{system} may intersect with each other(see the domain $D$ in Figure \ref{fig1} or Figure \ref{fig3}), or rather, we cannot rule out this case, which may bring  more difficulties. It is also possible that $\Gamma(u)$ is tangent to $\Gamma(v)$ as shown in Figure \ref{fig2}. Even the most anticipated situation may arise, that is $\Gamma(u)$ is away from $\Gamma(v)$ as shown in Figure \ref{fig4}.\\
 \begin{minipage}{\textwidth}
\centering
\includegraphics[width=8cm,height=3cm]{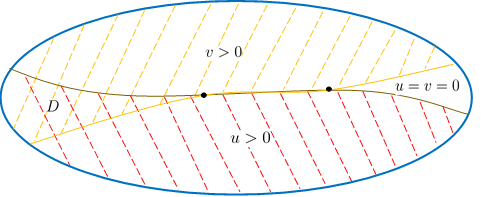}
\captionof{figure}{$\Omega(u)$ intersects $\Omega(v)$ in $D$}\label{fig1}
\end{minipage}
\begin{minipage}{\textwidth}
\centering
\includegraphics[width=8cm,height=3.5cm]{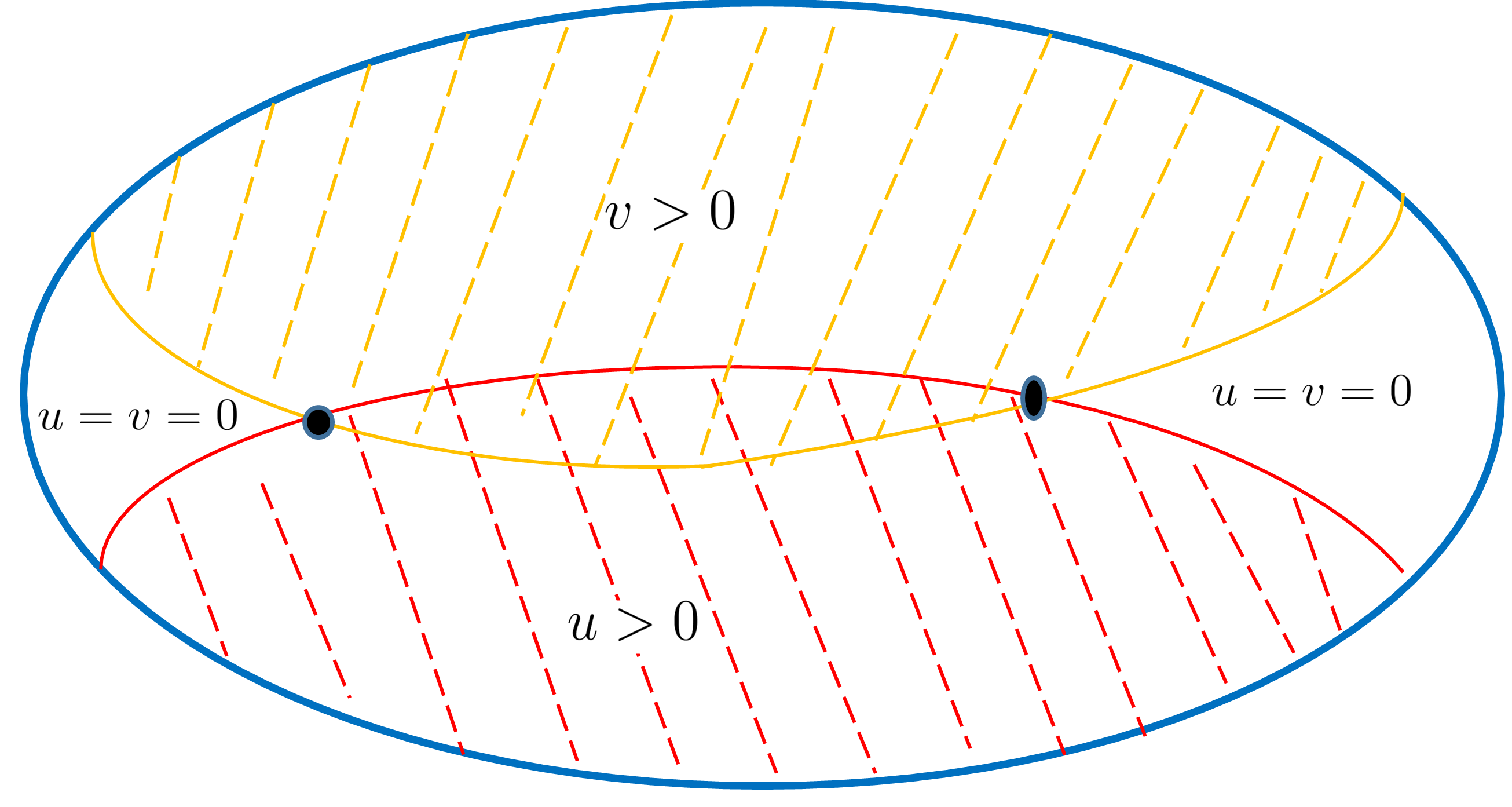}
\captionof{figure}{$\Omega(u)$ intersects $\Omega(v)$}\label{fig3}
\end{minipage}
\begin{minipage}{\textwidth}
\centering
\includegraphics[width=8cm,height=3cm]{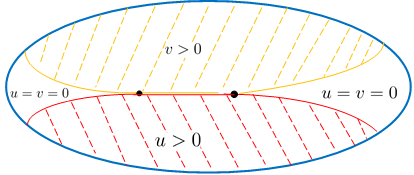}
\captionof{figure}{$\Omega(u)$ is tangent to  $\Omega(v)$}\label{fig2}
\end{minipage}
\begin{minipage}{\textwidth}
\centering
\includegraphics[width=8cm,height=3cm]{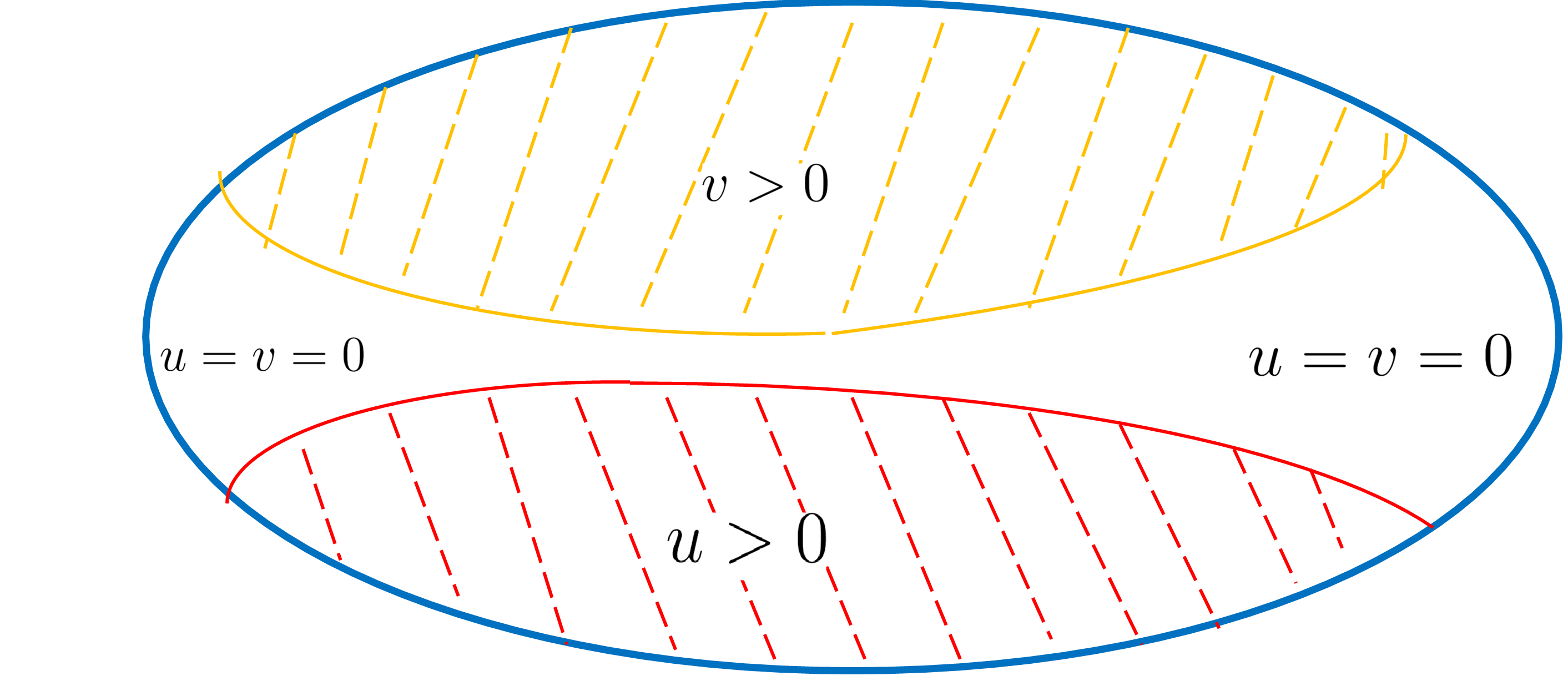}
\captionof{figure}{$\Omega(u)$ is away from $\Omega(v)$}\label{fig4}
\end{minipage}

In Subsection \ref{sec5.1}, the properties are established such as nondegeneracy of solutons, porosity, Lebesgue measure and Hausdorff measure of the free boundary $\Gamma(u,v)$, etc. It is noteworthy that the nondegeneracy of solutions plays an important role in proving the porosity, Lebesgue measure, Hausdorff measure of free boundary and the subsequent application of blowup analysis to study the regularity of free boundaries. In order ro establish the nondegeneracy of solutions for the case of $\lambda<0$, we impose two \textbf{assumptions}:
\begin{align}
\lambda<0,\ 1+\lambda v(z)>0\ {\rm for\ any}\ z\in\Gamma(u)\tag{$\Gamma_{u,\lambda}$},\label{gammaulambda}\\
\lambda<0,\ 1+\lambda u(z)>0\ {\rm for\ any}\ z\in\Gamma(v)\tag{$\Gamma_{v,\lambda}$}.\label{gammavlambda}
\end{align}
\begin{rem}\label{rem1.6}
The maximum principle is a very momentous tool in establishing the nondegeneracy of solutions. Obviously, $\lambda>0$ can make solutions of system \eqref{system} utilize the  maximum principle.  The case $\lambda<0$  seems to be unable to  directly apply the  maximum principle. However, the assumptions \eqref{gammaulambda} and \eqref{gammaulambda} get rid of such doubts effectively.

In fact, the cases of $\Gamma(u)$ is tangent to $\Gamma(v)${\rm(}see Figure {\rm\ref{fig2})} and $\Gamma(u)$ is away from $\Gamma(v)${\rm(}see Figure {\rm\ref{fig4})}  satisfy  clearly the assumptions \eqref{gammaulambda} and \eqref{gammaulambda}. But the other cases{\rm(}such as Figure {\rm\ref{fig1}} or Figure {\rm\ref{fig3})}, as mentioned above, may require assumptions \eqref{gammaulambda} and \eqref{gammaulambda} to apply the maximum principle.
\end{rem}

Moreover,  in Subsection \ref{sec5.1},  we consider the blowups of solutions on the free boundary, the unique type of blowup. Furthermore, we examine the properties of free boundary divided into, as described following, singular points and regular points.

\begin{myDef}\label{def1.7}

If $x^0\in\Gamma(u)$:
\begin{itemize}
  \item $x^0$ is a regular point if one and therefore every blowup of $u$ at $x^0$ is a half-space solution $\frac{1+\lambda v(x^0)}{2}[(x\cdot e)^+]^2, x\in \Rn$, where $e$ is a unit vector.
  \item $x^0$ is a singular point if one and therefore every blowup of $u$ at $x^0$ is polynomial solution $\frac{1+\lambda v(x^0)}{2}(x\cdot Ax), x\in \Rn$, where $A$ is an $n\times n$ symmetric matrix with ${\rm tr}A=1$.
\end{itemize}
If $x^0\in\Gamma(v)$:
\begin{itemize}
  \item $x^0$ is a regular point if one and therefore every blowup of $u$ at $x^0$ is a half-space solution $\frac{1+\lambda u(x^0)}{2}[(x\cdot e)^+]^2, x\in \Rn$, where $e$ is a unit vector.
  \item $x^0$ is a singular point if one and therefore every blowup of $u$ at $x^0$ is polynomial solution $\frac{1+\lambda u(x^0)}{2}(x\cdot Ax), x\in \Rn$, where $A$ is an $n\times n$ symmetric matrix with ${\rm tr}A=1$.
\end{itemize}

\end{myDef}

 The regularity of free boundary near regular points is considered in Subsection \ref{subsec5.4}, the results are following.
\begin{thm}[Regularity of free boundaries]\label{thmC1Alpha}
Let $(u,v)$ be a solution of system \eqref{system} satisfying $0$ is a regular point.
\begin{description}
\item[(i)] If  $0\in\Gamma(u)$,  $\lambda\geq0$ or \eqref{gammaulambda} holds.  Then there exists $\rho>0$ such that $\Gamma(u)\cap B_\rho$ is a $C^{1,\alpha}$ graph.
\item[(ii)] If  $0\in\Gamma(v)$,  $\lambda\geq0$ or \eqref{gammavlambda} holds.  Then there exists $\rho>0$ such that $\Gamma(v)\cap B_\rho$ is a $C^{1,\alpha}$ graph.
\end{description}
\end{thm}

The structure of  the singular set  can be found in Subsection \ref{subsec5.5}, we obtain the following results.
	\begin{thm}[Singular set of free boundaries]\label{thm9.7}
		Let $(u,v)$ be a solution of system \eqref{system} in an open set $\Omega$ in $\Rn$. Then
		\begin{itemize}
			\item As $\lambda\geq0$ or \eqref{gammaulambda} holds. For any $d\in \{0,1,\cdots,n-1\}$, the set
			\begin{eqnarray*}
				\Sigma^d(u)=\{x^0\in\Sigma(u),d_{x^0}(u)=d\}
			\end{eqnarray*}
			is contained in a countable union of $d$-dimensional $C^1$ manifold.
			\item As $\lambda\geq0$ or \eqref{gammavlambda}  holds.  For any $d\in \{0,1,\cdots,n-1\}$, the set
			\begin{eqnarray*}
				\Sigma^d(v)=\{x^0\in\Sigma(v),d_{x^0}(v)=d\}
			\end{eqnarray*}
			is contained in a countable union of $d$-dimensional $C^1$ manifold.
		\end{itemize}
	\end{thm}

In Appendix, we establish the  properties of the set $\mathcal{K}$, the proofs of Proposition \ref{lemHausdorff} about Hausdorff measure of free boundaries, Proposition \ref{prolimit} with regard to rescaling of solutions, Theorem  \ref{weissmonotonocity} concerning Weiss type monotonicity formulas, Theorem  \ref{thm9.4} involving Monneau type monotonicity formulas.

 A more general form of system \eqref{system} is following
\begin{eqnarray}
\left\{
  \begin{array}{ll}
\Delta \textbf{U}=(I+\lambda \textbf{U}\textbf{A})^T\boldsymbol{\chi}_{\{\textbf{U}>0\}},\ &{\rm in}\ \Omega,\\
\textbf{U}\geq0, \ &{\rm in}\ \Omega,\\
\textbf{U}=\textbf{g},\ &{\rm on}\ \partial\Omega,
  \end{array}
\right.\label{systemMul}
\end{eqnarray}
where $\textbf{U}=(u_1,u_2,\cdots, u_m)\in W^{1,2}(\Omega,\mathbb{R}^m)$ is a vector function: $\Omega\to \mathbb{R}^m$. The vector $I=(1,1,\cdots,1)^T$ is composed of $m$ number one. The $m\times m$ matrices $\textbf{A}, \boldsymbol{\chi}_{\{\textbf{U}>0\}}$ are
\begin{eqnarray*}
\textbf{A}:=\left(
            \begin{array}{cccccc}
            0 &   0&  0& \cdots&0 &1 \\
            1 &   0& 0&\cdots &0 & 0\\
            0 &   1& 0& \cdots & 0& 0\\
            0&    0& 1 &\cdots &0 &0\\
            \vdots& \vdots& \vdots&  \ddots&\vdots & \vdots\\
            0& 0 & 0 &   \cdots& 1& 0
            \end{array}
          \right),\ \
\boldsymbol{\chi}_{\{\textbf{U}>0\}}:=\left(
            \begin{array}{ccccc}
            \chi_{\{u_1>0\}} &   0&   \cdots&0 \\
            0 &    \chi_{\{u_2>0\}}& \cdots  & 0\\
            \vdots& \vdots&   \ddots& \vdots\\
            0& 0 &    \cdots&  \chi_{\{u_m>0\}}
            \end{array}
          \right).
\end{eqnarray*}
The expressions $\textbf{U}\geq0$ and $\textbf{U}=\textbf{g}$ stand for  $u_i\geq0$ and $u_i=g_i$  respectively for all $i=1,2,\cdots,m$ with $\textbf{g}=(g_1,g_2,\cdots, g_m)\geq0$ on $\partial\Omega$ and $\textbf{g}\in W^{1,2}(\Omega,\mathbb{R}^m)$. The results of system \eqref{system} established in the present paper can be extended to the multiple obstacle system \eqref{systemMul}.

	\section{Preliminaries}\label{secpre}
	We collect some notations and  pre-knowledges, stating the results that will be used in the following sections. We first define some notations.
	\begin{eqnarray*}
		&&\|u\|_{L^p(\Omega)}:=\left(\int_\Omega|u|^pdx\right)^{1/p}\ {\rm with}\ 1\leq p<\infty.\\
		&&\displaystyle\|u\|_{L^\infty(\Omega)}:=\mathop{\rm ess\ sup}\limits_{x\in\Omega}|u(x)|.\\
		&&u^{\pm}:=\max\{\pm u,0\}.\\
		&&e_i:=(0,\cdots, 1, \cdots,0)\ {\rm with}\ 1\ {\rm in\ the}\ i{\rm{-th}}\ {\rm position}.\\
		&&B_r(x_0):=\{x\in \Rn| |x-x_0|<r\},\ B_r:=B_r(0).\\
		&&\omega_n:=\frac{\pi^{n/2}}{\Gamma(\frac{n}{2}+1)} {\rm is\ the\ volume\ of\ unit\ ball}\ B_1\ {\rm in}\ \Rn\ {\rm with\ Gamma\ function}\ \Gamma.\\
		&&{\rm diam}(E):=\sup_{x,y\in E}|x-y|.\\
		&&\mathcal{C}_\delta:=\{x\in \Rn:x_n>\delta|x'|\},\ x'=(x_1,\cdots,x_{n-1}).\\
		&&K(\delta,s,h):=\left\{|x'|<\delta,\ -s\leq x_n\leq h\right\}\ {\rm for}\ \delta, s, h>0.\\
		&&\Sigma(u), \Sigma(v)\ {\rm is\ the\ set\ of\ all\ singular\ points\ of}\ u\ {\rm and}\ v\ {\rm repectively}.\\
		&& \nu\ {\rm is\ the\ outer\ normal\ on\ a\ given\ surface}.
	\end{eqnarray*}

	An crucial observation of the Alt-Caffarelli-Friedman(Abbr. ACF) monotonicity formula originally introduced in \cite{1984AltCaffarelliFriedman}, see also \cite[Theorem 2.4]{2012PetrosyanShahgholianUraltseva} enabling us to obtain easily the following version of ACF monotonicity formula centering on $x^0\in \mathbb{R}^n$. Which will play an important role in establish the regularity of solutions, see the proof of Theorem \ref{thmopimalregularity}.
	\begin{lem}\label{thmACF}
		Let $u_\pm$ be a pair of continuous functions such that
		\begin{eqnarray*}
			u_\pm\geq0,\ \Delta u_\pm\geq0,\ u_+\cdot u_-=0\ {\rm in}\ B_1(x^0).
		\end{eqnarray*}
		Then the functional
		\begin{eqnarray*}
			r\to \Phi(r,u_+, u_-,x^0):=\frac{1}{r^4}\int_{B_r(x^0)}\frac{|\nabla u_+|^2}{|x-x^0|^{n-2}}dx\int_{B_r(x^0)}\frac{|\nabla u_-|^2}{|x-x^0|^{n-2}}dx
		\end{eqnarray*}
		is nondecreasing for $0<r<1$.
	\end{lem}
We introduce the definition of Porosity, which allows us to demonstrate the Lebesgue measure  and Hausdorff dimension of free boundaries.
\begin{myDef}{\rm \cite[Porosity]{2012PetrosyanShahgholianUraltseva}}\label{defporosity}
		We say that a measurable set $E\subset \Rn$ is porous with a porosity constant $0<\delta<1$ if every ball $B=B_r(x)$ contains a smaller ball $B'=B_{\delta r}(y)$ such that
\begin{eqnarray*}
B_{\delta r}(y)\subset B_r(x)\setminus E.
\end{eqnarray*}
	\end{myDef}
	\begin{myDef}{\rm \cite[Minimal diameter]{2012PetrosyanShahgholianUraltseva}}\label{def}
		The minimal diameter of bounded set $E\subset \Rn$, denoted $\min \diam(E)$, is the infimum of distances between pairs of parallel planes enclosing $E$.
	\end{myDef}
According the definition of minimal diameter above, we introduce the definition of thickness function which allows us to give a sufficient condition to establish the uniform results about regularity of free boundary, furthermore, we can establish the higher regularity of free boundaries near the regular points.
	\begin{myDef}{\rm \cite[Thickness function]{2012PetrosyanShahgholianUraltseva}}\label{def2.3}  For $u$ satisfying $|D^2 u|_{L^\infty(B_R(x^0))}\leq M$ with constant $M$ being independent on radius $R$, the thickness function of $\Omega^c(u)$ at $x^0$ is defined by
		\begin{eqnarray*}
			\delta(\rho,u,x^0)=\frac{\min\diam(\Omega^c(u)\cap B_\rho(x^0))}{\rho},\ 0<\rho<R.
		\end{eqnarray*}
	\end{myDef}

	\section{The existence   of weak solutions}\label{sec2}
	In present section, we focus on the existence  of weak solutions. We find a weak solution of system \eqref{system} via constrained minimization for the  functional
	\begin{eqnarray*}
		J(u,v)=\frac{1}{2}\int_\Omega |\nabla u|^2+|\nabla v|^2dx+\int_\Omega u+v+\lambda uv dx
	\end{eqnarray*}
	over the admissible set
	\begin{eqnarray*}
		\mathcal{K}:=\left\{(u,v)\in\W\times\W, (u-g_1,v-g_2)\in W^{1,2}_0(\Omega)\times W^{1,2}_0(\Omega), u,v\geq0\ {\rm a.e.\ in}\ \Omega\right\}.
	\end{eqnarray*}
	
	Now, we begin to show that the minimizer of $J(u,v)$ over $\mathcal{K}$ solves the system \eqref{system} in the weak sense. The following key lemma plays an important role in the proof of desired results.
	\begin{lem}\label{lem3.3}
		A pair of functions $(u,v)\in\WW$ minimized the functional $J$ over $\mathcal{K}$ if and only if $(u,v)$ is a minimizer of the functional
		\begin{eqnarray*}
			\tilde{J}(u,v)=\tilde{J}(u,v)=\frac{1}{2}\int_\Omega |\nabla u|^2+|\nabla v|^2dx+\int_\Omega u^++v^++\lambda u^+v^+ dx
		\end{eqnarray*}
		over $\widetilde{\mathcal{K}}:=\left\{(u,v)\in\WW, (u-g_1,v-g_2)\in\WWZ\right\}$.
	\end{lem}
	\pf $``\Leftarrow"$ Let $(u,v)\in \widetilde{\mathcal{K}}$ minimize $\tilde{J}$ over $\KW$.
	
	We first claim that $(u,v)=(u^+,v^+)$. It is easy to see that $\tilde{J}(u,v)\geq \tilde{J}(u^+,v^+)$, obviously, $(u^+,v^+)\in \KW$. Since $(u,v)$ is a minimizer of $\tilde{J}$ over $\KW$, then $(u^+,v^+)$ is a minimizer, that is $\tilde{J}(u,v)=\tilde{J}(u^+,v^+)$. Furthermore, $\int_\Omega |\nabla (u,v)|^2dx=\int_\Omega |\nabla (u^+,v^+)|^2dx$ and $\nabla (u,v)=(0,0)$ on $\{u\leq0,v\leq 0\}$, which  is equivalent to $\nabla(u^-,v^-)=(0,0)$ a.e. in $\Omega$. The latter means that $(u^-,v^-)$ is locally constant, and $(u^-,v^-)\in W^{1,2}_0(\Omega)\times W^{1,2}_0(\Omega)$ implies $(u^-,v^-)=(0,0)$. From the above, it can be seen that $(u,v)=(u^+,v^+)$.
	
	If $(w_1, w_2)\in\K$, obviously, $(w_1, w_2)\in\KW$ and $w_1,w_2\geq0$. Then
	\begin{eqnarray*}
		J(u,v)=J(u^+,v^+)=\tilde{J}(u^+,v^+)=\tilde{J}(u,v)\leq \tilde{J}(w_1, w_2)=J(w_1, w_2).
	\end{eqnarray*}
	It follows that $(u,v)$ is a minimizer of $J$ over $\K$.
	
	$``\Rightarrow"$ If $(u,v)\in\K$ is a minimizer of $J$ over $\K$. For any $(w_1,w_2)\in\KW$,
	\begin{eqnarray*}
		\tilde{J}(u,v)=J(u,v)\leq J(w_1^+,w_2^+)=\tilde{J}(w_1^+,w_2^+)\leq \tilde{J}(w_1,w_2).
	\end{eqnarray*}
	Thus $(u,v)$ minimizes $\tilde{J}$ over $\KW$. \hbx

	{\bf Proof of Theorem \ref{thm3.2}.}
	The general form of Young's inequality is
	\begin{eqnarray}
		|a\pm b|^2\leq (1+\epsilon)a^2+C(\epsilon)b^2\ \ {\rm for\ any}\ \epsilon>0, a,b\in \mathbb{R}.\label{3.2}
	\end{eqnarray}
	Based on \eqref{33.1} and  Young's inequality \eqref{3.2}, on the set $\mathcal{K}$,
	\begin{eqnarray}
		\left.
		\begin{array}{ll}
			\displaystyle \int_\Omega uvdx&\leq \frac{1}{2}\int_\Omega|u|^2+|v|^2dx\\[0.8em]
			&\displaystyle = \frac{1}{2}\int_\Omega|u-g_1+g_1|^2+|v-g_2+g_2|^2dx\\[0.8em]
			&\displaystyle \leq \frac{1}{2}\int_\Omega(1+\sigma)|u-g_1|^2+(1+\sigma)|v-g_2|^2+C(\sigma)(|g_1|^2+|g_2|^2)dx\\[0.8em]
			&\displaystyle \leq \frac{(1+\sigma)}{2\lambda_1}\int_\Omega |\nabla(u-g_1)|^2+|\nabla(v-g_2)|^2dx +\frac{C(\sigma)}{2}\int_\Omega (|g_1|^2+|g_2|^2)dx\\[0.8em]
			&\displaystyle  = \frac{(1+\sigma)^2}{2\lambda_1}\int_\Omega (|\nabla u|^2+|\nabla v|^2)+C(\sigma,\lambda_1)\int_\Omega(|\nabla g_1|^2+|g_1|^2+|\nabla g_2|^2+|g_2|^2)dx.
		\end{array}
		\right.\label{3.3-3}
	\end{eqnarray}
	For $\lambda>-\lambda_1$, using \eqref{3.3-3}, taking and fixing $\sigma$ such that $1+\chi_{\{\lambda<0\}}\frac{\lambda(1+\sigma)^2}{\lambda_1}\geq0$, then
	\begin{eqnarray*}
		J(u,v)&\geq&  \frac{1}{2}\left(1+\chi_{\{\lambda<0\}}\frac{\lambda(1+\sigma)^2}{\lambda_1}\right)\int_\Omega |\nabla u|^2+|\nabla v|^2dx\\
		&&+\chi_{\{\lambda<0\}}\lambda C(\sigma,\lambda_1)\int_\Omega(|\nabla g_1|^2+|g_1|^2+|\nabla g_2|^2+|g_2|^2)dx\\
		&\geq&  \chi_{\{\lambda<0\}}\lambda C(\sigma,\lambda_1)\int_\Omega(|\nabla g_1|^2+|g_1|^2+|\nabla g_2|^2+|g_2|^2)dx.
	\end{eqnarray*}
	It follows that $J$ is bounded from below over $\mathcal{K}$ for any $\lambda>-\lambda_1$. It is easy, as we said,  and relative standard to know that there exists a minimizer $(\bar{u},\bar{v})\in \mathcal{K}$ such that $\displaystyle J(\bar{u},\bar{v})=\inf_{(u,v)\in \mathcal{K}}J(u,v)$ combining with Proposition \ref{pro3.1},  cf. \cite[Proposition 2.1]{2018Figalli}.
	
	The following content will mainly focus on proving the solvability of the system \eqref{system} and which is attributed to the minimizer. Considering the minimizers of $\tilde{J}$ over $\KW$ and a family of regularized problems
	\begin{eqnarray}
		\left\{
		\begin{array}{ll}
			\Delta u_\epsilon=(1+\lambda\Phi_\epsilon(v_\epsilon))\chi_\epsilon(u_\epsilon)\ {\rm in}\ \Omega,\\
			\Delta v_\epsilon =(1+\lambda\Phi_\epsilon(u_\epsilon))\chi_\epsilon(v_\epsilon)\ {\rm in}\ \Omega,\\
			u_\epsilon=g_1, v_\epsilon=g_2,\ {\rm on}\ \partial\Omega
		\end{array}
		\right.\label{system1}
	\end{eqnarray}
	for $0<\epsilon<1$, where $\chi_\epsilon(s)$ is a smooth approximation of the Heaviside function $\chi_{\{s>0\}}$ such that
	\begin{eqnarray*}
		\chi'_\epsilon\geq0,\ \chi_\epsilon(s)=0\ {\rm for}\ s\leq-\epsilon,\ \chi_\epsilon(s)=1\ {\rm for}\ s\geq\epsilon.
	\end{eqnarray*}
	A solution $(u_\epsilon,v_\epsilon)$ to the problem \eqref{system1} can be obtained by minimizing the functional
	\begin{eqnarray*}
		J_\epsilon(u,v)=\frac{1}{2}\int_\Omega |\nabla u|^2+|\nabla v|^2dx+\int_\Omega \Phi_\epsilon(u)+\Phi_\epsilon(v)+\lambda\Phi_\epsilon(u)\Phi_\epsilon(v)dx,
	\end{eqnarray*}
	where $\Phi_\epsilon(s)=\int_{-\infty}^s\chi_\epsilon(t) dt\geq0$. It is easy to see that $J_\epsilon\in C^1(W^{1,2}(\Omega)\times W^{1,2}(\Omega), \mathbb{R})$. We show that $J_\epsilon(u,v)$ is bounded from below, using \eqref{3.3-3}, \eqref{3.3-3} and the Sobolev embedding theorem,
	\begin{eqnarray*}
		J_\epsilon(u,v)&=&\frac{1}{2}\int_\Omega |\nabla u|^2+|\nabla v|^2dx+\int_\Omega \Phi_\epsilon(u)+\Phi_\epsilon(v)+\lambda \Phi_\epsilon(u)\Phi_\epsilon(v)dx\\
		&\geq&\frac{1}{2}\int_\Omega |\nabla u|^2+|\nabla v|^2dx+\lambda\chi_{\{\lambda<0\}}\int_\Omega \epsilon(|u|+|v|)+ |u||v|dx+\lambda\chi_{\{\lambda<0\}}C(\Omega)\\
		&\geq&\frac{1}{2}\left(1+\chi_{\{\lambda<0\}}\frac{\lambda(1+\sigma)^2}{\lambda_1}\right)\int_\Omega |\nabla u|^2+|\nabla v|^2dx\\
		&&+\lambda\chi_{\{\lambda<0\}} C(|\nabla u|_{L^2(\Omega)}+|\nabla v|_{L^2(\Omega)})+\lambda\chi_{\{\lambda<0\}}C(\Omega)\\
		&&+\lambda\chi_{\{\lambda<0\}}C(\sigma,\lambda_1)\int_\Omega(|\nabla g_1|^2+|g_1|^2+|\nabla g_2|^2+|g_2|^2)dx\\
		&&+\lambda\chi_{\{\lambda<0\}}\int_\Omega(|g_1|+|g_2|)dx+ \lambda\chi_{\{\lambda<0\}} C(\|\nabla g_1\|_{L^2(\Omega)}+\|\nabla g_2\|_{L^2(\Omega)}) .
	\end{eqnarray*}
	For $\lambda>-\lambda_1$, the desired result that $J_\epsilon(u,v)$ is bounded from below follows from taking and fixing $\sigma$ such that $1+\chi_{\{\lambda<0\}}\frac{\lambda(1+\sigma)^2}{\lambda_1}>0$. Since  $\KW$ is a Banach space, applying Ekeland's variational principle(cf. \cite[Corollary 5.3]{2008struwe}), we can obtain $(u_\epsilon, v_\epsilon)\in\KW$ such that $\displaystyle J_\epsilon(u_\epsilon,v_\epsilon)=\inf_{(u,v)\in\KW}J_\epsilon(u,v)$. To be more precise, the minimizer $(u_\epsilon, v_\epsilon)$ solves the system \eqref{system1} in the weak sense. We claim that $\{(u_\epsilon, v_\epsilon)\}_\epsilon$ is uniformly bounded in $\WW$. Since $(u_\epsilon, v_\epsilon)$ is a weak solution of system \eqref{system1}, then
	\begin{eqnarray}
		\left\{
		\begin{array}{ll}
			\displaystyle 0=\int_\Omega |\nabla(u_\epsilon-g_1)|^2dx+\int_\Omega \nabla g_1\cdot\nabla(u_\epsilon-g_1)+(1+\lambda\Phi_\epsilon(v_\epsilon))\chi_\epsilon(u_\epsilon)(u_\epsilon-g_1)dx,\\[1em]
			\displaystyle 0=\int_\Omega |\nabla(v_\epsilon-g_2)|^2dx+\int_\Omega \nabla g_2\cdot\nabla(v_\epsilon-g_2)+(1+\lambda\Phi_\epsilon(u_\epsilon))\chi_\epsilon(v_\epsilon)(v_\epsilon-g_2)dx.
		\end{array}
		\right.\label{3.5}
	\end{eqnarray}
	We first show that $u_\epsilon$ has the desirous boundedness. Based on Young inequality, according the following a sequence of estimates
	\begin{eqnarray}
		\left.
		\begin{array}{ll}
			&\displaystyle -\lambda\int_\Omega \Phi_\epsilon(v_\epsilon) \chi_\epsilon(u_\epsilon)(u_\epsilon-g_1)dx\\[0.8em]
			&\displaystyle \leq \chi_{\{\lambda<0\}}\lambda\int_{\{-\epsilon\leq u_\epsilon\leq0\}\cap\Omega}\Phi_\epsilon(v_\epsilon)
			\chi_\epsilon(u_\epsilon)g_1dx-\chi_{\{\lambda>0\}}\lambda\int_{\{-\epsilon\leq u_\epsilon\leq0\}\cap\Omega}\Phi_\epsilon(v_\epsilon)
			\chi_\epsilon(u_\epsilon)u_\epsilon dx\\[0.8em]
			&\ \ \displaystyle +\chi_{\{\lambda>0\}}\lambda\int_{\{-\epsilon\leq u_\epsilon\leq0\}\cap\Omega}(|v_\epsilon|+1)|g_1|dx   +\chi_{\{\lambda>0\}}\lambda\int_{\{u_\epsilon>0\}\cap\Omega} (|v_\epsilon|+1)|g_1|dx\\[0.8em]
			&\ \ \displaystyle -\chi_{\{\lambda<0\}} \lambda\int_{\{u_\epsilon>0\}\cap\Omega}\Phi_\epsilon(v_\epsilon)\chi_\epsilon(u_\epsilon)(u_\epsilon-g_1) dx\\[0.8em]
			&\displaystyle \leq \chi_{\{\lambda<0\}}|\lambda|\int_{\Omega}(|v_\epsilon|+1)
			|g_1|dx+\epsilon\chi_{\{\lambda>0\}}\lambda\int_\Omega(|v_\epsilon-g_2|+|g_2|+1)
			dx\\[0.8em]
			&\ \ \displaystyle +\chi_{\{\lambda>0\}}\lambda\int_{\Omega}(|v_\epsilon-g_2|+|g_2|+1)|g_1|dx-\chi_{\{\lambda<0\}} \lambda\int_{\Omega}(|v_\epsilon|+1)(|u_\epsilon|+|g_1|)dx\\[0.8em]
			&\displaystyle \leq C(\lambda, \Omega, |g_1|_{L^2(\Omega)},|g_2|_{L^2(\Omega)})+ C(\lambda,\Omega,|g_1|_{L^2(\Omega)})|\nabla (v_\epsilon-g_2)|_{L^2(\Omega)}\\[0.8em]
			&\ \ \displaystyle +C(\lambda,\Omega,|g_2|_{L^2(\Omega)}) |\nabla (u_\epsilon-g_1)|_{L^2(\Omega)}-\chi_{\{\lambda<0\}} \frac{\lambda}{2\lambda_1}\int_{\Omega}|\nabla(v_\epsilon-g_2)|^2+|\nabla (u_\epsilon-g_1)|^2dx,
		\end{array}
		\right.\label{3.6}
	\end{eqnarray}
	using \eqref{3.5}, we have
	\begin{eqnarray*}
		\left.
		\begin{array}{ll}
			&\displaystyle \int_\Omega |\nabla(u_\epsilon-g_1)|^2dx\\[0.8em]
			&\displaystyle =-\int_\Omega \nabla g_1\cdot\nabla(u_\epsilon-g_1)-\int_\Omega(1+\lambda\Phi_\epsilon(v_\epsilon))
			\chi_\epsilon(u_\epsilon)(u_\epsilon-g_1)dx\\[0.8em]
			&\displaystyle \leq C(\|\nabla g_1\|_{L^2(\Omega)}) \|\nabla (u_\epsilon-g_1)\|_{L^2(\Omega)}+C(\Omega)\|\nabla (u_\epsilon-g_1)\|_{L^2(\Omega)}\\[0.8em]
			&\ \ \ \displaystyle +C(\lambda, \Omega, \|g_1\|_{L^2(\Omega)},\|g_2\|_{L^2(\Omega)})+ C(\lambda,\Omega,\|g_1\|_{L^2(\Omega)})\|\nabla (v_\epsilon-g_2)\|_{L^2(\Omega)}\\[0.8em]
			&\ \ \ \ \displaystyle +C(\lambda,\Omega,\|g_2\|_{L^2(\Omega)}) \|\nabla (u_\epsilon-g_1)\|_{L^2(\Omega)}-\chi_{\{\lambda<0\}} \frac{\lambda}{2\lambda_1}\int_{\Omega}|\nabla(v_\epsilon-g_2)|^2+|\nabla (u_\epsilon-g_1)|^2dx\\[0.8em]
			&\displaystyle \leq C(\lambda,\Omega,\|g_2\|_{L^2(\Omega)},\|\nabla g_1\|_{L^2(\Omega)}) \|\nabla (u_\epsilon-g_1)\|_{L^2(\Omega)}\\[0.8em]
			&\ \ \ \ \displaystyle +C(\lambda, \Omega, \|g_1\|_{L^2(\Omega)},\|g_2\|_{L^2(\Omega)})+ C(\lambda,\Omega,\|g_1\|_{L^2(\Omega)})\|\nabla (v_\epsilon-g_2)\|_{L^2(\Omega)}\\[0.8em]
			&\ \ \ \ \displaystyle -\chi_{\{\lambda<0\}} \frac{\lambda}{2\lambda_1}\int_{\Omega}|\nabla(v_\epsilon-g_2)|^2+|\nabla (u_\epsilon-g_1)|^2dx,
		\end{array}
		\right.
	\end{eqnarray*}
	which implies that
	\begin{eqnarray}
		\left.
		\begin{array}{ll}
			&\displaystyle \left(1+ \frac{\lambda\chi_{\{\lambda<0\}}}{2\lambda_1}\right)\int_\Omega |\nabla(u_\epsilon-g_1)|^2dx\\[0.8em]
			&\displaystyle \leq C(\lambda,\Omega,\|g_2\|_{L^2(\Omega)},\|\nabla g_1\|_{L^2(\Omega)}) \|\nabla (u_\epsilon-g_1)\|_{L^2(\Omega)}+C(\lambda, \Omega, \|g_1\|_{L^2(\Omega)},\|g_2\|_{L^2(\Omega)})\\[0.8em]
			&\ \ \ \displaystyle + C(\lambda,\Omega,\|g_1\|_{L^2(\Omega)})\|\nabla (v_\epsilon-g_2)\|_{L^2(\Omega)} -\chi_{\{\lambda<0\}} \frac{\lambda}{2\lambda_1}\int_{\Omega}|\nabla(v_\epsilon-g_2)|^2dx.
		\end{array}
		\right.\label{3.7}
	\end{eqnarray}
	Copying the proof above on $v_\epsilon$, we can obtain that
	\begin{eqnarray}
		\left.
		\begin{array}{ll}
			&\displaystyle \left(1+ \frac{\lambda\chi_{\{\lambda<0\}}}{2\lambda_1}\right)\int_\Omega |\nabla(v_\epsilon-g_2)|^2dx\\
			&\displaystyle \leq C(\lambda,\Omega,\|g_1\|_{L^2(\Omega)},\|\nabla g_2\|_{L^2(\Omega)}) \|\nabla (v_\epsilon-g_2)\|_{L^2(\Omega)}+C(\lambda, \Omega, \|g_1\|_{L^2(\Omega)},\|g_2\|_{L^2(\Omega)})\\
			&\ \ \ \displaystyle + C(\lambda,\Omega,\|g_2\|_{L^2(\Omega)})\|\nabla (u_\epsilon-g_1)\|_{L^2(\Omega)} -\chi_{\{\lambda<0\}} \frac{\lambda}{2\lambda_1}\int_{\Omega}|\nabla(u_\epsilon-g_1)|^2dx.
		\end{array}
		\right.\label{3.8}
	\end{eqnarray}
	It follows from \eqref{3.7} and \eqref{3.8} that
	\begin{eqnarray}
		\left.
		\begin{array}{ll}
			&\displaystyle \left(1+ \frac{\lambda\chi_{\{\lambda<0\}}}{2\lambda_1}\right)\int_\Omega |\nabla(u_\epsilon-g_1)|^2dx\\
			&\displaystyle \leq C(\lambda,\lambda_1,\Omega,\|g_2\|_{L^2(\Omega)},\|\nabla g_1\|_{L^2(\Omega)}) \|\nabla (u_\epsilon-g_1)\|_{L^2(\Omega)}+C(\lambda, \lambda_1, \Omega, \|g_1\|_{L^2(\Omega)},\|g_2\|_{L^2(\Omega)})\\
			&\ \ \ \displaystyle + C(\lambda,\lambda_1,\Omega,\|g_1\|_{L^2(\Omega)},\|\nabla g_2\|_{L^2(\Omega)})\|\nabla (v_\epsilon-g_2)\|_{L^2(\Omega)}\\
			&\ \ \ \displaystyle +\frac{\lambda^2\chi_{\{\lambda<0\}}}{(2\lambda_1+\lambda\chi_{\{\lambda<0\}})2\lambda_1}\int_{\Omega}
			|\nabla(u_\epsilon-g_1)|^2dx.
		\end{array}
		\right.\label{3.9}
	\end{eqnarray}
	We  obtain from \eqref{3.8} that
	\begin{eqnarray}
		\left.
		\begin{array}{ll}
			\displaystyle \|\nabla (v_\epsilon-g_2) \|_{L^2(\Omega)}\leq& C(\lambda,\lambda_1,\Omega,\|g_1\|_{L^2(\Omega)},\|g_2\|_{L^2(\Omega)},\|\nabla g_2\|_{L^2(\Omega)})\\[0.8em]
			&\displaystyle + C(\lambda,\lambda_1,\Omega,\|g_2\|_{L^2(\Omega)})\|\nabla (u_\epsilon-g_1)\|_{L^2(\Omega)}^{1/2}\\[0.8em]
			&\displaystyle +C(\lambda,\lambda_1,\Omega)\|\nabla (u_\epsilon-g_1)\|_{L^2(\Omega)}.
		\end{array}
		\right.\label{3.10}
	\end{eqnarray}
	Combining with \eqref{3.9} and \eqref{3.10}, we get that
	\begin{eqnarray*}
		\left.
		\begin{array}{ll}
			&\displaystyle \left(1+ \frac{\lambda\chi_{\{\lambda<0\}}}{2\lambda_1}-\frac{\lambda^2\chi_{\{\lambda<0\}}}{(2\lambda_1+\lambda\chi_{\{\lambda<0\}})2\lambda_1}\right)\int_\Omega |\nabla(u_\epsilon-g_1)|^2dx\\[0.8em]
			&\displaystyle \leq C(\lambda,\lambda_1,\Omega,\|g_2\|_{L^2(\Omega)},\|g_1\|_{L^2(\Omega)},\|\nabla g_1\|_{L^2(\Omega)},\|\nabla g_2\|_{L^2(\Omega)}) \|\nabla (u_\epsilon-g_1)\|_{L^2(\Omega)}\\[0.8em]
			&\ \ \ +C(\lambda,\lambda_1,\Omega,\|g_1\|_{L^2(\Omega)},\|g_2\|_{L^2(\Omega)},\|\nabla g_2\|_{L^2(\Omega)})\|\nabla (u_\epsilon-g_1)\|_{L^2(\Omega)}^{1/2}\\[0.8em]
			&\ \ \ +C(\lambda,\lambda_1,\Omega,\|g_1\|_{L^2(\Omega)},\|g_2\|_{L^2(\Omega)},\|\nabla g_2\|_{L^2(\Omega)}),
		\end{array}
		\right.
	\end{eqnarray*}
	Since $1+ \frac{\lambda\chi_{\{\lambda<0\}} }{2\lambda_1}-\frac{\lambda^2\chi_{\{\lambda<0\}}}{(2\lambda_1+\lambda\chi_{\{\lambda<0\}})2\lambda_1}>0$ from $\lambda>-\lambda_1$ and $\lambda_1>0$, we get the desired result,
	\begin{eqnarray*}
		\|\nabla (u_\epsilon-g_1)\|_{L^2(\Omega)}\leq  C(\lambda,\lambda_1,\Omega,\|g_1\|_{L^2(\Omega)},\|g_2\|_{L^2(\Omega)},\|\nabla g_1\|_{L^2(\Omega)},\|\nabla g_2\|_{L^2(\Omega)}).
	\end{eqnarray*}
	Similarly, one has
	\begin{eqnarray*}
		\|\nabla (v_\epsilon-g_2)\|_{L^2(\Omega)}\leq  C(\lambda,\lambda_1,\Omega,\|g_1\|_{L^2(\Omega)},\|g_2\|_{L^2(\Omega)},\|\nabla g_1\|_{L^2(\Omega)},\|\nabla g_2\|_{L^2(\Omega)}).
	\end{eqnarray*}
	Applying again the Sobolev embedding theorem, we know that $\{(u_\epsilon, v_\epsilon)\}_\epsilon$ is uniformly bounded in $\WW$. As a consequence, we obtain that there exists $(u,v)\in \WW$ such that over a subsequence $\epsilon=\epsilon_k\to0$,
	\begin{eqnarray*}
		u_\epsilon\to u, v_\epsilon\to v\ {\rm weakly\ in}\ \W,\\
		u_\epsilon\to u,\ v_\epsilon\to v\ {\rm strongly\ in}\ L^2(\Omega).
	\end{eqnarray*}
	Moreover, since $(u_\epsilon-g_1,v_\epsilon-g_2)\in W^{1,2}_0(\Omega)\times W^{1,2}_0(\Omega)$, then  $(u,v) \in \KW$. Next, recalling that $(u_\epsilon, v_\epsilon)$ is weak solutions of regularized problem \eqref{3.4}. By the standard bootstrap argument, we get that $(u_\epsilon, v_\epsilon)\in W^{2,p}_{\rm loc}(\Omega)\times W^{2,p}_{\rm loc}(\Omega), 1\leq p<+\infty$ and it follows \cite[Chapter 9]{2001GilbargTrudinger} or \cite[Theorem 1.1]{2012PetrosyanShahgholianUraltseva} that $\{(u_\epsilon, v_\epsilon)\}_\epsilon$ is uniformly bounded in $W^{2,p}_{\rm loc}(\Omega)\times W^{2,p}_{\rm loc}(\Omega)$ on $\epsilon$. It follows that there exists a subsequence such that
	\begin{eqnarray*}
		u_\epsilon \rightharpoonup u,\ v_\epsilon\rightharpoonup v\ {\rm weakly\ in}\ W^{2,p}_{\rm loc}(\Omega).
	\end{eqnarray*}
	Then $(u,v)\in W^{2,p}_{\rm loc}(\Omega)\times W^{2,p}_{\rm loc}(\Omega)$. For $\epsilon=\epsilon_k\to0$, since
	\begin{eqnarray*}
		\int_\Omega |\nabla u|^2dx\leq \liminf_{\epsilon\to0}\int_\Omega|\nabla u_\epsilon|^2dx, \int_\Omega |\nabla v|^2dx\leq \liminf_{\epsilon\to0}\int_\Omega|\nabla v_\epsilon|^2dx.
	\end{eqnarray*}
	By Lebesgue theorem, we can know that
	\begin{eqnarray*}
		\int_\Omega u^+dx=\lim_{\epsilon\to0}\int_\Omega \Phi_{\epsilon}(u_\epsilon)dx,\ \int_\Omega v^+dx=\lim_{\epsilon\to0}\int_\Omega \Phi_{\epsilon}(v_\epsilon)dx,\  \int_\Omega v^+u^+dx=\lim_{\epsilon\to0}\int_\Omega \Phi_{\epsilon}(u_\epsilon)\Phi_{\epsilon}(v_\epsilon)dx.
	\end{eqnarray*}
	It follows that
	\begin{eqnarray*}
		\tilde{J}(u,v)\leq\liminf_{\epsilon\to0}J_\epsilon(u_\epsilon, v_\epsilon)\leq \liminf_{\epsilon\to0}J_\epsilon(w,\tilde{w})=\tilde{J}(w,\tilde{w})
	\end{eqnarray*}
	for any $(w,\tilde{w})\in\KW$. Applying Lemma \ref{lem3.3}, we get that $(u,v)$ is the minimizer of $J$ over $\mathcal{K}$, which implies that $u,v\geq0$ a.e. in $\Omega$. Finally we verify that $(u,v)$ solves the problem \eqref{system} in the sense of distributions. Since $(u,v)\in W^{2,p}_{\rm loc}(\Omega)\times W^{2,p}_{\rm loc}(\Omega)$, we readily have that $\Delta u\in L^p_{\rm loc}(\Omega)$  and we have to verify that
	\begin{eqnarray}
		\left\{
		\begin{array}{ll}
			\Delta u=(1+\lambda v)\chi_{\{u>0\}},\ {\rm a.e.\ in}\ \Omega,\\
			\Delta v=(1+\lambda u)\chi_{\{v>0\}},\ {\rm a.e.\ in}\ \Omega.
		\end{array}
		\right.\label{6}
	\end{eqnarray}
	To this end, the fact that $(u_\epsilon,v_\epsilon)\in W^{2,p}_{\rm loc}(\Omega)\times W^{2,p}_{\rm loc}(\Omega)$ leads to the equation
	\begin{eqnarray*}
		\left\{
		\begin{array}{ll}
			\Delta u_\epsilon=(1+\lambda\Phi_\epsilon(v_\epsilon))\chi_\epsilon(u_\epsilon),\ {\rm a.e.\ in}\ \Omega,\\
			\Delta v_\epsilon=(1+\lambda\Phi_\epsilon(u_\epsilon))\chi_\epsilon(v_\epsilon),\ {\rm a.e.\ in}\ \Omega
		\end{array}
		\right.
	\end{eqnarray*}
	holds in the strong sense. The arbitrariness of $p$ implies $u_\epsilon\to u,\ v_\epsilon\to v$ in $C^{1,\alpha}_{\rm loc}(\Omega)$ hold as  $\epsilon\to0$. Then the locally uniformly convergence and $u,v\geq0$ a.e. in $\Omega$ imply that
	\begin{eqnarray}
		\left\{
		\begin{array}{ll}
			\Delta u=1+\lambda v,\ {\rm a.e.\ in\ the\ open\ set}\ \{u>0\},\\
			\Delta v=1+\lambda u,\ {\rm a.e.\ in\ the\ open\ set}\ \{v>0\}.
		\end{array}
		\right.\label{4}
	\end{eqnarray}
	Since $(u,v)\in W^{2,p}_{\rm loc}(\Omega)\times W^{2,p}_{\rm loc}(\Omega)$, then \
	\begin{eqnarray}
		\left\{
		\begin{array}{ll}
			\Delta u=0,\ {\rm a.e.\ on}\ \{u=0\},\\
			\Delta v=0,\ {\rm a.e.\ on}\ \{v=0\}.
		\end{array}
		\right.\label{5}
	\end{eqnarray}
	It follows from \eqref{4} and \eqref{5} that \eqref{6} holds. \hbx

	\section{Optimal regularity of solutions}\label{sec4}
	
	In this section, we prove the  optimal regularity of solutions. By the standard bootstrap argument, we know that $(u,v)\in C^{1,\alpha}_{\rm loc}(\Omega)\times C^{1,\beta}_{\rm loc}(\Omega)$ with any  $\alpha,\beta\in(0,1)$. Furthermore, we will give the proof of the regularity of solutions for system \eqref{system}, that is Theorem \ref{thmopimalregularity}.

{\bf Proof of Theorem \ref{thmopimalregularity}.} To obtain the desired result, we consider individual equation in the system separately.   Since $u\in W^{2,p}_\loc(\Omega)$ with $p>1$, we see that at any Lebesgue point $x^0$ of $D^2 u$, the function $u$ is twice differentiable. Then fix such a point $x^0\in K\Subset \Omega$, where $u$ is twice differentiable and define
	\begin{eqnarray*}
		w(x)=\partial_eu(x)
	\end{eqnarray*}
	for a unit vector $e$ being orthogonal to $\nabla u(x^0)$(if $\nabla u(x^0)=0$, take an arbitrary unit vector $e$). Our aim is to obtain a uniform estimate for $\partial_{x_je}u(x^0)=\partial_{x_j}w(x^0),\ j=1,2,\cdots, n$. By construction, $w(x^0)=0$ because $e$ is orthogonal to $\nabla u(x^0)$. Hence, the Taylor expansion at $x^0$ is
	\begin{eqnarray*}
		w(x)=\xi\cdot (x-x^0)+o(|x-x^0|),\ \xi=\nabla w(x^0).
	\end{eqnarray*}
	Now, if $\xi=0$ then $\partial_{x_j}w(x^0)=0$ for $j=1,2,\cdots, n$ and there nothing to estimate. If $\xi\not=0$, we consider the cone
	\begin{eqnarray*}
		\mathcal{C}:=\left\{x\in \Rn: \xi\cdot (x-x^0)>\frac{|\xi||x-x^0|}{2}\right\},
	\end{eqnarray*}
	which has the property that
	\begin{eqnarray*}
		\mathcal{C}\cap B_r(x^0)\subset \{w>0\},\ -\mathcal{C}\cap B_r(x^0)\subset \{w<0\}
	\end{eqnarray*}
	for sufficiently small $r>0$. Considering also the rescalings
	\begin{eqnarray*}
		w_r(x)=\frac{w(rx+x^0)}{r},\ x\in B_1.
	\end{eqnarray*}
	Note that $w_r(x)\to w_0(x):=\xi\cdot x$ uniformly in $B_1$ and $\nabla w_r\to \nabla w_0$ in $L^p(B_1), p>n$. Since
	\begin{eqnarray*}
		\int_{B_1}|\nabla w_r(x)-\xi|^pdx=\frac{1}{r^n}\int_{B_r(x^0)}|\nabla w(y)-\nabla w(x^0)|^pdx\to 0,
	\end{eqnarray*}
	where the right-hand side goes to zero as $r\to0$ because $x^0$ is a Lebesgue point for $\nabla w$. Then for constant $C(n)$, using Lemma \ref{thmACF}  we obtain
	\begin{eqnarray*}
		(C(n))^2|\xi|^4&=&\int_{\mathcal{C}\cap B_1}\frac{|\nabla w_0(x)|^2}{|x|^{n-2}}dx\int_{-\mathcal{C}\cap B_1}\frac{|\nabla w_0(x)|^2}{|x|^{n-2}}dx\\
		&=&\lim_{r\to0}\int_{\mathcal{C}\cap B_1}\frac{|\nabla w_r(x)|^2}{|x|^{n-2}}dx\int_{-\mathcal{C}\cap B_1}\frac{|\nabla w_r(x)|^2}{|x|^{n-2}}dx\\
		&=&\lim_{r\to0}\frac{1}{r^4}\int_{\mathcal{C}\cap B_r(x^0)}\frac{|\nabla w(x)|^2}{|x-x^0|^{n-2}}dx\int_{-\mathcal{C}\cap B_r(x^0)}\frac{|\nabla w(x)|^2}{|x-x^0|^{n-2}}dx\\
		&\leq& \lim_{r\to0} \Phi(r, w^+,w^-,x^0).
	\end{eqnarray*}
	Let $E_\delta:=\{w>0\}\cap K_\delta$, then, formally, for $x\in E_\delta$,
	\begin{eqnarray}\label{3.1-1}
		\Delta(w^+)=\partial_e\Delta u=\partial_e(1+\lambda v(x))=\lambda e\cdot\nabla v\geq-|\lambda|\|\nabla v\|_{L^\infty(K_{2\delta})}.
	\end{eqnarray}
	To justify this computation, observe that $\Delta(w^+)\geq-|\lambda|\|\nabla v\|_{L^\infty(K_{2\delta})}$ in $K_\delta$ is equivalent to the inequality
	\begin{eqnarray}
		-\int_{K_\delta}\nabla(w^+)\nabla \eta dx\geq-|\lambda|\|\nabla v\|_{L^\infty(K_{2\delta})}\int_{K_\delta}\eta dx\label{3.1}
	\end{eqnarray}
	for any nonnegative $\eta\in C_0^\infty(K_\delta)$. Suppose that $\supp\ \eta\subset\{w>\zeta\}\cap K_\delta$ with $\zeta>0$, then writing the equation
	\begin{eqnarray*}
		-\int_\Omega\nabla u\nabla \eta dx=\int_\Omega(1+\lambda v)\eta dx
	\end{eqnarray*}
	with $\eta=\eta(x)$ and $\eta=\eta(x-he)$, we obtain an equation for the incremental quotient
	\begin{eqnarray*}
		v_{(h)}(x)=\frac{u(x+he)-u(x)}{h}.
	\end{eqnarray*}
	Namely,
	\begin{eqnarray}
		-\int_\Omega\nabla v_{(h)}\nabla \eta dx=\frac{\lambda}{h} \int_\Omega [v(x+he)-v(x)]\eta dx\label{3.3}
	\end{eqnarray}
	for small $h>0$. Note that $u(x+he)>u(x)$ on $\supp\ \eta\subset\{w>\zeta\}\cap K_\delta$. Moreover, on $\{w>\zeta\}\cap K_\delta$,
	\begin{eqnarray*}
		\lambda v(x+he)-\lambda v(x)=\lambda \partial_ev(x)h\geq-|\lambda|\cdot\|\nabla v\|_{L^\infty(K_{2\delta})}
	\end{eqnarray*}
	for small $h$. Letting $h\to0$ and $\zeta\to0$, we arrive at
	\begin{eqnarray}
		-\int_\Omega \nabla w\nabla \eta dx\geq-|\lambda|\|\nabla v\|_{L^\infty(K_{2\delta})}\int_\Omega \eta dx\label{3.3}
	\end{eqnarray}
	for any arbitrary $\eta\geq0$ with $\supp\ \eta\subset\{w>0\}\cap K_\delta$. Now we prove the desired result \eqref{3.1}. Since \eqref{3.3} holds and $C_0^\infty(\{w>0\}\cap K_\delta)$ is dense in $W^{1,2}_0(\{w>0\}\cap K_\delta)$, we get for any $\eta\geq0$ and $\eta\in W^{1,2}_0(\{w>0\}\cap K_\delta)$ that
	\begin{eqnarray*}
		\int_\Omega \nabla w\nabla\eta dx\leq|\lambda|\|\nabla v\|_{L^\infty(K_{2\delta})}\int_\Omega  \eta dx,
	\end{eqnarray*}
	plugging $\eta=\psi_\epsilon(w)\phi$ with nonnegative $\phi\in C_0^\infty(K_\delta)$ and $\psi_\epsilon$ satisfying
	\begin{eqnarray*}
		0\leq\psi_\epsilon\leq1, \psi'_\epsilon\geq0, \psi_\epsilon(t)=0\ {\rm for}\ t\leq\frac{\epsilon}{2}, \psi_\epsilon(t)=1\ {\rm for}\ t\geq\epsilon.
	\end{eqnarray*}
	Then
	\begin{eqnarray*}
		\int_\Omega (|\nabla w|^2 \psi'_\epsilon(w)\phi+\psi_\epsilon(w)\nabla w \nabla \phi)dx\leq |\lambda|\|\nabla v\|_{L^\infty(K_{2\delta})}\int_\Omega \psi_\epsilon(w)\phi dx.
	\end{eqnarray*}
	Let $\epsilon\to0^+$, we can verify the correctness of \eqref{3.1}.
	
	Coping the proof above, then
	\begin{eqnarray}\label{3.4}
		\Delta(w^-)\geq-|\lambda| \|\nabla v\|_{L^\infty(K_{2\delta})}.
	\end{eqnarray}
	Consequently, based on the \cite[Theorem 2.12\'{}]{2012PetrosyanShahgholianUraltseva},
	\begin{eqnarray*}
		\Phi(r,w^+,w^-,x^0)\leq C(n)\left(\delta^2|\lambda|^2\|\nabla v\|_{L^\infty(E_{2\delta})}^2
		+\frac{\|w^+\|_{L^2(B_\delta(x^0))}^2+\|w^-\|_{L^2(B_\delta(x^0))}^2}{\delta^{n+2}}\right)^2
	\end{eqnarray*}
	with $0<r\leq\frac{\delta}{2}$. It follows from the definition of $w$ that
	\begin{eqnarray*}
		(C(n))^2|\xi|^4\leq \lim_{r\to0}\Phi(r,w^+,w^-,x^0)&\leq& C(n)\left(\delta^2|\lambda|^2\|\nabla v\|_{L^\infty(K_{2\delta})}^2+\frac{\|\nabla u\|_{L^\infty(K_\delta)}^2}{\delta^2}\right)^2\\
		&\leq& C(n,\delta,\lambda)\left(\|\nabla v\|_{L^\infty(K_{2\delta})}+\|\nabla u\|_{L^\infty(K_\delta)}\right)^4.
	\end{eqnarray*}
	Furthermore,
	\begin{eqnarray*}
		|\xi|\leq C(n,\delta,\lambda)\left(\|\nabla v\|_{L^\infty(K_{2\delta})}+\|\nabla u\|_{L^\infty(K_\delta)}\right).
	\end{eqnarray*}
	Recalling that $\xi=\nabla\partial_e u(x^0)$, we get that
	\begin{eqnarray*}
		|\nabla\partial_e u(x^0)|\leq C(n,\delta,\lambda)\left(\|\nabla v\|_{L^\infty(K_{2\delta})}+\|\nabla u\|_{L^\infty(K_\delta)}\right).
	\end{eqnarray*}
	This does not give the desired estimate on $|D^2 u|$ yet, since $e$ merely satisfies the condition $e\cdot\nabla u(x^0)=0$, unless $\nabla u(x^0)=0$. For the case of $\nabla u(x^0)\not=0$ we may choose the coordinate system so that $\nabla u(x^0)$ is parallel to $e_n$. We take subsequently $e=e_1,e_2,\cdots, e_n$ in the estimate above, we arrive at
	\begin{eqnarray*}
		|\partial_{x_ix_j}u(x^0)|\leq C(n,\delta,\lambda) \left(\|\nabla v\|_{L^\infty(K_{2\delta})}+\|\nabla u\|_{L^\infty(K_\delta)}\right), i=1,2,\cdots, n-1, j=1,2,\cdots, n.
	\end{eqnarray*}
	To obtain the desired result in the missing direction $e_n$,  the original equation implies that
	\begin{eqnarray*}
		|\partial_{x_nx_n}u(x^0)|&\leq& |\Delta u(x^0)|+|\partial_{x_1x_1}u(x^0)|+\cdots+|\partial_{x_{n-1}x_{n-1}}u(x^0)|\\
		&\leq& C(n,\delta,\lambda)\left(1+\|v\|_{L^\infty(K)}+\|\nabla v\|_{L^\infty(K_{2\delta})}+\|\nabla u\|_{L^\infty(K_\delta)}\right).
	\end{eqnarray*}
	Consequently, the conclusion for $u$ is obtained obviously and the proof of results for $v$ just repeat the process above.\hbx
	
		\section{Asymptotic behavior of weak solutions} \label{secbehavior}
	Now, we consider the asymptotic behavior of solutions of obstacle system \eqref{system} as $\lambda\to0$ and give the proof of Theorem \ref{thmasy}.
	
{\bf proof of Theorem \ref{thmasy}.}\ Let $(u_\lambda,v_\lambda)$ be the solution of system, then
	\begin{eqnarray*}
		&&\int_\Omega \nabla u_\lambda\nabla(u_\lambda-g_1)dx+\nabla v_\lambda\nabla(v_\lambda-g_2)dx\\
		&=&-\int_\Omega (1+\lambda v_\lambda)\chi_{\{u_\lambda>0\}}(u_\lambda-g_1)+(1+\lambda u_\lambda)\chi_{\{v_\lambda>0\}}(v_\lambda-g_2)dx.
	\end{eqnarray*}
	It follows from H\"{o}lder inequalities that
	\begin{eqnarray*}
		&&\int_\Omega |\nabla u_\lambda|^2 +|\nabla v_\lambda|^2 dx\\
		&=&\int_\Omega \nabla u_\lambda\nabla g_1dx+\nabla v_\lambda\nabla g_2dx-\int_\Omega (1+\lambda v)\chi_{\{u_\lambda>0\}}(u_\lambda-g_1)\\
		&&\ \ \ \ \ +(1+\lambda u_\lambda)\chi_{\{v_\lambda>0\}}(v_\lambda-g_2)dx\\
		&\leq& \|\nabla u_\lambda\|_{L^2(\Omega)}\|\nabla g_1\|_{L^2(\Omega)}
		+\|\nabla v_\lambda\|_{L^2(\Omega)}\|\nabla g_2\|_{L^2(\Omega)}\\
		&&\ \ \ +\int_\Omega [1+|\lambda|(|v_\lambda-g_2|+|g_2|)]|u_\lambda-g_1|+[1+|\lambda|(|u_\lambda-g_1|+|g_1|]|v_\lambda-g_2|dx\\
		&\leq& \|\nabla u_\lambda\|_{L^2(\Omega)}\|\nabla g_1\|_{L^2(\Omega)}
		+\|\nabla v_\lambda\|_{L^2(\Omega)}\|\nabla g_2\|_{L^2(\Omega)}\\
		&&\ \ \ +C\|\nabla u_\lambda-\nabla g_1\|_{L^2(\Omega)}+\|\nabla v_\lambda-\nabla g_2\|_{L^2(\Omega)}\\
		&&\ \ \ \ +\int_\Omega 2|\lambda||v_\lambda-g_2||u_\lambda-g_1|+|\lambda||g_2||u_\lambda-g_1|+|\lambda||g_1||v_\lambda-g_2|dx\\
		&\leq& \|\nabla u_\lambda\|_{L^2(\Omega)}\|\nabla g_1\|_{L^2(\Omega)}
		+\|\nabla v_\lambda\|_{L^2(\Omega)}\|\nabla g_2\|_{L^2(\Omega)}\\
		&&\ \ \ +C\|\nabla u_\lambda-\nabla g_1\|_{L^2(\Omega)}+\|\nabla v_\lambda-\nabla g_2\|_{L^2(\Omega)}\\
		&&\ \ \ \ +\lambda_1|\lambda|\|\nabla v_\lambda-\nabla g_2\|_{L^2(\Omega)}^2+\lambda_1|\lambda|\|\nabla u_\lambda-\nabla g_1\|_{L^2(\Omega)}^2\\
		&&\ \ \ \  +|\lambda|C\|g_2\|_{L^2(\Omega)}\|\nabla u_\lambda-\nabla g_1\|_{L^2(\Omega)}+|\lambda|C\|g_1\|_{L^2(\Omega)}\|\nabla v_\lambda-\nabla g_2\|_{L^2(\Omega)}.
	\end{eqnarray*}
	As $\lambda$ small, we obtain
	\begin{eqnarray*}
		\left(1-\lambda_1|\lambda|\right)\left(\int_\Omega |\nabla u_\lambda|^2 +|\nabla v_\lambda|^2 dx\right)^{\frac12}\leq C(g_1,g_2,\lambda_1, n,\Omega).
	\end{eqnarray*}
	And applying \eqref{33.1},
	\begin{eqnarray*}
		\int_\Omega |u_\lambda-g_1|^2\leq \lambda_1\int_\Omega |\nabla u_\lambda-\nabla g_1|^2\ {\rm and}\   \int_\Omega |v_\lambda-g_2|^2\leq \lambda_1\int_\Omega |\nabla v_\lambda-\nabla g_2|^2.
	\end{eqnarray*}
	Therefore, for $\lambda$ small,
	\begin{eqnarray*}
		\left(\int_\Omega |\nabla u_\lambda|^2+|u_\lambda|^2+|\nabla v_\lambda|^2+|v_\lambda|^2 dx\right)^2\leq C
	\end{eqnarray*}
	with $C$ being independent on  $u,v$. By the standard bootstrap argument, we get that $(u_{\lambda}, v_\lambda)\in W^{2,p}_{\rm loc}(\Omega)\times W^{2,p}_{\rm loc}(\Omega),1\leq p<+\infty$ and it follows \cite[Chapter 9]{2001GilbargTrudinger} or \cite[Theorem 1.1]{2012PetrosyanShahgholianUraltseva} that $\{(u_\lambda, v_\lambda)\}_\lambda$ is uniformly bounded in $W^{2,p}_{\rm loc}(\Omega)\times W^{2,p}_{\rm loc}(\Omega)$ on $\lambda$. And there exists  a sequence such that
	\begin{eqnarray*}
		u_{\lambda_n}\to u_0, v_{\lambda_n}\to v_0\ \ {\rm in}\ W^{2,p}_{\rm loc}(\Omega)
	\end{eqnarray*}
and
\begin{eqnarray*}
		u_{\lambda_n}\to u_0, v_{\lambda_n}\to v_0\ \ {\rm in}\ C^{1,\alpha}_{\loc}(\Omega)
	\end{eqnarray*}
	as $\lambda_n \to0$ with $\alpha\in(0,1)$. Reduplicating the proof of Proposition \ref{prolimit}(\textbf{iii}), it is not difficult to know that $(u_0,v_0)$ is a solution of problem \eqref{classicalobstacle}. \hbx

	\section{Regularity of the free boundary}\label{sec5}
	In this section, we main consider the properties, including regularity, of the free boundary  $\Gamma(u,v)=\Gamma(u)\cup\Gamma(v)$ of system \eqref{system}, where
	\begin{eqnarray*}
		\Gamma(u)=\Omega\cap \partial\Omega(u)=\Omega\cap \partial\Lambda(u), \Gamma(v)=\Omega\cap \partial\Omega(v)=\Omega\cap\partial\Lambda(v)
	\end{eqnarray*}
	with $\Omega(u)=\{u>0\},\Lambda(u)=\{u=0\},\Omega(v)=\{v>0\},\Lambda(v)=\{v=0\}$.
	
	\subsection{Preliminary analysis of the free boundary}\label{sec5.1}

\subsubsection*{Nondegeneracy}
	We start with an important nondegenerary property.	
	
	\begin{lem}\label{lemnondegeneracy}
		Let $(u,v)$ be the solution of system \eqref{system}, $x^0\in \Omega(u)\cup\Gamma(u)\cup \Omega(v)\cup \Gamma(v)$ and $B_r(x^0)\subset \Omega$. Then
		\begin{itemize}
			\item if $x^0\in \Omega(u)\cup\Gamma(u)$, there holds
			\begin{eqnarray*}
				\sup_{\partial B_r(x^0)}u\geq u(x^0)+\frac{r^2}{2mn}
			\end{eqnarray*}
			under the one of assumptions:

{\rm Case (i).}\  $\lambda\geq0$ and $r$ small enough;

{\rm Case (ii).}\ \eqref{gammaulambda} for $m>1$ large and there exists $r_0>0$ small enough such that $r<r_0$.

			\item if $x^0\in \Omega(v)\cup \Gamma(v)$, there holds
			\begin{eqnarray*}
				\sup_{\partial B_r(x^0)}v\geq u(x^0)+\frac{r^2}{2mn}
			\end{eqnarray*}
			under the  one of assumptions:

{\rm Case (1).}\  $\lambda\geq0$\ and $r$ small enough;

{\rm Case (2).}\  \eqref{gammavlambda} for $m>1$ large and there exists $r_0>0$ small enough such that $r<r_0$.
		\end{itemize}
	\end{lem}
	\pf (I)\ If $x^0\in \Omega(u)$, consider the auxiliary function
	\begin{eqnarray*}
		w(x)=u(x)-u(x^0)-\frac{|x-x^0|^2}{2nm}.
	\end{eqnarray*}
	Since $u$ is a solution of system \eqref{system}, we can get that  in $B_r(x^0)\cap \Omega(u)$,
	\begin{eqnarray*}
		\Delta w=\Delta u-\frac{1}{m}=\frac{m-1}{m}+\lambda v
	\end{eqnarray*}

	{ (i)}\ If $\lambda\geq0$, then $\Delta w\geq0$.

	{ (ii)}\ If $\lambda<0$,  for $v(x^0)<-\frac{1}{\lambda}$, let $m$ be large such that $v(x^0)<-\frac{m-1}{\lambda m}$,  there exists $r_0>0$ small enough depending on $\lambda,m$ such that  $\Delta w\geq0$ as $r<r_0$.

	The obvious conclusion $w(x^0)=0$ and the maximum principle imply that
	\begin{eqnarray*}
		\sup_{\partial(B_r(x^0)\cap\Omega(u))} w\geq0.
	\end{eqnarray*}
	Besides, $w(x)=-u(x^0)-\frac{|x-x^0|^2}{2nm}<0$ on $\partial\Omega(u)$. Therefore, we must have
	\begin{eqnarray*}
		\sup_{\partial B_r(x^0)\cap\Omega(u)}w\geq0,
	\end{eqnarray*}
	which shows directly  that
	\begin{eqnarray*}
		\sup_{\partial B_r(x^0)\cap\Omega(u)}u(x)\geq u(x^0)+\frac{r^2}{2nm}.
	\end{eqnarray*}
{(II)}\ For the case of $x^0\in \partial\Omega(u)$, we take a sequence of $x^j\in \Omega(u)$ such that $x^j\to x^0$ as $j\to\infty$.

	{ (i)}\ If $\lambda\geq0$, then $\Delta w\geq0$. Furthermore, $\displaystyle \sup_{\partial B_r(x^j)\cap\Omega(u)}u(x)\geq u(x^j)+\frac{r^2}{2mn}$.

	{(ii)}\ If $\lambda<0$, as $v(x^0)<-\frac{1}{\lambda}$. For $m$ large, since there exists a $r_0$ such that $\frac{m-1}{m}+\lambda v>0$ in $B_{r_0}(x^0)$. Then for $j$ large, we have that  $\frac{m-1}{m}+\lambda v\geq0$ on $B_{r}(x^j)\subset B_{r_0}(x^0)$ as $r<\frac{r_0}{2}$ and $m$ large. Thus $\Delta w\geq0$ for any $r<\frac{r_0}{2}$, which leads to that $\displaystyle \sup_{\partial B_r(x^j)\cap\Omega(u)}u(x)\geq u(x^j)+\frac{r^2}{2mn}$ for large $m$.
	
	Therefore, taking the limit $j\to\infty$ gives that
	\begin{eqnarray*}
		\sup_{\partial B_r(x^0)\cap\Omega(u)}u(x)\geq u(x^0)+\frac{r^2}{2mn}.
	\end{eqnarray*}
The results for $v$ can be proved by the analogous arguments.\hbx
	
	\begin{rem}\label{rem5.4}
		Under the assumptions of Lemma {\rm \ref{lemnondegeneracy}}, then $u$ or $v$ is subharmonic in $B_r(x^0)$, furthermore supremum over $\partial B_r(x^0)$ can be replaced by over $B_r(x^0)$.
	\end{rem}

We can establish the nondegeneracy of the gradient.
	\begin{cor}\label{cor1}
		Let $(u,v)$ be the solution of system \eqref{system}, $x^0\in \Omega(u)\cup\Gamma(u)\cup \Omega(v)\cup \Gamma(v)$ and $B_r(x^0)\subset \Omega$. Then
		\begin{itemize}
			\item if $x^0\in \Omega(u)\cup\Gamma(u)$, there holds
			\begin{eqnarray*}
				\sup_{B_r(x^0)}|\nabla u|\geq \frac{r}{mn}
			\end{eqnarray*}
			under the one of assumptions: (i) $\lambda\geq0$ for any $m>1$ and  $r$ small enough; (ii)\ \eqref{gammaulambda} for $m>1$ large and there exists $r_0>0$ small enough such that $r<r_0$.
			\item if $x^0\in \Omega(v)\cup \Gamma(v)$, there holds
			\begin{eqnarray*}
				\sup_{B_r(x^0)}|\nabla v|\geq  \frac{r}{mn}
			\end{eqnarray*}
			under the  one of assumptions: (i) $\lambda\geq0$ for any $m>1$ and  $r$ small enough; (ii)\ \eqref{gammavlambda} for $m>1$ large and there exists $r_0>0$ small enough such that $r<r_0$.
		\end{itemize}
	\end{cor}
	\pf The results are easy to obtain from  Lemma \ref{lemnondegeneracy}, Remark \ref{rem5.4} and computations
	\begin{eqnarray*}
		\frac{r^2}{2mn}\leq \sup_{B_r(x^0)}(u(x)-u(x^0))=\sup_{B_r(x^0)}\frac{1}{2}\int_0^1\nabla u(x^0+t(x-x^0))\cdot(x-x^0)dt\leq \frac{r}{2}\sup_{B_r(x^0)}|\nabla u|
	\end{eqnarray*}
	for $x^0\in\Gamma(u)\cup\Omega(u)$ and
	\begin{eqnarray*}
		\frac{r^2}{2mn}\leq \sup_{B_r(x^0)}(v(x)-v(x^0))=\sup_{B_r(x^0)}\frac{1}{2}\int_0^1\nabla v(x^0+t(x-x^0))\cdot(x-x^0)dt\leq \frac{r}{2}\sup_{B_r(x^0)}|\nabla v|
	\end{eqnarray*}
	for $x^0\in\Gamma(v)\cup\Omega(v)$. \hbx
	
\subsubsection*{Lebesgue and Hausdorff measures of the free boundary}

	Now we give a   porosity for the free boundary $\Gamma(u,v)$.
	\begin{lem}\label{lem1porosity}
		Let $(u,v)$ be a solution of system \eqref{system} in an open set $\Omega\subset\Rn$. Then
		\begin{itemize}
			\item $\Gamma(u)$ is locally porous as $\lambda\geq0$ or \eqref{gammaulambda} holds.
			\item $\Gamma(u)$ is locally porous as $\lambda\geq0$ or \eqref{gammavlambda} holds.
		\end{itemize}
	\end{lem}
	\pf Let $K\Subset \Omega, x^0\in\Gamma(u)$ and $B_r(x^0)\subset K$. We consider the two cases.
	
	(i)\ If $\lambda\geq0$. Using the nondegeneracy in Lemma \ref{lemnondegeneracy} of the solution with $m=2$, one can find $y\in \overline{B_{r/2}(x^0)}$ such that
	\begin{eqnarray*}
		|u(y)|\geq \frac{r^2}{8n}.
	\end{eqnarray*}
	Now, the fact $|D u|_{L^\infty(K)}<\infty$ implies that
	\begin{eqnarray*}
		\displaystyle \inf_{B_{\delta r}(y)}u\geq\left(\frac{1}{4n}r-|D u|_{L^\infty(K)}\delta\right)r
	\end{eqnarray*}
	if $\delta^2\leq \min\left\{\frac{1}{2},\frac{r}{8n|D u|_{L^\infty(K)}}\right\}$. This leads to
	\begin{eqnarray*}
		B_{\delta r}(y)\subset B_r(x^0)\cap\Omega(u)\subset B_r(x^0)\backslash\Gamma(u).
	\end{eqnarray*}
	Hence the porosity condition is satisfied for any ball centered at $\Gamma(u)$ and contained in $K$.
	
	(ii)\ If $\lambda<0$, $v(x^0)<-\frac{1}{\lambda}$ and $x^0\in\Gamma(u)$. By Lemma \ref{lemnondegeneracy}, we get that for $m$ large, there exists $r_0$ such that as $s<r_0$, there holds,
	\begin{eqnarray*}
		\sup_{\partial B_s(x^0)}u\geq \frac{s^2}{2mn}.
	\end{eqnarray*}
	For the ball $B_r(x^0)\subset K$, as $k$ large we have $r/k<r_0$. Thus, there exists $y\in \overline{B_{r/k}(x^0)}$ such that
	\begin{eqnarray*}
		|u(y)|\geq \frac{r^2}{4kmn}.
	\end{eqnarray*}
	Now, the fact $|D u|_{L^\infty(K)}<\infty$ implies that
	\begin{eqnarray*}
		\displaystyle \inf_{B_{\delta r}(y)}u\geq\left(\frac{1}{4kmn}r-|D u|_{L^\infty(K)}\delta\right)r
	\end{eqnarray*}
	if $\delta\leq \min\left\{\frac{1}{k},\frac{r}{8kmn|D u|_{L^\infty(K)}}\right\}$. This leads to
	\begin{eqnarray*}
		B_{\delta r}(y)\subset B_r(x^0)\cap\Omega(u)\subset B_r(x^0)\backslash\Gamma(u).
	\end{eqnarray*}
	Hence the porosity condition is satisfied for any ball centered at $\Gamma(u)$ and contained in $K$.

	Therefore, $\Gamma(u)$ is locally porous in $\Omega$. And similarly, the result that $\Gamma(v)$ is locally porous in $\Omega$ holds.\hbx
	
	Now we establish the Lebesgue measure of $\Gamma(u,v)$ based on the porosity.
	\begin{cor}\label{corGamma}\ Let $(u,v)$ be a solution of system, then
		\begin{itemize}
			\item $\Gamma(u)$  has a Lebesgue measure zero  as $\lambda\geq0$ or \eqref{gammaulambda} holds.
			\item $\Gamma(v)$  has a Lebesgue measure zero as $\lambda\geq0$ or \eqref{gammavlambda} holds.
		\end{itemize}
	\end{cor}
	\pf The statement follows immediately from the local porosity of $\Gamma(u), \Gamma(v)$ and the \cite[Proposition 3.9]{2012PetrosyanShahgholianUraltseva}.\hbx
	
	Now we give a lemma about boundedness of density from below for $\Omega(u)$ and $\Omega(v)$.
	\begin{lem}\label{lemdensity}
		Let $(u,v)$ be a solution of system in $\Omega$, then
		\begin{itemize}
			\item $$\frac{|B_r(x^0)\cap\Omega(u)|}{|B_r|}\geq \beta_1,$$
			if $ x^0\in\Gamma(u)$ and $\lambda\geq0$ or \eqref{gammaulambda} holds.
			\item $$\frac{|B_r(x^0)\cap\Omega(v)|}{|B_r|}\geq \beta_1,$$
			if $ x^0\in\Gamma(u)$ and  $\lambda\geq0$ or \eqref{gammavlambda} holds.
		\end{itemize}
		provided $B_r(x^0)\subset\Omega$, where $\beta_1,\beta_2$ are positive constants, which depend only on $|D u|_{L^\infty(\overline{B_r(x^0)})}$ and $|D v|_{L^\infty(\overline{B_r(x^0)})}$ respectively as well as $n$.
	\end{lem}
	\pf Let $x^0\in\Gamma(u)$, the proof of Lemma \ref{lem1porosity} implies that
	\begin{eqnarray*}
		\frac{|B_r(x^0)\cap\Omega(u)|}{|B_r|}\geq \delta^n.
	\end{eqnarray*}
	Similarly, there exists $\widetilde{\delta}>0$ such that
	\begin{eqnarray*}
		\frac{|B_r(x^0)\cap\Omega(v)|}{|B_r|}\geq \widetilde{\delta}^n.
	\end{eqnarray*}
	These complete the proof.
	
In fact, the porosity of the free boundary implies not only its Lebesgue measure is zero but  also that it actually has a Hausdorff dimension less than $n$. We will further use the results of Lemma \ref{lemdensity} and Besicovitch covering lemma to prove that the free boundary of system \eqref{system} has locally finite $n-1$ dimensional Hausdorff measure. The detailed  proof given in Appendix.
	
	\begin{pro}\label{lemHausdorff} If $(u,v)$ is a pair of solution for the system \eqref{system} in an open set $\Omega$. Then
		\begin{itemize}
			\item $\Gamma(u)$ is a set of finite $(n-1)$-dimensional Hausdorff measure locally in $\Omega$  as $\lambda\geq0$  or \eqref{gammaulambda} holds.
			\item $\Gamma(v)$ is a set of finite $(n-1)$-dimensional Hausdorff measure locally in $\Omega$ as $\lambda\geq0$  or \eqref{gammavlambda} holds.
		\end{itemize}
	\end{pro}

\subsubsection*{Rescalings and blowups}\label{sec5.2}
	
The further analysis of the free boundary is based on the study of so-called blowup.
	
	If $(u,v)\in C_{\rm loc}^{1,1}(\Omega)\times  C_{\rm loc}^{1,1}(\Omega)$  and $r>0$, taking  $B_R(x^0)\Subset \Omega$ with $x^0\in\Gamma(u,v)$, then the rescaling of $u$ or $v$ at $x^0$ for $x\in B_{R/r}$ defined by
	\begin{eqnarray*}
		\left\{
		\begin{array}{ll}
			u_r(x)=u_{x^0,r}(x)=\frac{u(x^0+r x)-u(x^0)}{r^2},\ {\rm if}\ x^0\in\Gamma(u),\\[0.5em]
			v_r(x)=v_{x^0,r}(x)=\frac{v(x^0+r x)-v(x^0)}{r^2}, \ {\rm if}\ x^0\in\Gamma(v)
		\end{array}
		\right.
	\end{eqnarray*}
	satisfy
	\begin{eqnarray}
		\left.
		\begin{array}{ll}
			|D^2u_r(x)|=|D^2 u(x^0+r x)|\leq \|D^2 u\|_{L^\infty(B_R(x^0))},\ {\rm if}\ x^0\in\Gamma(u),\\[0.5em]
			|D^2v_r(x)|=|D^2 v(x^0+r x)|\leq \|D^2 v\|_{L^\infty(B_R(x^0))},\ {\rm if}\ x^0\in\Gamma(v),
		\end{array}
		\right.\label{5.3.2-0}\\
		|\nabla u_r(x)|\leq \|D^2 u(x)\|_{L^\infty(B_R(x^0))}|x|,\ x\in B_{R/r},\ {\rm if}\ x^0\in\Gamma(u),\nonumber\\
		|\nabla v_r(x)|\leq \|D^2 u(x)\|_{L^\infty(B_R(x^0))}|x|,\ x\in B_{R/r},\ {\rm if}\ x^0\in\Gamma(v)\nonumber
	\end{eqnarray}
	and
	\begin{eqnarray*}
		|u_r(x)|\leq \frac{\|D^2 u(x)\|_{L^\infty(B_R(x^0))}}{2}|x|^2,\  x\in B_{R/r}, {\rm if}\ x^0\in \Gamma(u),\\
		|v_r(x)|\leq \frac{\|D^2 u(x)\|_{L^\infty(B_R(x^0))}}{2}|x|^2,\  x\in B_{R/r}, {\rm if}\  x^0\in \Gamma(v).
	\end{eqnarray*}
	Hence, the Arzel\`{a}-Ascoli theorem implies there exists a subsequence $r=r_i\to0$ such that
	\begin{eqnarray}\label{5.2}
		\left.
		\begin{array}{ll}
			u_{r_i}\to u_0\ {\rm in}\ C^{1,\alpha}_\loc(\Rn),\ 0<\alpha<1,\ {\rm if}\ x^0\in \Gamma(u),\\
			v_{r_i}\to v_0\ {\rm in}\ C^{1,\alpha}_\loc(\Rn),\ 0<\alpha<1,\ {\rm if}\ x^0\in \Gamma(v),
		\end{array}
		\right.
	\end{eqnarray}
	where $u_0,v_0\in C^{1,1}_\loc(\Rn)$. And it follows from \eqref{5.3.2-0} that
	\begin{eqnarray*}
		\left.
		\begin{array}{ll}
			\|D^2u_0\|_{L^\infty(\Rn)}\leq \|D^2 u\|_{L^\infty(B_R(x^0))},\ {\rm if}\ x^0\in\Gamma(u),\\[0.5em]
			\|D^2v_0\|_{L^\infty(\Rn)}\leq \|D^2 v\|_{L^\infty(B_R(x^0))},\ {\rm if}\ x^0\in\Gamma(v).
		\end{array}
		\right.
	\end{eqnarray*}
	
	\begin{pro}\label{prolimit}
		Let $(u,v)$ be a solution of system \eqref{system} and $(u_{r_i},v_{r_i})$ satisfy \eqref{5.2}. Assume either $\lambda\geq0$ or the facts \eqref{gammaulambda}  and \eqref{gammavlambda} hold,
		then following states hold.
		\begin{description}
			\item[(i)] If $x^0\in\Omega$, then we have the implications
			\begin{eqnarray*}
				\left\{
				\begin{array}{ll}
					u_0(x^0)>0\Rightarrow u_{r_i}>0\ {\rm on}\ B_\delta(x^0),\ i\geq i^0,\   x^0\in\Gamma(u),\\
					v_0(x^0)>0\Rightarrow v_{r_i}>0\ {\rm on}\ B_\delta(x^0),\ i\geq i^0,\ x^0\in\Gamma(v)
				\end{array}
				\right.
			\end{eqnarray*}
			for some $\delta>0$ and sufficiently large $i^0$.
			\item[(ii)] If $B_\delta(x^0)\subset \Omega$, then we have
			\begin{eqnarray*}
				|u_0|=0\ {\rm in}\ B_\delta(x^0)\Rightarrow |u_{r_i}|=0\ {\rm on}\ B_{\delta/2}(x^0),\ i\geq i^0,\ {\rm if}\ x^0\in\Gamma(u),\\
				|v_0|=0\ {\rm in}\ B_\delta(x^0)\Rightarrow |v_{r_i}|=0\ {\rm on}\ B_{\delta/2}(x^0),\ i\geq i^0,\ {\rm if}\ x^0\in\Gamma(v)
			\end{eqnarray*}
			for sufficiently large $i^0$.
			
			\item[(iii)] $(u_0,v_0)$ satisfy the system
			\begin{eqnarray*}
				\left\{
				\begin{array}{ll}
					\Delta u_0=(1+\lambda v(x^0))\chi_{\{u_0>0\}},\ {\rm if}\ x^0\in \Gamma(u),\\
					\Delta v_0=(1+\lambda u(x^0))\chi_{\{v_0>0\}},\ {\rm if}\  x^0\in \Gamma(v).
				\end{array}
				\right.
			\end{eqnarray*}
			\item[(iv)] If $x^{i_k}\in\Gamma(u_{i_k})$ and $x^{i_k}\to x^0\in \Omega$ as $i_k\to \infty$, then $x^0\in\Gamma(u_0)$. And if $x^{i_k}\in\Gamma(v_{i_k})$ and $x^{i_k}\to x^0\in \Omega$ as $i_k\to \infty$, then $x^0\in\Gamma(v_0)$.
			\item[(v)]
			\begin{eqnarray*}
				\left\{
				\begin{array}{ll}
					u_{r_{i_k}}\to u_0,\ {\rm if}\ x^0\in \Gamma(u),\\
					u_{r_{i_k}}\to v_0,\ {\rm if}\  x^0\in \Gamma(v)
				\end{array}
				\right.
			\end{eqnarray*}
			strongly in $W^{2,p}_\loc(\Omega)$ for any $1<p<\infty$.
			\item[(vi)] $\chi_{\{u_{r_{i_k}}>0\}}(x)\to \chi_{\{u_0>0\}}(x)$ a.e. in $\Omega$ if $x^0\in \Gamma(u)$ and $\chi_{\{v_{r_{i_k}}>0\}}(x)\to \chi_{\{v_0>0\}}(x)$ a.e. in $\Omega$ if $x^0\in \Gamma(v)$.
			
		\end{description}
		
	\end{pro}
	
	\pf See Appendix.
	
Based on Proposition \ref{prolimit}, the explicit expression about blowups of solutions will be derived.
	\begin{lem}\label{thmstrongblowup}
		Let $(u,v)\in C^{1,1}_\loc(\Omega)\times C^{1,1}_\loc(\Omega)$ be solution of system \eqref{system}.
		
		If $x^0\in \Omega$ satisfies $u(x^0)=0$, either $\lambda\geq0$ or \eqref{gammaulambda} holds.
		Then any strong limit $\displaystyle u_0=\lim_{i\to\infty}u_{x^0,r_i}$ in $W^{2,p}(\Rn), p>n$ over a sequence $r_i\to 0^+$ of rescaling $u_{x^0,r}(x)$ is either a monotone function of one independent variable or a homogenous quadratic polynomial.
		
		And if $x^0\in \Omega$ satisfies $v(x^0)=0$, either $\lambda\geq0$ or \eqref{gammavlambda} holds.
		Then conclusions replacing $u$ by $v$ above hold.
	\end{lem}
	\pf  For a unit vector $e$ and $r>0$, we define
	\begin{eqnarray*}
		\Phi(r,(\partial_e u)^+,(\partial_e u)^-,x^0)=\frac{1}{r^4}\int_{B_r(x^0)}\frac{|\nabla(\partial_e u)^+|^2)dx}{|x-x^0|^{n-2}}\int_{B_r(x^0)}\frac{|\nabla(\partial_e u)^-|^2)dx}{|x-x^0|^{n-2}}.
	\end{eqnarray*}
	The fact $u_{x^0,r_i}\to u_0$ strongly in $W^{2,p}$ implies that
	\begin{eqnarray*}
		\displaystyle \Phi(r,(\partial_e u_0)^+,(\partial_e u_0)^-,0)=\lim_{i\to\infty}\Phi(r,(\partial_e u_{x^0,r_i})^+,(\partial_e u_{x^0,r_i})^-,0).
	\end{eqnarray*}
	On the other hand, since
	\begin{eqnarray*}
		\displaystyle \Phi(rr_i,(\partial_e u)^+,(\partial_e u)^-,x^0)=\Phi(r,(\partial_e u_{x^0,r_i})^+,(\partial_e u_{x^0,r_i})^-,0).
	\end{eqnarray*}
	By the proof of Theorem \ref{thmopimalregularity}, we obtain that $\Delta(\partial_e u)^{\pm}\geq-|\lambda||\nabla v|_{L^\infty(B_{r+\delta}(x^0))}$ in $B_r(x^0)$ with $\delta=\frac{1}{2}\dist(B_r(x^0),\partial\Omega)$. Besides,
	\begin{eqnarray*}
		|(\partial_e u)^{\pm}|=\frac{1}{2}||\partial_e u|\pm \partial_e u|\leq |\nabla u(x)|\leq \|u\|_{C^{1,1}(\overline{B_r(x^0)})}|x-x^0|
	\end{eqnarray*}
	in $B_r(x^0)$. Applying \cite[Theorem 2.13]{2012PetrosyanShahgholianUraltseva}, we know that the limit $\Phi(0^+,(\partial_e u)^+,(\partial_e u)^-,x^0)$ exists. As a consequence, we obtain that
	\begin{eqnarray*}
		\displaystyle \Phi(r,(\partial_e u_0)^+,(\partial_e u_0)^-,0)=\lim_{i\to\infty}\Phi(r,(\partial_e u_{x^0,r_i})^+,(\partial_e u_{x^0,r_i})^-,0)=\Phi(0^+,(\partial_e u)^+,(\partial_e u)^-,x^0)
	\end{eqnarray*}
	is constant for $r>0$. It follows from  $\Delta u_0=(1+\lambda v(x^0))\chi_{\{u_0>0\}}$ in $\Rn$ and the proof   of Theorem \ref{thmopimalregularity} that $(\partial_e u_0)^{\pm}$ is subharmonic, since $1+\lambda v(x^0)$ is constant. Now from the case equality in the ACF monotonicity formula, either one or the other of the following statements holds:
	\begin{description}
		\item[(i)] One of the functions $(\partial_e u_0)^+$ and $(\partial_e u_0)^-$ vanishes identically in $\Rn$;
		\item[(ii)] $(\partial_e u_0)^\pm=k_{\pm}(x\cdot\omega)^\pm, x\in \Rn$, for some constants $k_\pm=k_\pm(e)>0$ and a unit vector $\omega=\omega(e)$.
	\end{description}
	To finish the proof we need to consider two possibilities. \\
	1) Suppose $\{u_0=0\}$ has a positive Lebesgue measure. Then the case (ii) above cannot hold for any direction $e$. Consequently, (i) holds for any direction $e$, which is equivalent to $u_0$ being a monotone function of one variable, see \cite[Exercise 3.5]{2012PetrosyanShahgholianUraltseva}.\\
	2) Suppose now that $|\{u_0=0\}|=0$. Then $u_0$ satisfies $\Delta u_0=1+\lambda v(x^0)$ a.e. in $\Rn$. Since $|\{u_0=0\}|=0$, we obtain that $\Delta u_0=1+\lambda v(x^0)$  in $\Rn$ combining with the regularity theory of elliptic PDE. Thus $\Delta D^2u_0=0$, and combining with the global boundedness of $D^2u_0$ in $\Rn$ from \eqref{5.3.2-0} and Liouvile theorem, we know that $u_0$ is quadratic polynomial. Finally, since $u_0(0)=0$, $|\nabla u_0(0)|=0$ and $u_0\geq0$, we have $u_0=\frac{1+\lambda v(x^0)}{2}x^T A x$ for some $A\geq0$ and it follows from $\Delta u_0=1+\lambda v(x^0)$ that ${\rm tr} A=1$. The proof of results about $v$ just copy the process above. \hbx
	
	Now, we prove that all blowup limits at a point $x^0$ satisfying $u(x^0)=0$ possess the same description, that is the unique type of blowups.
	
	\begin{lem}\label{thmunique} Let $u$ and $x^0$ be as in Lemma {\rm\ref{thmstrongblowup}} and  $\lambda\geq0$ or  \eqref{gammaulambda} holds. Then either all possible strong blowups of $u$ at $x^0$ are monotone function of one variable or all of them are homogeneous polynomials. The results  also hold for $v$.
	\end{lem}
	\pf If there exists a blowup $u_0$ which is a monotone function of one variable, then
	\begin{eqnarray*}
		\Phi(0^+,(\partial_e u)^+,(\partial_e u)^-,x^0)= \Phi(r,(\partial_e u_0)^+,(\partial_e u_0)^-,0)=0
	\end{eqnarray*}
	for any direction $e$. There exists a unit $e_0$ and monotone function $\phi$ such that $u_0=\phi(e\cdot e_0)$. Since
	\begin{eqnarray*}
		\partial_e u_0=\phi'(x\cdot e_0)e_0\cdot e,
	\end{eqnarray*}
	then $(\partial_e u_0)^+=0$ or $(\partial_e u_0)^-=0$. This, in turn, implies that for any other blowup $u_0'$,
 \begin{eqnarray*}
\Phi(0^+,(\partial_e u)^+,(\partial_e u)^-,x^0)= \Phi(r,(\partial_e u_0')^+,(\partial_e u_0')^-,0)=0
 \end{eqnarray*}
for any direction $e$ and $r>0$. The latter is possible only if $u_0'$ is monotone function of one variable. Coping the proof of $u$, we can confirm the results of $v$. \hbx
	
	Now we prove the classification of blowups.
	\begin{pro}\label{thmclassification}
		Let $(u_0,v_0)$ be a blowup of a solution of system \eqref{system} with $x^0\in \Gamma(u)$ or $x^0\in\Gamma(v)$. As $\lambda\geq0$  or \eqref{gammaulambda} holds. Then $u_0$ has one of the following forms
		\begin{itemize}
			\item Polynomial solutions $u_0(x)=\frac{1+\lambda v(x^0)}{2}(x\cdot Ax), x\in \Rn$, where $A$ is an $n\times n$ symmetric matrix with ${\rm tr}A=1$;
			\item Half-space solutions $u_0(x)=\frac{1+\lambda v(x^0)}{2}[(x\cdot e)^+]^2, x\in \Rn$, where $e$ is a unit vector.
		\end{itemize}
		As $\lambda\geq0$ or \eqref{gammavlambda} holds. Then $v_0$ has one of the following forms
		\begin{itemize}
			\item Polynomial solutions $v_0(x)=\frac{1+\lambda u(x^0)}{2}(x\cdot Ax), x\in \Rn$, where $A$ is an $n\times n$ symmetric matrix with ${\rm tr}A=1$;
			\item Half-space solutions $v_0(x)=\frac{1+\lambda u(x^0)}{2}[(x\cdot e)^+]^2, x\in \Rn$, where $e$ is a unit vector.
		\end{itemize}
	\end{pro}
	\pf By Lemmas \ref{thmstrongblowup} and \ref{thmunique}, we get that $u_0$ is a homogeneous polynomial or monotone function of one variable. If $u_0$ is a homogeneous polynomial, since $u_0(0)=0$, one has $u_0(x)\geq0$ in $\Rn$ and $u_0(x)=\frac{1+\lambda v(x^0)}{2}x^TAx$ with $A$ being an $n\times n$ symmetric matrix and ${\rm tr}A=1$. If $u_0$ is monotone function, since $u_0(0)=0=|\nabla u_0(0)|=0$, we obtain a monotone function $\phi$ such that $u_0(x)=\phi(x\cdot e)$ and $\phi(0)=\phi'(0)=0$. As $x\cdot e\geq0$, then $\phi(x\cdot e)\geq0$ and as $x\cdot e\leq0$ then $\phi(x\cdot e)=0$ since $u_0(x)\geq0$. According to $\Delta \phi(x\cdot e)=\phi''(t)=1+\lambda v(x^0)$ for $x\cdot e\geq0$, we get that $\phi'(t)=\int_0^t \phi''(s)ds=(1+\lambda v(x^0))t$ for $t\geq0$. And $\phi(t)=\int_0^t\phi'(s)ds=\frac{1+\lambda v(x^0)}{2}t^2$, furthermore, $\phi(s)=\frac{1+\lambda v(x^0)}{2}[(x\cdot e)^+]^2$. The proof of results for $v$ is similar.\hbx

Taking inspiration from \cite{2012PetrosyanShahgholianUraltseva}, as an application of the Weiss-type monotonicity formula in Appendix, we will provide a classification of regular and singular points from an energy perspective.
	\begin{pro}\label{proenegy}
		Let $(u,v)$ be a solution of system \eqref{system}. Then the balanced energy are upper semicontinuous functions of $x^0\in\Gamma(u)$ and $x^0\in\Gamma(v)$ which take only two values.  More precisely there exists a dimensional constant $\alpha_n>0$ such that
\begin{equation*}
			\omega(x^0,u_0):=W(r,u_0,v(x^0),0)\in\left\{\!\left(\!\frac{1\!+\!\lambda v(x^0)}{2}\!\right)^2\!\alpha_n,\!\left(\!\frac{1+\lambda v(x^0)}{2}\!\right)^2\!\frac{\alpha_n}{2}\right\}\ {\rm for\ any}\ x^0\in\Gamma(u), \!
\end{equation*}
\begin{equation*}
			\omega(x^0,v_0):=W(r,v_0,u(x^0),0)\in\!\left\{\!\left(\!\frac{1\!+\!\lambda u(x^0)}{2}\!\right)^2\!\alpha_n,\left(\frac{1+\lambda u(x^0)}{2}\!\right)^2\!\frac{\alpha_n}{2}\right\}\ {\rm for\ any}\ x^0\in\Gamma(v)\ \
\end{equation*}
		with
		\begin{eqnarray*}
			\alpha_n=\frac{\mathcal{H}^{n-1}(\partial B_1)}{n(n+2)}.
		\end{eqnarray*}
	\end{pro}
	\pf For $x^0\in\Gamma(u)$, the upper semicontinuity follows from that
	\begin{eqnarray*}
		W(r,u,v,\cdot)=:\omega_r(\cdot,u_0)\searrow \omega(\cdot,u_0)\ {\rm as}\ r\searrow0,
	\end{eqnarray*}
	and the function $\omega_r(\cdot,u_0)$ is continuous for any $r>0$. The results for $x^0\in\Gamma(v)$ are similar.

	\begin{lem}\label{lem2}
		Let $(u_0,v_0)$ be a homogeneous global solution of system \eqref{system}. Then there exists a constant $\alpha_n>0$ satisfying
		\begin{eqnarray*}
			W(r,u_0,v(x^0),0)=\left\{
			\begin{array}{ll}
				\left(\frac{1+\lambda v(x^0)}{2}\right)^2\frac{\alpha_n}{2},\ &{\rm if}\ u_0\ {\rm is\ a\ half\ space\ solution},\\[1em]
				\left(\frac{1+\lambda v(x^0)}{2}\right)^2\alpha_n,\ &{\rm if}\ u_0\ {\rm is\ a\ polynomial\ solution},
			\end{array}
			\right.
		\end{eqnarray*}
		\begin{eqnarray*}
			W(r,v_0,u(x^0),0)=\left\{
			\begin{array}{ll}
				\left(\frac{1+\lambda u(x^0)}{2}\right)^2\frac{\alpha_n}{2},\ &{\rm if}\ v_0\ {\rm is\ a\ half\ space\  solution},\\[1em]
				\left(\frac{1+\lambda u(x^0)}{2}\right)^2\alpha_n,\ &{\rm if}\ v_0\ {\rm is\ a\ polynomial\ solution}.
			\end{array}
			\right.
		\end{eqnarray*}
	\end{lem}
	\pf Since $\Delta u_0=(1+\lambda v(x^0))$ in $\Omega(u_0)$, we obtain that
	\begin{eqnarray*}
		&&W(r,u_0,v(x^0),0)\\
&=&W(1,u_0,v(x^0),0)\\
		&=&\int_{B_1}\frac{1}{2}|\nabla u_0|^2+u_0+\lambda v(x^0)u_0 dx-\int_{\partial B_1}u_0^2d \mathcal{H}^{n-1}\\
		&=&\int_{B_1}\frac{-1}{2}\Delta u_0 u_0+u_0+\lambda v(x^0)u_0 dx+\int_{\partial B_1}\frac{1}{2}\frac{\partial u_0}{\partial\nu}d \mathcal{H}^{n-1}-\int_{\partial B_1}u_0^2d \mathcal{H}^{n-1}\\
		&=&\int_{B_1}\frac{1+\lambda v(x^0)}{2} u_0dx+\frac{1}{2}\int_{\partial B_1}(x\cdot\nabla u_0-2u_0)d \mathcal{H}^{n-1}\\
		&=&\int_{B_1}\frac{1+\lambda v(x^0)}{2} u_0dx.
	\end{eqnarray*}
	The fact is obvious,
	\begin{eqnarray*}
		\alpha_n=\int_{B_1}x_n^2dx=\frac{1}{n}\int_{B_1}|x|^2dx=\frac{\mathcal{H}^{n-1}(\partial B_1)}{n(n+2)}.
	\end{eqnarray*}
	Then for polynomial solution $u_0(x)=\frac{1+\lambda v(x^0)}{2}(x\cdot Ax)$, we can compute that
	\begin{eqnarray*}
		W(r,\frac{1+\lambda v(x^0)}{2}(x\cdot Ax),v(x^0),0)&=&\left(\frac{1+\lambda v(x^0)}{2}\right)^2\int_{B_1}x\cdot Ax dx\\
		&=&\left(\frac{1+\lambda v(x^0)}{2}\right)^2\alpha_n {\rm tr}A\\
		&=&\left(\frac{1+\lambda v(x^0)}{2}\right)^2\alpha_n.
	\end{eqnarray*}
	If $u_0(x)=\frac{1+\lambda v(x^0)}{2}[(x\cdot e)^+]^2$, we have
	\begin{eqnarray*}
		W(r,\frac{1+\lambda v(x^0)}{2}[(x\cdot e)^+]^2,v(x^0),0)&=&\left(\frac{1+\lambda v(x^0)}{2}\right)^2 \int_{B_1}[(x\cdot e)^+]^2dx \\ &=&\left(\frac{1+\lambda v(x^0)}{2}\right)^2\frac{\alpha_n}{2}.
	\end{eqnarray*}
After similar calculations, we can verify the conclusions of $v$. \hbx
	\begin{pro}\label{proclass}
		Let $(u,v)$ be a solution of system. Then\\
		if $x^0\in\Gamma(u)$,
		\begin{center}
			$x^0$ is regular $\Leftrightarrow$ $\omega(x^0,u_0)=\left(\frac{1+\lambda v(x^0)}{2}\right)^2\frac{\alpha_n}{2}$,\\
			$x^0$ is singular $\Leftrightarrow$ $\omega(x^0,u_0)=\left(\frac{1+\lambda v(x^0)}{2}\right)^2\alpha_n$.
		\end{center}
		If $x^0\in\Gamma(v)$,
		\begin{center}
			$x^0$ is regular $\Leftrightarrow$ $\omega(x^0,v_0)=\left(\frac{1+\lambda u(x^0)}{2}\right)^2\frac{\alpha_n}{2}$,\\
			$x^0$ is singular $\Leftrightarrow$ $\omega(x^0,v_0)=\left(\frac{1+\lambda u(x^0)}{2}\right)^2\alpha_n$,
		\end{center}
	\end{pro}

	\subsection{The regularity of free boundary near regular points}\label{subsec5.4}
	\begin{lem}\label{lem3}
		Let $(u,v)$ be the solution of system \eqref{system}.
		
		\begin{itemize}
			\item If  $0\in\Gamma(u)$. As $\lambda\geq0$ or \eqref{gammaulambda} holds.  For a half space solution $u_0(x)=\frac{1+\lambda v(0)}{2}(x_n^+)^2$ satisfying
			\begin{eqnarray*}
				\|u_{r}-u_0\|_{L^\infty(B_1)}\leq \epsilon
			\end{eqnarray*}
			with small $\epsilon$. Then $u_r>0$ in $\{x_n>\sqrt{2}\epsilon\}\cap B_1$ and $u_r=0$ in $\{x_n\leq -2\sqrt{mn\epsilon}\}\cap B_{1/2}$. In particular, $\Gamma(u_r)\cap B_{1/2}\subset \{|x_n|\leq 2\sqrt{mn\epsilon}\}$, where $m$ as in  Lemma {\rm\ref{lemnondegeneracy}}.
			\item If  $0\in\Gamma(v)$. As $\lambda\geq0$ or \eqref{gammavlambda} holds.  For a half space solution $v_0(x)=\frac{1+\lambda u(0)}{2}(x_n^+)^2$ satisfying
			\begin{eqnarray*}
				\|v_{r}-v_0\|_{L^\infty(B_1)}\leq \epsilon
			\end{eqnarray*}
			with small $\epsilon$. Then $v_r>0$ in $\{x_n>\sqrt{2}\epsilon\}\cap B_1$ and $v_r=0$ in $\{x_n\leq -2\sqrt{mn\epsilon}\}\cap B_{1/2}$. In particular, $\Gamma(v_r)\cap B_{1/2}\subset \{|x_n|\leq 2\sqrt{mn\epsilon}\}$, where $m$ as in  Lemma {\rm\ref{lemnondegeneracy}}.
		\end{itemize}
	\end{lem}
	\pf The first statement is evident $u_r\geq u_0-\epsilon=\frac{1+\lambda v(0)}{2}(x_n^+)^{\frac12}-\epsilon>(1+\lambda v(0))\epsilon^{\frac12}-\epsilon\geq0$ with the assumptions on $\lambda$ when $\epsilon$ small. So we will prove only the second one. If $\Omega(u_r)\cap B_{1/2}\cap\{x_n<0\}$ is empty, then the second one holds obviously. Otherwise, take any $x^1\in\Omega(u_r)\cap B_{1/2}\cap\{x_n<0\}$ and let $s=-x^1_n>0$. Then by   Lemma \ref{lemnondegeneracy}, there exists $m>0$,
	\begin{eqnarray*}
		\sup_{B_s(x^1)}u_r\geq u_r(x^1)+\frac{1}{2mn}s^2.
	\end{eqnarray*}
	Now, since $B_s(x^1)\subset \{x_n<0\}, u_0=0$ in $B_s(x^1)$ and consequently $|u_r|\leq\epsilon$ there. Therefore $\frac{s^2}{2mn}\leq 2\epsilon\Leftrightarrow s\leq2\sqrt{mn\epsilon}$, which leads to that $\Omega(u_r)\cap B_{1/2}\subset \{x_n>-2\sqrt{mn\epsilon}\}$. The case of $v$ is similar.\hbx
	
	Now we establish the flatness property of free boundary $\Gamma(u,v)$ at regular points.
	\begin{pro}\label{pro} Let $(u,v)$ be a solution of system \eqref{system} such that $0$ is a regular point. Then
		
		\begin{description}
			\item[(i)] If  $0\in\Gamma(u)$, either  $\lambda\geq0$ or \eqref{gammaulambda} holds. For any $\sigma>0$, there exists $r^7_\sigma=r_\sigma(u)>0$ such that $\Gamma(u)$ is $\sigma$-flat in $B_r$ for any $r<r^7_\sigma$ in the sense that $\Gamma(u)\cap B_r\subset \{|x\cdot e|<\sigma r\}$ for some direction $e=e_{r,\sigma}$.
			\item[(ii)] If  $0\in\Gamma(v)$, either   $\lambda\geq0$ or \eqref{gammavlambda} holds.  For any $\sigma>0$, there exists $r^8_\sigma=r_\sigma(u)>0$ such that $\Gamma(v)$ is $\sigma$-flat in $B_r$ for any $r<r^8_\sigma$ in the sense that $\Gamma(v)\cap B_r\subset \{|x\cdot e|<\sigma r\}$ for some direction $e=e_{r,\sigma}$.
		\end{description}
		
	\end{pro}
	\pf Let $0\in\Gamma(u)$. Suppose the statement fails for some $\sigma>0$. Then there exists a sequence $r_i\to 0$ such that $\Gamma(u)$ is not $\sigma$-flat in $B_{r_i/2}$. Without loss of generality we may assume that $u_{r_i}\to u_0$ in $C_\loc^{1,\alpha}(\Rn)$, where $u_0(x)=\frac{1+\lambda v(0)}{2}(x_n^+)^2$. Then for $i\geq i_\epsilon$, we have $\|u_{r_i}-u_0\|_{L^\infty(B_1)}\leq \epsilon$ and therefore by Lemma \ref{lem3},
	\begin{eqnarray*}
		\Gamma(u_{r_j})\cap B_{1/2}\subset \{|x_n|\leq 2\sqrt{mn\epsilon}\},
	\end{eqnarray*}
	which will cause $\Gamma(u_{r_j})$ to be $\sigma$-flat in $B_{1/2}$ if $2\sqrt{m n\epsilon}<\frac{\sigma}{2}$, a contradiction. The results of  $v$ can  be proved using the analysis above.\hbx
	
	In Lemma \ref{thmunique}, we prove the unique type of blowup limits. Now we will go a step further in the unique of blowup limits at regular points.

	\begin{lem}\label{lem5} Let $(u,v)$ be a solution of system and $0$ is a regular point.
		\begin{description}
			\item[(i)] If  $0\in\Gamma(u)$, either   $\lambda\geq0$ or \eqref{gammaulambda} holds. Then the blowup limit at $0$ is unique. That is,
			\begin{eqnarray*}
				u_r(x)\to u_0(x)=\frac{1+\lambda v(0)}{2}(x_n^+)^2\ {\rm as}\ r\to0.
			\end{eqnarray*}
			\item[(ii)] If  $0\in\Gamma(v)$, either   $\lambda\geq0$ or \eqref{gammavlambda} holds.  Then the blowup limit at $0$ is unique. That is,
			\begin{eqnarray*}
				v_r(x)\to v_0(x)=\frac{1+\lambda u(0)}{2}(x_n^+)^2\ {\rm as}\ r\to0.
			\end{eqnarray*}
		\end{description}
	\end{lem}
	\pf Without loss of generality we may assume that $u_{r_i}\to u_0:=\frac{1+\lambda v(0)}{2}(x_n^+)^2$ over some sequence $r_i\to 0$. Suppose that now $u_{r'_i}\to u'_0=\frac{1+\lambda v(0)}{2}((x\cdot e')^+)^2$ over another sequence $r'_j\to0$. We claim that there exists $r_\delta>0$ such that as soon as $r'_i<r_\delta$, then
	\begin{eqnarray*}
		\partial_{\tilde{e}}u_{r'_i}\geq0,\  \ \forall \tilde{e}\in \mathcal{C}_\delta\cap\partial B_1.
	\end{eqnarray*}
	In fact, since $\tilde{e}_n\geq \frac{\delta}{2}$ with $e=(\tilde{e}_1,\cdots,\tilde{e}_n)$. We obtain, in $B_1$,
	\begin{eqnarray*}
		\delta^{-1}\partial_{\tilde{e}} u_0-u_0&=&\delta^{-1}\nabla u_0\cdot\tilde{e}-u_0\\
		&=&\delta^{-1}(1+\lambda v(0))x_n^+\tilde{e}_n-\frac{1+\lambda v(0)}{2}(x_n^+)^2\\
		&=&\frac{1+\lambda v(0)}{2}x_n^+\left(\frac{2}{\delta}\tilde{e}_n-x_n^+\right)\\
		&\geq&0,
	\end{eqnarray*}
	where using the fact $\tilde{e}\in C_\delta\cap \partial B_1$. Since $u_{r_i}\to u_0$ in $C^{1,\alpha}_\loc(\Rn)$, for  $\epsilon<\frac{\delta}{64mn}$ with  suitable $m>1$ depending on the sign of $\lambda$, we have $\displaystyle \sup_{B_1}|u_{r_i}-u_0|<\epsilon, \sup_{B_1}|\nabla u_{r_i}-\nabla u_0|<\epsilon$ as $i$ large. It follows that, in $B_1$,
	\begin{eqnarray}
		\delta^{-1}\partial_e u_{r_i}-u_{r_i}\geq-\delta^{-1}|\nabla u_{r_i}-\nabla u_0|-|u_{r_i}-u_0|
		\geq-2\delta^{-1}\epsilon.\label{6.1}
	\end{eqnarray}
	Now we claim that $\delta^{-1}\partial_e u_{r_i}-u_{r_i}\geq0$ in $B_{1/2}, \forall e\in C_\delta\cap \partial B_1$ with choosing some $i$. Suppose the claim is false, let $y_i\in B_{1/2}\cap \Omega(u_{r_i})$(if $y_i\in (\Omega(u_{r_i}))^c$ the conclusion is obvious) be such that $\delta^{-1}\partial_e u_{r_i}(y_i)-u_{r_i}(y_i)<0$. Consider then the auxiliary function
	\begin{eqnarray*}
		w(x)=\delta^{-1}\partial_e u_{r_i}(x)-u_{r_i}+\frac{1}{4mn}|x-y_i|^2.
	\end{eqnarray*}
	Now we prove that $w$ is super-harmonic in $B_{1/2}(y_i)\cap \Omega(u_{r_i})$ for large $i$ and suitable $m>1$ depending on the sign of $\lambda$.
	
	Assume $0\in \Omega(v)$. Then in $D:=B_{1/2}(y_i)\cap \Omega(u_{r_i})$, as $i$ large, formally,
	\begin{eqnarray*}
		\Delta w=\delta^{-1}\lambda \partial_e(v(r_i x))-\lambda v(r_i x)-\frac{1}{2}=\delta^{-1}\lambda \nabla v(r_i x)(r_ie)-\lambda v(r_i x)-\frac{2m-1}{2m}.
	\end{eqnarray*}
	
	As $\lambda\geq0$, since $v(0)>0$, there exists $i_1$ such that $v(r_ix)\geq \frac{v(0)}{2}$ for $i\geq i_1$. It follows from $v\in C^{1,1}_\loc(\Omega)$ that for $i$ large, $|\nabla v(r_ix)|\leq C$, $x\in B_{1/2}(y_i)$. Then there exists $i_2(\delta, v)$ such that $\delta^{-1}\nabla v(r_ix)r_ie\leq \frac{v(0)}{4}$ as $i\geq i_2(\delta, v)$. Thus, we can find $i_0:=i_0(\delta,v)$ such that, as $i\geq i_0$, $m>1$,
	\begin{eqnarray}
		\lambda\left(\delta^{-1}\partial_e(v(r_i x))-v(r_i x)\right)-\frac{2m-1}{2m}\leq0\label{5.8}
	\end{eqnarray}
	is obvious.
	
	As $\lambda<0$, since $\displaystyle \lim_{i\to\infty}\delta^{-1}\lambda \nabla v(r_i x)(r_ie)=0, \lim_{i\to\infty}\lambda v(r_i x)<1$, hence, there exists $i_3:=i_3(\delta,v)$ such that as $i\geq i_3$, choosing a suitable $m>1$ such that \eqref{5.8} holds.
	
	Assume that $v(0)=0$, then formally $\Delta w=\lambda \delta^{-1}\nabla v(r_i x)(r_ie)-\lambda v(r_i x)-\frac{2m-1}{2m}$. It follows from $\displaystyle \lim_{r_i\to0}(\delta^{-1}\partial_e (v(r_ix))-v(r_ix))=0$ that there exists $i_4:=i(v)$ such that for $i\geq i_4$, suitable $m$ depending on the sign $\lambda$,
	\begin{eqnarray*}
		\delta^{-1}\nabla v(r_i x)(r_ie)-v(r_i x)-\frac{2m-1}{2m}\leq0.
	\end{eqnarray*}
	Therefore, there exists $m$ such that $w$ is super-harmonic in $B_{1/2}(y_i)\cap \Omega(u_{r_i})$ as $i$ large. Combining with $w(y_i)<0$ and $w\geq0$ on $\partial\Omega(u_{r_i})$, by the minimum principle, $w$ has a negative infimum on $\partial B_{1/2}(y_i)\cap\Omega(u_{r_i})$, i.e.
	\begin{eqnarray*}
		\displaystyle \inf_{\partial B_{1/2}(y_i)\cap \Omega(u_{r_i})} w<0.
	\end{eqnarray*}
	This can be rewritten as
	\begin{eqnarray*}
		\displaystyle \inf_{\partial B_{1/2}(y_i)\cap \Omega(u_{r_i})} (\delta^{-1}\partial_e u_{r_i}(x)-u_{r_i})<-\frac{1}{16m n}.
	\end{eqnarray*}
	This contradicts \eqref{6.1} with $\epsilon<\frac{1}{64 m n}$. Thus there exists $r_\delta(\lambda, v)$ such that as $r_i<r_\delta(\lambda, v)$, we have $\partial_e u_{r_i}(x)\geq0$ in $B_{1/2}$, $\forall e\in \mathcal{C}_\delta\cap \partial B_1$ and we obtain that $\partial_e u_{r_i'}(x)\geq0$ in $B_{1/2}$, $\forall e\in \mathcal{C}_\delta\cap \partial B_1$ and $r_i'<r_\delta(\lambda, v)$. Therefore passing to a limit we obtain that $\partial_e u'_0(x)\geq0$, $\forall e\in \mathcal{C}_\delta\cap \partial B_1$. This is equivalent to saying that $e'_0\cdot e\geq0$. Letting $\delta\to0$, we see that $e'_0$ necessarily coincides with $e_n$ and therefore $u'_0=u_0$. The results for the case $v$ can be obtained by  the similar proof above. \hbx
	
	\begin{pro}[Lipschitz regularity]\label{lippro}
		Let $(u,v)$ be a solution of system \eqref{system} and $0$ be a regular point.
		\begin{description}
			\item[(i)] If $0\in\Gamma(u)$ and $u_0=\frac{1+\lambda v(0)}{2}(x_n^+)^2$ is a blowup at $0$. As $\lambda\geq0$ or \eqref{gammaulambda} holds.  Then there exists $\rho=\rho(u,v,\lambda)>0$ and a Lipschitz function $f:B'_\rho\to R$ such that
			\begin{eqnarray*}
				\Omega(u)\cap B_\rho=\{x\in B_\rho: x_n>f(x')\},\\
				\Gamma(u)\cap B_\rho=\{x\in B_\rho: x_n=f(x')\}.
			\end{eqnarray*}
			Moreover, for $\delta\in(0,1]$, if $r_\delta^9:=r_\delta(\lambda,v)$ is as in the proof of Lemma {\rm\ref{lem5}}, then $|\nabla_{x'}f|\leq\delta$ a.e. on $B'_{r_\delta^9/2}$.
			\item[(ii)] If $0\in\Gamma(v)$ and $u_0=\frac{1+\lambda u(0)}{2}(x_n^+)^2$ is a blowup at $0$. As $\lambda\geq0$ or \eqref{gammavlambda} holds. Then there exists $\rho=\rho(v,u,\lambda)>0$ and a Lipschitz function $f:B'_\rho\to R$ such that
			\begin{eqnarray*}
				\Omega(v)\cap B_\rho=\{x\in B_\rho: x_n>f(x')\},\\
				\Gamma(v)\cap B_\rho=\{x\in B_\rho: x_n=f(x')\}.
			\end{eqnarray*}
			Moreover, for $\delta\in(0,1]$, if $r_\delta^{10}:=r_\delta(\lambda,u)$ is as in the proof of Lemma {\rm\ref{lem5}}, then $|\nabla_{x'}f|\leq\delta$ a.e. on $B'_{r_\delta^{10}/2}$.
		\end{description}
	\end{pro}
	\pf Let $0\in\Gamma(u)$, set $u_{r_k}$ is so that $u_0$ is a blowup limit. Then for $i$ large, $\displaystyle \sup_{B_1}|u_{r_i}-u_0|\leq\epsilon,\ \sup_{B_1}|\nabla u_{r_i}-\nabla u_0|\leq\epsilon$. Since $0\in\Gamma(u)$, then $u_{r_i}(0)=0$. For any $e\in \mathcal{C}_\delta\cap \partial B_1$, we can obtain that $\delta^{-1}\partial_e u_0-u_0\geq0$ from the proof of Lemma \ref{lem5} and we know that for $i$ large $\delta^{-1}\partial_e u_{r_i}-u_{r_i}\geq0$ in $B_{1/2}$ as $e\in \mathcal{C}_\delta\cap\partial B_1$. We immediately get
	\begin{eqnarray*}
		u_{r_i}\geq0\ {\rm in}\ \mathcal{C}_\delta\cap B_{1/2},\ \ \  u_{r_i}\leq0\ {\rm in}\ -\mathcal{C}_\delta\cap B_{1/2}.
	\end{eqnarray*}
	On the other hand, $u_{r_i}\geq0$ in $B_{1/2}$ and therefore $u_{r_i}=0$ in $-C_\delta\cap B_{1/2}$. Now we assume that $0\in \Gamma(u_{r_i})\cap B_{1/2}$ and which implies $u_{r_i}>0$ in $\mathcal{C}_\delta\cap B_{1/2}$. Indeed, otherwise there exists $y\in \mathcal{C}_\delta\cap B_{1/2}$ with $u_{r_i}(y)=0$, and from the argument above, $u_{r_i}=0$ in $(y-\mathcal{C}_\delta)\cap B_{1/2}$. The latter set, however, is a neighborhood of $0$, which contradicts the fact that $0\in \Gamma(u_{r_i})$. Hence $u_{r_i}>0$ in $\mathcal{C}_\delta\cap B_{1/2}$, since
	\begin{eqnarray*}
		\sup_{B_1}|u_{r_i}-u_0|<\epsilon,\ \sup_{B_1}|\nabla u_{r_i}-\nabla u_0|\leq\epsilon
	\end{eqnarray*}
	for large $i$. If $\epsilon<\frac{\delta}{64mn}$ with suitable $m$ depending on the sign of $\lambda$, see Lemma \ref{lem5}, $\delta\in(0,1]$, we get that $\partial_e u_{r_i}\geq0$ in $B_{1/2}$ for $e\in \mathcal{C}_\delta\cap\partial B_1$ as $r_i<r_\delta^9$. The obvious fact $u_{r_i}(0)=0$ shows that the set $\{u_{r_i}=0\}$ is nonempty. Then there exists a function  $g(x')$ on $B'_{1/2}$ such that $\{u_{r_i}>0\}=\{x\in B_{1/2}|\ x_n>g(x')\}$, which is equivalent to $\{u(x)>0\}=\{y\in B_{\frac{1}{2r_i}}| y>f(x')\}$ with $f(x')=r_ig(\frac{x'}{r_i})$, and $\Gamma(u)=\{x\in B_\rho:\ x_n=f(x')\}$. Moreover, since $\partial_eu_{r_i}\geq0$ in $B_{1/2}$, we have $\partial_e u\geq0$ in $B_{r_i/2}$, using \cite[Exercise 4.1]{2012PetrosyanShahgholianUraltseva}, we get that $|\nabla_{x'}g|\leq\delta$, which implies that $|\nabla_{x'}f(x')|\leq\delta$. We can obtain the results of $v$ replacing $u$ by $v$. \hbx
	
	\begin{lem}\label{lem4}
		Let $(u,v)$ be a solution of system \eqref{system} such that $0$ is a regular point.
		\begin{description}
			\item[(i)] If  $0\in\Gamma(u)$, as $\lambda\geq0$ or \eqref{gammaulambda} holds. Then there exists $\rho>0$ such that all points in $\Gamma(u)\cap B_\rho$ are also regular.
			\item[(ii)] If  $0\in\Gamma(v)$, as $\lambda\geq0$ or \eqref{gammavlambda} holds. Then there exists $\rho>0$ such that all points in $\Gamma(v)\cap B_\rho$ are also regular.
		\end{description}
	\end{lem}

	\pf By the proof of Lemma \ref{lem5}, we get that $\partial_e u_0-u_0\geq0$ in $B_1$ for $e\in \mathcal{C}_\delta\cap \partial B_1$, furthermore, $\partial_e u_{r_i}\geq0$ in $B_{1/2}$ for $e\in \mathcal{C}_\delta\cap \partial B_1$. Thus $\partial_e u(x)\geq0$ in $B_{\frac{r_i}{2}}$ for any $e\in \mathcal{C}_\delta\cap \partial B_1$ for large $i$. Now fixed $i_0$ such that the inequality above holds, that is $\partial_e u(x)\geq0$ in $B_{r_0}$ for any $e\in \mathcal{C}_\delta\cap \partial B_1$. Thus for any $z\in \Gamma(u)\cap B_{r_0}$, $x\in \Rn$, then $r_ix+z\to z$ as $i\to\infty$, Thus
	\begin{eqnarray*}
		\partial_e(u_{z,r_i}(x)) =\frac{\partial_e u(r_ix+z)r_i}{r_i^2}\geq0
	\end{eqnarray*}
	for $x\in \Rn, e\in  \mathcal{C}_\delta\cap \partial B_1$. Which leads to the limit $u_0^z$ of $u_{z,r_i}$ satisfies that $\partial_e u_0^z(x)\geq0$ in $\Rn$ for any $z\in B_{r_0}\cap \Gamma(u)$ and $e\in  \mathcal{C}_1$. We claim that $u_0^z$ must ba a half-space solution, since no homogeneous quadratic polynomial can be monotone in a cone of direction $\mathcal{C}_1$. Thus the lemma holds taking $\rho=r_0$. Replacing $u$ with $v$, we can prove the conclusion of $v$. \hbx
	
	\begin{thm}[$C^1$ regularity]\label{thmC1}
		Let $(u,v)$ be a solution of system \eqref{system} satisfying $0$ is a regular point.
		\begin{description}
			\item[(i)] If  $0\in\Gamma(u)$, either    $\lambda\geq0$ or \eqref{gammaulambda} holds.  Then there exists $\rho>0$ such that $\Gamma(u)\cap B_\rho$ is a $C^1$ graph.
			\item[(ii)] If  $0\in\Gamma(v)$, either   $\lambda\geq0$ or \eqref{gammavlambda} holds.  Then there exists $\rho>0$ such that $\Gamma(v)\cap B_\rho$ is a $C^1$ graph.
		\end{description}
	\end{thm}
	\pf We first consider the case $0\in\Gamma(u)$, we may assume that $u_0=\frac{1+\lambda v(0)}{2}(x_n^+)^2$ is the blowup at the origin. Then by Proposition \ref{lippro}, $\Gamma(u)\cap B_\rho$ is a Lipschitz graph $x_n=f(x')$ with $f(0)=0$ and $|f(x')|\leq\delta |x'|$ in $B'_{r_\delta/2}$. Since we can choose $\delta>0$ arbitrarily small, we immediately obtain the existence of tangent plane to $\Gamma(u)$ at the origin with $e_n$ being the normal vector.
	
	Now by Lemma \ref{lem4}, every free boundary point $z\in \Gamma(u)\cap B_\rho$ is regular and therefore $\Gamma(u)$ has a tangent plane at all these points. Let $\nu_z$ denote the unit normal vector at $z\in \Gamma(u)$ pointing into $\Omega(u)$. Then the proof of   Lemma \ref{lem5} implies that $\nu_z\cdot e\geq0$ for any $e\in \mathcal{C}_\delta$ if $z\in \Gamma(u)\cap B_{r_\delta}$. This means that $\nu_z$ is from the dual of the cone $\mathcal{C}_\delta$, that is $\nu_z\in \mathcal{C}_{1/\delta}$. In particular, we obtain
	\begin{eqnarray*}
		|\nu_z-e_n|\leq C\delta,\ z\in \Gamma(u)\cap B_{r_\delta}.
	\end{eqnarray*}
	The latter means that $\Gamma(u)$ is $C^1$ at the origin. It is obviously true for any $z\in \Gamma(u)\cap B_\rho$.   This completes the proof by proving the results of $v$ based the proof of $u$.\hbx

By carefully examining the monotonicity of the solutions  for a cone of directions, we can also prove that the second derivative of the solutions is continuous up to the free boundary. The results are presented as follows.	
	\begin{pro}\label{prohighregular}
		Let $(u,v)$ be a solution of system \eqref{system} with $0$ a regular point.
		\begin{description}
			\item[(i)] If  $0\in\Gamma(u)$ and $u_0(x)=\frac{1+\lambda v(0)}{2}(x_n^+)^2$ is a blowup of $u$ at $0$, assume $\lambda\geq0$ or \eqref{gammaulambda} holds. Then
			\begin{eqnarray*}
				\lim_{x\in \Omega(u)\atop x\to 0}\partial_{x_i x_j}u(x)=\delta_{in}\delta_{jn}=\partial_{x_ix_j}\left(\frac{1+\lambda v(0)}{2}x_n^2\right).
			\end{eqnarray*}
			\item[(ii)] If  $0\in\Gamma(v)$ and $v_0(x)=\frac{1+\lambda u(0)}{2}(x_n^+)^2$ is a blowup of $v$ at $0$, assume $\lambda\geq0$ or \eqref{gammavlambda} holds. Then
			\begin{eqnarray*}
				\lim_{x\in \Omega(v)\atop x\to 0}\partial_{x_i x_j}v(x)=\delta_{in}\delta_{jn}=\partial_{x_ix_j}\left(\frac{1+\lambda u(0)}{2}x_n^2\right).
			\end{eqnarray*}
		\end{description}
		
	\end{pro}
	\pf Let $0\in\Gamma(u)$ and $x_j\in\Omega(u)$ be a sequence converging to $0$, $d_j=\dist(x^j,\Gamma(u))$ and $y^j\in\Gamma(u)\cap\partial B_{d_j}(x^j)$. Also let $\xi^j=(x^j-y^j)/{d_j}$. Considering then the rescalings of $u$ at $y^j$,
	\begin{eqnarray*}
		w_j(x)=u_{y^j,d_j}(x)=\frac{u(y^j+d_j x)}{d_j^2}.
	\end{eqnarray*}
	Since
	\begin{eqnarray*}
		D^2w_j(\xi^j)=D^2u(y^j+d_j\xi^j)=D^2u(x^j).
	\end{eqnarray*}
	We will show that
	\begin{eqnarray*}
		D^2 w_j(\xi^j)\to D^2\left(\frac{1+\lambda v(0)}{2}x_n^2\right).
	\end{eqnarray*}
	According to Section {\rm\ref{sec5.2}}, note that we have  uniform $C^{1,1}$ estimates locally in $\Rn$ for the family $\{w_j\}$ and thus we may assume that $w_j\to w_0\in C^{1,1}_\loc(\Rn)$ in $C^{1,\alpha}_\loc(\Rn)$. The fact $|\xi_i|=1$ implies that we may assume that $\xi^j\to \xi^0$. Since the functions $w_j$ satisfies $\Delta w_j=1+v(y^j+d_j x)$ in $B_1(\xi^j)$, we claim that $\Delta w_0=1+v(0)$ in $B_{2/3}(\xi^0)$, which is equivalent to
	\begin{eqnarray}
		-\int_{B_{2/3}(\xi^0)}\nabla w_0\nabla \varphi dx=\int_{B_{2/3}(\xi^0)}(1+v(0))\varphi dx\label{6.2}
	\end{eqnarray}
	for any $\varphi\in C_0^\infty(B_{2/3}(\xi^0))$. Firstly,
	\begin{eqnarray*}
		\int_{B_{2/3}(\xi^0)}\Delta w_j \varphi dx=\int_{B_{2/3}(\xi^0)}(1+v(y^j+d_jx))\varphi dx
	\end{eqnarray*}
	for $j$ large. As $j$ large, then $B_{2/3}(\xi^0)\subset B_1(\xi^j)$ and $\varphi(x)=0$ for $B_1(\xi^j)\setminus B_{2/3}(\xi^0)$, it follows that
	\begin{eqnarray*}
		\int_{B_1(\xi^j)}\Delta w_j\varphi dx=\int_{B_1(\xi^j)}(1+v(y^j+d_j x))\varphi dx.
	\end{eqnarray*}
	Thus  passing to the limit,  \eqref{6.2} holds. It is easy to see that
	\begin{eqnarray*}
		\Delta(w_j-w_0)=\lambda v(y^j+d_j x)-\lambda v(0)\ {\rm in}\ B_{2/3}(\xi^0).
	\end{eqnarray*}
	Using   \cite[Theorem 1.1]{2012PetrosyanShahgholianUraltseva}, we have, for $p>1$,
	\begin{eqnarray*}
		\|w_j-w_0\|_{W^{2,p}(B_{1/2}(\xi^0))}\leq C\left(\|w_j-w_0\|_{L^1(B_{2/3}(\xi^0))}+|\lambda|\|v(y^j+d_jx)-v(0)\|_{L^p(B_{2/3}(\xi^0))}\right).
	\end{eqnarray*}
	It follows that
	\begin{eqnarray*}
		w_j\to w_0\ {\rm in}\ W^{2,p}(B_{1/2}(\xi^0)).
	\end{eqnarray*}
	For the equation $\Delta w_j=1+v(y^j+d_j x)$ in $B_{2/3}(\xi^0)$ with $j$ being large, using   \cite[Corollary 2.16]{2022FR}, we get that $w_j\in C^2(B_{1/2}(\xi^0))$ and in particular that $D^2 w_j(\xi^j)\to D^2 w_0(\xi^0)$. Thus we will be done once we show that
	\begin{eqnarray}
		w_0(x)=\frac{1+\lambda v(0)}{2}(x_n^+)^2\ {\rm and}\ \xi^0=e_n.\label{5.10}
	\end{eqnarray}
	In order to prove \eqref{5.10} above, we recall that by Lemma \ref{lem5}, for $\delta>0$, there exists $r_\delta:=r_\delta(\lambda,v)>0$ such that $\partial_e u\geq0$ in $B_{r_\delta}$ for any direction $e\in \mathcal{C}_\delta$. This immediately implies that for any $R>0$, we have $\partial_e w_j\geq0$ in $B_R$ provided $j$ is sufficiently large. Passing to the limit, we obtain therefore that
	\begin{eqnarray*}
		\partial_e w_0\geq0\ {\rm in}\ \Rn
	\end{eqnarray*}
	for any direction $e$ with $e\cdot e_n>0$. This directly implies that $w_0(x)=\phi(x_n)$ for a monotone increasing function of $x_n$, since $0$ is a free boundary point, the only such function is $\frac{1+\lambda v(0)}{2}(x_n^+)^2$ as claimed. Finally, $|\xi^0|=1$ and $\Delta w_0=1+\lambda v(0)$ in $B_1(\xi^0)$, the only such point for $w_0(x)=\frac{1+\lambda v(0)}{2}(x_n^+)^2$ is $\xi^0=e_n$. The results for $v$ are similar. \hbx
	
\subsubsection*{Higher regularity  of the free boundary}\label{sec5.3}
The regularity results of the free boundary discussed above depend on specific solutions. Similar to  \cite[Theorem 5.1]{2012PetrosyanShahgholianUraltseva}, we can prove that the global solutions of system \eqref{system} are convex. Next, we will start from the convexity of the global solutions to prove the uniform Lipschitz and $C^{1,\alpha}$ regularity of the free boundary. Now we give a  result that the solutions of system can be approximated by convex global solutions in $\Rn$.
	\begin{lem}\label{lem7.1} Let $(u,v)$ be a solution of system \eqref{system} with  $0\in\Gamma(u,v)$. For any ball $B_R\Subset \Omega$, fix $\sigma\in(0,1]$ and $\forall \epsilon>0$, there exist  radius $r^1:=r(\sigma,\|D^2 u\|_{L^\infty(B_R)},n)$, $r^2:=r(\sigma,\|D^2 v\|_{L^\infty(B_R)},n)$ such that if
		\begin{description}
			\item[(i)]
			$0\in\Gamma(u_{r_i}), r<r^1, \delta(1,u_{r_i})\geq\sigma$.
			Then there exists a convex  global solution $u_0$ of equation
			\begin{eqnarray*}
				\Delta w=1+\lambda v(0) \ {\rm in}\ \Rn
			\end{eqnarray*}
			with $\|D^2 u_0\|_{L^\infty(\Rn)}\leq \|D^2 u\|_{L^\infty(B_R)}$ such that
			\begin{enumerate}
				\item $\|u_{r_i}-u_0\|_{C^1(B_1)}\leq \epsilon$;
				\item there exists a ball $B_\rho(x)\subset B_1$ of radius $\rho=\frac{\sigma}{2n}$ such that $B_\rho(x)\subset \Omega^c(u_0)$.
			\end{enumerate}
			\item[(ii)]
			$0\in\Gamma(v_{r_i}),  r<r^2, \delta(1,v_{r_i})\geq\sigma$. Then there exists a convex  global solution $v_0$ of equation
			\begin{eqnarray*}
				\Delta w=1+\lambda u(0) \ {\rm in}\ \Rn
			\end{eqnarray*}
			with $\|D^2 v_0\|_{L^\infty(\Rn)}\leq \|D^2 v\|_{L^\infty(B_R)}$ such that
			\begin{enumerate}
				\item $\|v_{r_i}-v_0\|_{C^1(B_1)}\leq \epsilon$;
				\item there exists a ball $B_\rho(x)\subset B_1$ of radius $\rho=\frac{\sigma}{2n}$ such that $B_\rho(x)\subset \Omega^c(v_0)$.
			\end{enumerate}
		\end{description}
		
	\end{lem}
	\pf Fix $\epsilon>0$ and argue by contradiction. For the case $0\in\Gamma(u)$, if no such $r^1$ exists, then one can find a sequence $r_i\to \infty$ and scaling functions $u_{r_i}$ such that
	\begin{eqnarray}
		\|u_{r_i}-u_0\|_{C^1(B_1)}>\epsilon\ \ {\rm for\ any}\ u_0\ {\rm with}\ \|D^2 u_0\|_{L^\infty(\Rn)}\leq \|D^2 u\|_{L^\infty(B_R)}. \label{7.1}
	\end{eqnarray}
	Then we can obtain the uniform estimates $|u_{r_i}(x)|\leq \frac{\|D^2 u\|_{L^\infty(B_R)}}{2}|x|^2,\ |x|\leq R/{r_i}$. We may assume $u_{r_i}$ converges in $C_\loc^1(\Rn)$ to a global solution $\tilde{u}_0$ with $\|D^2 \tilde{u}_0\|_{L^\infty(\Rn)}\leq \|D^2 u\|_{L^\infty(B_R)}$, but then
	\begin{eqnarray*}
		\|u_{r_i}-\tilde{u}_0\|_{C^1(B_1)}<\epsilon\ {\rm for}\ i\geq i_\epsilon,
	\end{eqnarray*}
	contradicting \eqref{7.1}.
	
	Applying the John's ellipsoid lemma \cite[Lemma 5.14]{2012PetrosyanShahgholianUraltseva}, to the convex set $D=\Omega^c(\tilde{u}_0)\cap B_1$, we obtain that $D$ contains a ball $B$ of radius $\rho=\sigma/2n$. Indeed, if $E$ is one of its diameters smaller than $2\rho$, then $nE$ is one of its diameters smaller than $2n\rho=\sigma$ and $D$ is contained in a strip of width smaller than $\sigma$, a contradiction. Thus, conclusion $(i)-2$ is also satisfied for $\tilde{u}_0$, which means we have arrived at a contradiction with the assumption that no $r^1$ exists. The results for $v$ and $0\in\Gamma(v)$ can be confirmed from the process above easily.  \hbx

	\begin{lem}\label{thm8.1} Let $(u,v)$ be a solution of system \eqref{system} with $0\in\Gamma(u,v), B_R\Subset \Omega$.
		For every $\sigma\in(0,1]$. Then
		\begin{itemize}
			\item if $0\in\Gamma(u)$, $\delta(1,u_{r_i})\geq\sigma$, there exists $r_\sigma^3=r_\sigma(|D^2 u|_{L^\infty(B_R)},n)$ such that  as $r_i<r_\sigma^3$, then $\Gamma(u)\cap B_\rho$ is Lipschitz regular with a Lipschitz constant $L=L(\sigma,n,\|D^2 u\|_{L^\infty(B_R)})$ and $\rho$ depends on $n,\sigma,r_\sigma^3$.
			\item If $0\in\Gamma(v)$, $\delta(1,v_{r_i})\geq\sigma$, there exists $r_\sigma^4=r_\sigma(\|D^2 v\|_{L^\infty(B_R)},n)$ such that  as $r_i<r_\sigma^4$, then $\Gamma(v)\cap B_\rho$ is Lipschitz regular with a Lipschitz constant $L=L(\sigma,n,\|D^2 v\|_{L^\infty(B_R)})$ and $\rho$ depends on $n,\sigma,r_\sigma^4$.
		\end{itemize}
	\end{lem}
	\pf Fix a small $\epsilon=\epsilon(\sigma,n,\|D^2 u\|_{L^\infty(B_R)})>0$ and apply Lemma \ref{lem7.1}, so, if $r<r^1_\sigma$, we can find a convex global solution $u_0$ such that $\Omega^c(u_0)\cap B_1$ must contain a ball $B$ of radius $\rho=\frac{\sigma}{2n}$. Without loss of generality, we may assume that $B=B\rho(-se_n)$ for $\rho\leq s\leq 1-\rho$ and we will also have $B_{\rho/2}(-se_n)\subset \Omega^c(u_{r_i})\cap B_1$ using Proposition \ref{prolimit}, we now apply \cite[Lemma 6.1]{2012PetrosyanShahgholianUraltseva}, which gives that
	\begin{eqnarray*}
		C_0\partial_e u_0-u_0\geq0\ {\rm in}\ K(\rho/16,s,1/2)\ {\rm for\ any}\ e\in \mathcal{C}_{8/\rho}.
	\end{eqnarray*}
	Further, if $\epsilon$ is small enough, the approximation $\|u_{r_i}-u_0\|_{C^1(B_1)}\leq \epsilon$ implies that $C_0\partial_e u_{r_i}-u_{r_i}\geq -\frac{(\rho/16)^2}{8n}$ in $K(\rho/16,s,1/2)$. Using the proof of Lemma \ref{lem5}, we obtain $C_0\partial_e u_{r_i}-u_{r_i}\geq0$ in  $K(\rho/32,s,1/4)$(To be more accurate, one needs to apply  the proof of Lemma \ref{lem5} in every ball $B_{\rho/16}(te_n)$ for $t\in[-s,\frac{1}{4})$). The latter inequality can be rewritten as
	\begin{eqnarray*}
		\partial_e(e^{-C_0(x\cdot e)}u_{r_i})\geq0
	\end{eqnarray*}
	in $K(\rho/32,s,1/4)$. Taking $e=e_n$ and noting that $u_{r_i}=0$ on $B'_{\rho/32}\times \{-s\}$, we obtain after integration that
	\begin{eqnarray*}
		u_{r_i}\geq0\ {\rm in}\  K(\rho/32,s,1/4).
	\end{eqnarray*}
	Combining with the previous inequality, this gives
	\begin{eqnarray*}
		\partial_e u_{r_i}\geq 0 \ {\rm in}\  K(\rho/32,s,1/4)\ {\rm for\ any}\ e\in \mathcal{C}_{8/\rho}.
	\end{eqnarray*}
	The rest of the proof follows from the process of Proposition \ref{lippro}. Moreover, the proofs of $v$ are similar.\hbx
	
	\begin{lem}[Lipschitz regularity]\label{lem8.2}
		Let $(u,v)$ be a solution of system with $B_R\Subset \Omega$ and $0\in\Gamma(u,v)$.
		\begin{itemize}
			\item If $0\in \Gamma(u)$,  there exists a modulus of continuity $\sigma^1(r):=\sigma_{\|D^2 u\|_{L^\infty(B_R(x^0))},n}(r)$ such that as $\delta(r,u)\geq \sigma^1(r)$ for some value $r=r_0\in(0,1)$, then $\Gamma\cap B_\rho$ is a Lipschitz graph with a Lipschitz constant $L\leq L(n,\|D^2 u\|_{L^\infty(B_R(x^0))},r_0)$ with $\rho$ depending on $n, r_0, \sigma^1(r_0)$.
			\item If $0\in \Gamma(v)$, then there exists a modulus of continuity $\sigma^2(r):=\sigma_{\|D^2 v\|_{L^\infty(B_R(x^0))},n}(r)$ such that as $\delta(r,v)\geq \sigma^2(r)$ for some value $r=r_0\in(0,1)$, then $\Gamma\cap B_\rho$ is a Lipschitz graph with a Lipschitz constant $L\leq L(n,\|D^2 v\|_{L^\infty(B_R(x^0))},r_0)$  with $\rho$ depending on $n, r_0, \sigma^2(r_0)$.
		\end{itemize}
		
	\end{lem}
	\pf Let $0\in\Gamma(u)$. Note that in Lemma \ref{thm8.1} one can take the function $\sigma\mapsto r_\sigma^3$ to be monotone and continuous in $\sigma$ and such that $\displaystyle\lim_{\sigma\to 0^+}r_\sigma^3=0$.  Now if $\delta(r_0,u)\geq \sigma^1(r_0)$, then the rescaling $u_{r_0}$ satisfies
	\begin{eqnarray*}
		\|D^2 u_{r_0}\|_{L^\infty(B_{R/r_0})}\leq \|D^2 u\|_{L^\infty(B_R(x^0))},\ \delta(1,u_{r_0})\geq \sigma^1(r_0).
	\end{eqnarray*}
	Applying Lemma \ref{thm8.1}. Thus scaling back to $u$ the corresponding statement for the free boundary of $u$ can be established. The results of $0\in\Gamma(v)$ are analogous proof.\hbx

We will apply the boundary Hanack principle in \cite{2019AllenShaholian} to obtain the following results about $C^{1,\alpha}$ regularity.
	\begin{thm}[$C^{1,\alpha}$ regularity]\label{thm8.3} Let $(u,v)$ be a solution of system \eqref{system} with $B_R\Subset\Omega, 0\in\Gamma(u,v)$. For every $\sigma>0$.
		\begin{itemize}
			\item If $0\in\Gamma(u)$, there exists $r^5_\sigma=r^5_\sigma(\|D^2u\|_{L^\infty(B_R)},n)$ such that as  $\delta(1,u_r)\geq\sigma$, for $r<r^5_\sigma$, then $\Gamma(u)\cap B_\rho$ is a $C^{1,\alpha}$ graph with $\alpha=\alpha({\sigma,\|D^2u\|_{L^\infty(B_R)},n})\in (0,1)$, the $C^{1,\alpha}$ norm $c\leq c(\sigma,n,\|D^2u\|_{L^\infty(B_R)})$ and $\rho$ depends on $n,\sigma,r^5_\sigma, \|D^2u\|_{L^\infty(B_R)},\lambda, v$.
			\item If $0\in\Gamma(v)$, there exists $r^6_\sigma=r^6_\sigma(\|D^2v\|_{L^\infty(B_R)},n)$ such that as $\delta(1,v_r)\geq\sigma$, for $r<r^6_\sigma$, then $\Gamma(v)\cap B_\rho$ is a $C^{1,\alpha}$ graph with $\alpha=\alpha({\sigma,\|D^2v\|_{L^\infty(B_R)},n})\in (0,1)$, the $C^{1,\alpha}$ norm $c\leq c(\sigma,n,\|D^2v\|_{L^\infty(B_R)})$ and $\rho$ depends on $n,\sigma,r^6_\sigma, \|D^2v\|_{L^\infty(B_R)},\lambda, u$.
		\end{itemize}
	\end{thm}
	\pf Since $v\in C^{1,1}_\loc (\Omega)$, using the proof of \cite[$\S$1.4.2]{2019AllenShaholian}, the results hold.\hbx		
\begin{rem}\label{remthm1.6}
In fact, the proof of Theorem \ref{thmC1Alpha} can be established combining with Theorem \ref{thm8.3} and Proposition \ref{pro9.1} below.
\end{rem}

	\subsection{The structure of the singular set}\label{subsec5.5}
The main purpose of this section is to explore the structure of the singular set of system \eqref{system}. To this end, we first provide a characterization of the singular points from the terms of energy and geometry.
	\begin{pro}\label{pro9.1}
		Let $(u,v)$ be a solution of system \eqref{system}. Then for $x^0\in\Gamma(u)$, either  $\lambda\geq0$ or \eqref{gammaulambda} holds, the following statements are equivalent:
		\begin{description}
			\item[(1)] $x^0\in \Sigma(u)$;
			\item[(2)] $w(x^0)=\left(\frac{1+\lambda v(x^0)}{2}\right)^2\alpha_n$;
			\item[(3)] $\displaystyle \lim_{r\to 0}\delta(r,u,x^0)=0$.
		\end{description}
		Then for $x^0\in\Gamma(v)$, either $\lambda\geq0$ or \eqref{gammavlambda}  holds, the following statements are equivalent:
		\begin{description}
			\item[(4)] $x^0\in \Sigma(v)$;
			\item[(5)] $w(x^0)=\left(\frac{1+\lambda u(x^0)}{2}\right)^2\alpha_n$;
			\item[(6)] $\displaystyle \lim_{r\to 0}\delta(r,v,x^0)=0$.
		\end{description}
	\end{pro}
	\pf We first prove the case $x^0\in\Gamma(u)$, the results of $v$ can be establish by coping   the following process. The equivalence of ${\rm (1)}$ and ${\rm (2)}$ in Proposition \ref{proclass}.
	
	${\rm (1)}\Rightarrow {\rm (3)}$ For given $\sigma>0$, if one has $\delta(r,u,x^0)>\sigma$ for a small enough $r$, then $\partial_e u_{x^0,r}\geq0$ in $B_{\beta\sigma}$ for a cone of directions $e\in \mathcal{C}_{\beta\sigma}$ with $\beta>0$ being constant, which implies that any blowup limit $u_0$ will have the same property and therefore can not be polynomial. This means that $\delta(r,u,x^0)\to0$ as $r\to0$.
	
	${\rm (3)}\Rightarrow {\rm (1)}$. Considering a blowup $\displaystyle u_0=\lim_{j\to\infty}u_{x^0,r_j}$ for some $r_j\to0$. Since $\Delta u_{x^0,r_j}=1+v(x^0+r_j x)$ in $B_1$ except a strip of width $\delta(r_j,u,x^0)\to0$, we will obtain that
	\begin{eqnarray*}
		\Delta u_0=1+v(x^0)\ {\rm a.e.\ in}\ B_1
	\end{eqnarray*}
	and therefore $u_0$ must be polynomial. This implies $x^0\in\Sigma(u)$.\hbx

	\begin{lem}\label{lem9.2}
		Let $(u,v)$ be a solution of system \eqref{system} and $0\in\Sigma(u)\cup \Sigma(v)$.
		\begin{itemize}
			\item If $0\in\Gamma(u)$. Then there exists a modulus of continuity $\sigma^1(r)$, depending only on $n,\lambda,v$, such that
			\begin{eqnarray*}
				\delta(r,u)\leq\sigma^1(r)\ {\rm for\ any}\   0<r<1.
			\end{eqnarray*}
			\item If $0\in\Gamma(v)$. Then there exists a modulus of continuity $\sigma^2(r)$, depending only on $n,\lambda,u$, such that
			\begin{eqnarray*}
				\delta(r,v)\leq\sigma^2(r)\ {\rm for\ any}\   0<r<1.
			\end{eqnarray*}
		\end{itemize}

	\end{lem}
	\pf Let $0\in\Gamma(u)$, for any $\sigma\geq0$, there exists $r_0\in(0,1)$ such that $\delta(r_0,u)>\sigma(r_0)$, by the proof of Lemma \ref{thm8.1}, we get that
	\begin{eqnarray*}
		\partial_e u_0\geq0
	\end{eqnarray*}
	in $B_\rho$, which implies $u_0$ must be not a polynomial. This is a contradiction. The result of $v$ is obvious.\hbx
	
	Next, we consider the class
	\begin{eqnarray*}
		\mathcal{Q}^+(u,x^0):=\{q\mid\ q(x)\geq0, q(x)\ {\rm is\ honogeneous\ quadratic\ polynomial}|\Delta q=1+\lambda v(x^0) \},
	\end{eqnarray*}
	\begin{eqnarray*}
		\mathcal{Q}^+(v,x^0):=\{q\mid\ q(x)\geq0, q(x)\ {\rm is\ honogeneous\ quadratic\ polynomial}|\Delta q=1+\lambda u(x^0) \}.
	\end{eqnarray*}
	In particularly, $\mathcal{Q}^+(u):=\mathcal{Q}^+(u,0), \mathcal{Q}^+(v):=\mathcal{Q}^+(v,0)$. Essentially, the elements of $\mathcal{Q}^+(u,x^0)$,  $\mathcal{Q}^+(v,x^0)$ are the blowups of solutions of system at the singular points $x^0$. Clearly, polynomial $q\in \mathcal{Q}^+(u,x^0)$ or $q\in\mathcal{Q}^+(v,x^0)$ is themself solution to the single equation. Each member of $\mathcal{Q}^+(u,x^0)$ is uniquely represented as $q(x)=\frac{1+\lambda v(x^0)}{2}(x\cdot Ax)$ and the member of $\mathcal{Q}^+(v,x^0)$ is  $q(x)=\frac{1+\lambda u(x^0)}{2}(x\cdot Ax)$  where $A$ is a symmetric $n\times n$ matrix with $\tr A=1$.  Then note that the free boundary
	\begin{eqnarray*}
		\Gamma(q)=\ker A
	\end{eqnarray*}
	consists completely of singular points. The dimension $d=\dim\ker A$ can be any integer from $0$ to $n-1$.
	The case $d=n-1$ deserves a special attention. For definiteness, consider $q_0(x)=\frac{1+\lambda v(x^0)}{2}x_n^2$. Then the free boundary of $u$ is a hyperplane $\Pi=\{x_n=0\}$ of codimension one which is as smooth as it can be. Nevertheless, all points on $\Pi$ are singular, since $\Omega^c(q_0)=\Pi$ is ``thin"(see Proposition \ref{pro9.1}). In contrast, for the half space solution $h_0(x)=\frac{1+\lambda v(x^0)}{2}(x_n^+)^2$, the free boundary is still $\Pi$, however, this time it consists of regular points, since $\Omega^c(h_0)=\{x_m\leq 0\}$is ``thick". This means that the singular set can be as large as the set of regular points.
	
	We give a key property that the solutions of system can be approximated by polynomial solutions.
	\begin{lem}\label{lem9.3}
		Let $(u,v)$ be a solution of system \eqref{system} and $B_\rho\Subset\Omega$.
		\begin{itemize}
			\item If $0\in\Sigma(u)$. Then there exists a modulus of continuity $\sigma^1(r)$ depending only on $u,n,\lambda, v$ such that for any $0<r<\frac{\rho}{2}$, there exists a homogeneous quadratic polynomial $q^r\in \mathcal{Q}^+(u)$, such that
			\begin{eqnarray*}
				\|u-q^r\|_{L^\infty(B_{2r})}\leq \sigma^1(r) r^2,\\
				\|\nabla u-\nabla q^r\|_{L^\infty(B_{2r})}\leq \sigma^1(r) r.
			\end{eqnarray*}
			\item If $0\in\Sigma(v)$. Then there exists a modulus of continuity $\sigma^2(r)$ depending only on $v,n,\lambda, u$ such that for any $0<r<\frac{\rho}{2}$, there exists a homogeneous quadratic polynomial $q^r\in \mathcal{Q}^+(v)$, such that
			\begin{eqnarray*}
				\|v-q^r\|_{L^\infty(B_{2r})}\leq \sigma^2(r) r^2,\\
				\|\nabla v-\nabla q^r\|_{L^\infty(B_{2r})}\leq \sigma^2(r) r.
			\end{eqnarray*}
		\end{itemize}

	\end{lem}
	\pf Let $0\in\Gamma(u)$. Equivalently, we can show that for any $\epsilon>0$, there exists $r_\epsilon>0$ depending only on $\rho, u,\lambda, v$ such that as $0<r\leq r_\epsilon$, we can find $q^r\in \mathcal{Q}^+(u)$ satisfies
	\begin{eqnarray*}
		|u-q^r|\leq \epsilon r^2,\ |\nabla u-\nabla q^r|\leq \epsilon r\ {\rm in}\ B_{2r}.
	\end{eqnarray*}
	Arguing by contradiction, assume that this statement fails for some $\epsilon_0>0$, which implies that there exists a sequences $r_i\to0$ and  for any $q\in \mathcal{Q}^+(u)$, either one or the other of the following holds,
	\begin{eqnarray*}
		\|u_i-q^r\|_{L^\infty(B_{2r_i})}\geq \epsilon_0r_i^2,\ \
		\|\nabla u_i-\nabla q^r\|_{L^\infty(B_{2r_i})}\geq \epsilon_0 r_i.
	\end{eqnarray*}
	We now define $w_i(x)=(u_i)_{r_i}(x)=\frac{u_i(r_ix)}{r_i^2}, |x|<\frac{\rho}{r_i}$. using the similar arguments of Section \ref{sec5.2}, $\|D^2 w_i\|_{L^\infty(B_{\frac{\rho}{r_i}})}\leq \|D^2 u\|_{L^\infty(B_{\rho})}$ and over a subsequence, $w_i\to w_0$ in $C_\loc^{1,\alpha}(\Rn)$, we claim that $w_0\in \mathcal{Q}^+(u)$. Indeed, by Lemma \ref{lem9.2},
	\begin{eqnarray*}
		\delta(\rho,w_i)=\delta(\rho r_i,u)\leq\sigma^1(\rho r_i)\to 0\ {\rm for}\ \forall \rho>0.
	\end{eqnarray*}
	This means that $w_i$ satisfies $\Delta w_i=1+\lambda v(0)$ in $B_\rho$, except a strip of width $\sigma(\rho r_i)\to0$. Therefore, $w_0$ satisfies $\Delta w_0=1+\lambda v(0)$ a.e. in $\Rn$. The Liouvile's theorem shows $w_0\in \mathcal{Q}^+(u)$. Which is from the following proof,
	recall now the classical Harnack inequality for nonnegative harmonic functions: $\Delta w=0$ and $w\geq0$ in $B_R$, then
	\begin{eqnarray*}
		\displaystyle \sup_{B_{R/2}}w\leq C_n\inf_{B_{R/2}}w
	\end{eqnarray*}
	for some dimensional constant $C_n>0$.
	
	Defining $\displaystyle m_e:=\inf_{\Rn}\partial_{ee}w_0$. It follows by Harnack inequality applied to $w=\partial_{ee} w_0-m_e$ that
	\begin{eqnarray*}
		\displaystyle \sup_{B_{R/2}}(\partial_{ee} w_0-m_e)\leq C_n\inf_{B_{R/2}}(\partial_{ee} w_0-m_e).
	\end{eqnarray*}
	Let $R\to +\infty$ the right hand side tends to $0$, therefore
	\begin{eqnarray*}
		\displaystyle \sup_{B_{R/2}}(\partial_{ee} w_0-m_e)=0.
	\end{eqnarray*}
	It follows that $\partial_{ee} w_0$ is constant.
	
	On the other hand, from the assumptions on $u_i$ with $q=w_0$, we must have either one or the other of following
	\begin{eqnarray*}
		\|u_i-w_0\|_{L^\infty(B_{2r_i})}\geq \epsilon r_i^2,\ \|\nabla u_i-\nabla w_0\|_{L^\infty(B_{2r_i})}\geq \epsilon r_i.
	\end{eqnarray*}
	Since $w_0\in \mathcal{Q}^+(u)$, we have
	\begin{eqnarray*}
		\|w_j-w_0\|_{L^\infty(B_{2})}\geq \epsilon ,\ \|\nabla w_j-\nabla w_0\|_{L^\infty(B_{2})}\geq\epsilon,
	\end{eqnarray*}	which generates  a contradicts with the convergence $w_j\to w_0$ in $C_\loc^{1,\alpha}(\Rn)$. It follows from the proof of $u$ that the conclusion of $v$ are true.\hbx
	
	Now, we define
	\begin{eqnarray*}
		r\to M(r,u,q)=\frac{1}{r^{n+3}}\int_{\partial B_r}(u-q)^2d \mathcal{H}^{n-1},
	\end{eqnarray*}
	\begin{eqnarray*}
		r\to M(r,v,q)=\frac{1}{r^{n+3}}\int_{\partial B_r}(v-q)^2d \mathcal{H}^{n-1}.
	\end{eqnarray*}
Furthermore, we establish   Monneau type monotonicity formulas in Appendix.
	
\begin{pro}\label{cor9.5}
		Let $(u,v)$ be a solution of system \eqref{system} and $x^0\in\Sigma(u,v)$.
		\begin{itemize}
			\item If $x^0\in\Gamma(u)$, as $\lambda\geq0$ or \eqref{gammaulambda} holds. Then the blowup of $u$ at $0$ is unique.
			\item If $x^0\in\Gamma(v)$, as $\lambda\geq0$ or \eqref{gammavlambda} holds. Then the blowup of $v$ at $0$ is unique.
		\end{itemize}
	\end{pro}
	\pf Assume that $x^0\in\Gamma(u)$. Let $\displaystyle q_0=\lim_{r_j\to0^+}u_{x^0,r_j}$ be a blowup over a certain sequence $r_j\to0$. Since $q_0\in \mathcal{Q}^+(u,x^0)$, using Proposition \ref{thm9.4} with $q=q_0$. Based on the convergence $u_{x^0,r_j}\to q_0$ in $C_\loc^{1,\alpha}(\Rn)$, we have, in particular, that
	\begin{eqnarray*}
		M(r_j,u,v,x^0)=\int_{\partial B_1}(u_{x^0,r_j}-q_0)^2d \mathcal{H}^{n-1}\to0.
	\end{eqnarray*}
	The fact that $\displaystyle \lim_{r\to0}M(r,u,v,x^0)$ exists based on Proposition  \ref{thm9.4} and Remark \ref{rem3} imply
	\begin{eqnarray*}
	\int_{\partial B_1}(u_{r}-q)^2d \mathcal{H}^{n-1}\to0,\ \ {\rm as}\ \ r\to0^+.
	\end{eqnarray*}
	Therefore any convergent sequence of $u_{r}$ must converge to $q_0$. The unique about $v$ can be demonstrated using the proof above.  \hbx
	
	Now, we establish the continuous dependence of blowups, which plays an important role in characterizing the properties of singular point set.

	\begin{pro}\label{pro9.6} Let $(u,v)$ be a solution of system \eqref{system} in an open set $\Omega$ in $\Rn$.
		\begin{itemize}
			\item If $x^0\in\Gamma(u)$, as $\lambda\geq0$ or \eqref{gammaulambda} holds. Denoting the blowup of $u$ at $x^0$ by $q_{x^0}$. Then for any $K\Subset \Omega$, there exists a modulus of continuity $\sigma^3(r)$ depending only on $\lambda, n, u, v, \dist(K,\partial \Omega)$, such that
			\begin{eqnarray*}
				|u(x)-q_{x^0}(x)|\leq\sigma^3(|x-x^0|)|x-x^0|^2,\\
				|\nabla u(x)-\nabla q_{x^0}(x)|\leq\sigma^3(|x-x^0|)|x-x^0|
			\end{eqnarray*}
			provided $x^0\in K\cap \Sigma(u), x\in \Omega$. Besides, the mapping $x^0\to q_{x^0}$ is continuous from $\Sigma(u)$ to $\mathcal{Q}^+(u,x^0)$ with
			\begin{eqnarray*}
				\|q_{x^1}-q_{x^2}\|_{L^2{(\partial B_1)}}\leq \sigma^3(|x^1-x^2|)
			\end{eqnarray*}
			provided $x^1,x^2\in\Sigma(u)\cap K$.
			\item If $x^0\in\Gamma(v)$, as $\lambda\geq0$ or \eqref{gammavlambda}  holds. Denoting the blowup of $v$ at $x^0$ by $q_{x^0}$. Then for any $K\Subset \Omega$, there exists a modulus of continuity $\sigma^4(r)$ depending only on $\lambda, n, v, u, \dist(K,\partial \Omega)$, such that
			\begin{eqnarray*}
				|v(x)-q_{x^0}(x)|\leq\sigma^4(|x-x^0|)|x-x^0|^2,\\
				|\nabla v(x)-\nabla q_{x^0}(x)|\leq\sigma^4(|x-x^0|)|x-x^0|,
			\end{eqnarray*}
			provided $x^0\in K\cap \Sigma(v), x\in \Omega$. Besides, the mapping $x^0\to q_{x^0}$ is continuous from $\Sigma(v)$ to $\mathcal{Q}^+(v,x^0)$ with
			\begin{eqnarray*}
				\|q_{x^1}-q_{x^2}\|_{L^2{(\partial B_1)}}\leq \sigma^4(|x^1-x^2|)
			\end{eqnarray*}
			provided $x^1,x^2\in\Sigma(v)\cap K$.
		\end{itemize}

	\end{pro}
	\pf Let $x^0\in\Gamma(u)$. Without loss of generality assume that $D=B_1$ and $K=\overline{B_{1/2}}$, $x^0=0$. On the one hand, for $r>0$, let $q^r\in \mathcal{Q}^+(u,x^0)$ be the approximating polynomial as in Lemma \ref{lem9.3}, and for $\epsilon>0$, let $0<r_\epsilon<\epsilon^2$ such that $\sigma^2(r)\leq \epsilon$ in the same lemma as $r<r_\epsilon$. Then, applying Monneau's monotonicity formula for $u$ with $q=q^{r_\epsilon}$, we will have
	\begin{eqnarray*}
		&&\int_{\partial B_1}(u_{r}-q^{r_\epsilon})^2d \mathcal{H}^{n-1}-\int_0^r\frac{4}{t}F_1(t)ds+\int_0^r\tilde{f}_1(s)ds\\
		&\leq& \int_{\partial B_1}(u_{r_\epsilon}-q^{r_\epsilon})^2d \mathcal{H}^{n-1}-\int_0^{r_\epsilon}\frac{4}{t}F_1(t)ds+\int_0^{r_\epsilon}\tilde{f}_1(s)ds\\
		&\leq& C\epsilon^2+\int_0^{r_\epsilon}|f_1(s)|ds+\int_0^{r_\epsilon}\tilde{f}_1(s)ds\\
		&\leq& C\epsilon^2+C\epsilon^2\\
		&\leq& C\epsilon^2,
	\end{eqnarray*}
	where $F_1, \tilde{f}_1$ are in  Proposition  \ref{thm9.4}. Letting $r\to0^+$, we will obtain
	\begin{eqnarray*}
		\int_{\partial B_1}(q_0-q^{r_\epsilon})^2d \mathcal{H}^{n-1}\leq C\epsilon^2.
	\end{eqnarray*}
	Since both $q_0$ and $q^{r_\epsilon}$ are polynomials from $\mathcal{Q}^+(u,x^0)$, we also have
	\begin{eqnarray*}
		|q_0-q^{r_\epsilon}|\leq C\epsilon,\ \ |\nabla q_0-\nabla q^{r_\epsilon}|\leq C\epsilon\  \ {\rm in}\ B_1.
	\end{eqnarray*}
	Therefore, combining this with Lemma \ref{lem9.3}, we obtain
	\begin{eqnarray*}
		|u-q_0|\leq C\epsilon r_\epsilon^2,\ \ |\nabla u-\nabla q_0|\leq C\epsilon r_\epsilon\ \ {\rm in}\ B_{r_\epsilon}.
	\end{eqnarray*}
	Considering the mapping $r\to\sigma(r)$, which is the inverse of $\epsilon\to r_\epsilon$, the first part of the proposition is established.
	
	On the other hand, we may assume $x^2=0$. Now if $x^1\in B_{r_\epsilon/2}\cap \Sigma(u)$, then as above, we have
	\begin{eqnarray*}
		\int_{\partial B_1}(q_{x^1}-q^{r_\epsilon})^2d \mathcal{H}^{n-1}\leq C \epsilon^2.
	\end{eqnarray*}
	Hence
	\begin{eqnarray*}
		\int_{\partial B_1}(q_{x^1}-q_0)^2d \mathcal{H}^{n-1}\leq C \epsilon^2\ \  {\rm if}\ x^1\in B_{r_\epsilon/2}\cap \Sigma(u).
	\end{eqnarray*}
	This completes the proof. The results of $v$ can be obtained by similar proof.\hbx

	For a singular point $x^0\in\Sigma(u)$, let $q_{x^0}(u)=\frac{1+\lambda v(x^0)}{2}(x\cdot A_{x^0}^ux)$ be the blowup limit of $u$ at $x^0$, and define the dimension of $\Sigma(u)$ at $x^0$ by $d_{x^0}(u)={\rm dim}\ {\rm ker}\ A_{x^0}^u$.
	
	For a singular point $x^0\in\Sigma(v)$, let $q_{x^0}(v)=\frac{1+\lambda u(x^0)}{2}(x\cdot A_{x^0}^vx)$ be the blowup limit of $v$ at $x^0$, and define the dimension of $\Sigma(v)$ at $x^0$ by $d_{x^0}(v)={\rm dim}\ {\rm ker}\ A_{x^0}^v$.
	
Now we give the proof of Theorem \ref{thm9.7}

{\bf Proof of Theorem \ref{thm9.7}.}\ Assume that $x^0\in\Gamma(u)$. Let $K$ be a compact subset of $\Omega$ and $E=\Sigma(u)\cap K$. Then $E$ is also compact using the fact that $\Sigma(u)$ is closed, we claim that the family of polynomials $p_{x^0}=q_{x^0}(x-x^0), x^0\in E$, satisfies the assumptions of Whitney's extension theorem(cf. \cite{1934Whitney}, \cite[Lemma 7.10]{2012PetrosyanShahgholianUraltseva}) on $E$ with $ f(x^0)\equiv0$ and $m=2$. Indeed the condition $(i)$ holds of Whitney's extension theorem, we only need to verify the compatibility conditions $(ii)$ in \cite[Lemma 7.10]{2012PetrosyanShahgholianUraltseva} for $k=0,1,2$.
	
	1)  $k=0$. We need to show that
	\begin{eqnarray*}
		|q_{x^0}(x^1-x^0)|=o(|x^1-x^0|^2),\ {\rm for}\ x^0,x^1\in\Sigma(u)\cap K.
	\end{eqnarray*}
	This is easily verified from the estimate,
	\begin{eqnarray*}
		|u(x)-q_{x^0}(x-x^0)|\leq\sigma(|x-x^0|)|x-x^0|^2
	\end{eqnarray*}
	in Proposition \ref{pro9.6}, by plugging $x=x^1\in \Sigma(u)\cap K$ and noticing that $u(x^1)=0$.
	
	2) $k=1$. Then the compatibility condition is
	\begin{eqnarray*}
		|\nabla q_{x^0}(x^1-x^0)|=o(|x^1-x^0|),\ {\rm for}\ x^0,x^1\in\Sigma(u)\cap K.
	\end{eqnarray*}
	The such case is verified from the second estimate
	\begin{eqnarray*}
		|\nabla u(x)-\nabla q_{x^0}(x-x^0)|\leq\sigma(|x-x^0|)|x-x^0|
	\end{eqnarray*}
	in Proposition \ref{pro9.6}, by plugging $x=x^1\in\Sigma(u)\cap K$ and noticing that $|\nabla u(x^1)|=0$.
	
	3) $k=2$. In terms of coefficient matrix of $q_{x^0}$, we    show that
	\begin{eqnarray*}
		A_{x^0}-A_{x^1}=o(1).
	\end{eqnarray*}
	But this is equivalent to the continuity of the mapping $x^0\to q_{x^0}$, which is established again in Proposition \ref{pro9.6}.
	
	Now let $f$ be a $C^2$ function on $\Rn$, provided by Whitney's extension theorem. To complete the proof, observe that
	\begin{eqnarray*}
		\displaystyle \Sigma(u)\cap K\subset \{|\nabla f|=0\}=\bigcap_{i=1}^n\{\partial_{x_i}f=0\}\subset \bigcap_{i=1}^{n-d}\{\partial_{x_i}f=0\}.
	\end{eqnarray*}
	For $x^0\in\Sigma^d(u)\cap K$, we can arrange the coordinate axes so that the vectors $e_1,\cdots, e_{n-d}$ are the eigenvalues of $D^2 f(x^0)=A_{x^0}$, we then have that
	\begin{eqnarray*}
		\det D^2_{(x_1,\cdots,x_{n-d})}f(x^0)\not=0.
	\end{eqnarray*}
	Since $f$ is $C^2$, the implicit function theorem implies that $\displaystyle \bigcap_{i=1}^{n-d}\{\partial_{x_i}f=0\}$ is a $d$-dimensional $C^1$ manifold in a neighborhood of $x^0$.  The results of $v$ are valid from the similar argument above.  \hbx

\section{Appendix}\label{secAppendix}

\subsection{The properties of the set $\mathcal{K}$}

\setcounter{equation}{0}

\renewcommand\theequation{A.\arabic{equation}}	

Now we give some properties of the set $\mathcal{K}$.
	\begin{pro}\label{pro3.1}
		The set $\mathcal{K}$ is convex, closed and nonempty.
	\end{pro}
	\pf {\bf(i)}. Let $(u_1,v_1),(u_2,v_2)\in \mathcal{K}$, $\forall \theta\in[0,1]$, the desired convexity conclusion $\theta (u_1,v_1)+(1-\theta)(u_2,v_2)\in \mathcal{K}$ is obvious.
	
	{\bf(ii)}. Let $(u_k,v_k)\in \mathcal{K}$ and $(u_k,v_k)\to (u_0,v_0)$ in $\W\times\W$. There hold that
	\begin{description}[itemindent=5pt]
		\item[(a)] $(u_0,v_0)\in \W\times\W$.
		\item[(b)] $(u_{k}-g_1,v_{k}-g_2)\in W^{1,2}_0(\Omega)\times W^{1,2}_0(\Omega)$ and strong convergence implies $(u_0,v_0)|_{\partial\Omega}=(g_1, g_2)$ in the sense of trace.
		\item[(c)] By compact embedding $W^{1,2}(\Omega)\hookrightarrow L^2(\Omega)$, it holds $u,v\geq0$.
	\end{description}
	Therefore, $\mathcal{K}$ is closed.
	
	{\bf(iii).}\ We claim $(g_1^+,g_2^+)\in \mathcal{K}$. Since $(g_1,g_2)\in \W\times\W$, then $(g_1^+,g_2^+)\in \W\times\W$. The facts that $g_1,g_2\geq0$ imply $(g_1^+,g_2^+)|_{\partial\Omega}=(g_1, g_2)|_{\partial\Omega}$, it follows that $(g_1^+-g_1, g_2^+-g_2)\in W^{1,2}_0(\Omega)\times W^{1,2}_0(\Omega)$. \hbx

\subsection{Hausdorff measure of free boundary}

\setcounter{equation}{0}

\renewcommand\theequation{C.\arabic{equation}}

{\bf Proof of Proposition \ref{lemHausdorff}.} We first prove that $\Gamma(u)$ is a  set of finite $(n-1)$-dimensional Hausdorff measure locally in $\Omega$. Let $w_i=\partial_{x_i}u, i=1,2,\cdots,n, E_\epsilon=\{|\nabla u|<\epsilon\}\cap\Omega(u)$.
	
	(i)\ We first consider the case $\lambda\geq0$, observe that
	\begin{eqnarray*}
		\displaystyle 1\leq|\Delta u|^2\leq 2^{n-1}\sum_{i=1}^n|\nabla w_i|^2\ \ \ \ {\rm in} \ \Omega(u).
	\end{eqnarray*}
	Thus for an arbitrary compact set $K\Subset \Omega$, we have
	\begin{eqnarray*}
		\displaystyle |K\cap E_\epsilon|\leq 2^{n-1}\sum_{i=1}^n\int_{K\cap E_\epsilon} |\nabla w_i|^2dx\leq 2^{n-1}\sum_{i=1}^n\int_{K\cap \{|w_i|<\epsilon\}\cap\Omega(u)}|\nabla w_i|^2dx,
	\end{eqnarray*}
	where used the fact $\{|\nabla u|<\epsilon\}\subset \{|\partial_{x_i}u|<\epsilon\}$. Taking $0<\delta <\frac{1}{3}\dist(K,\partial\Omega)$, set $K_\delta:=\{x|\dist(x,K)\le\delta\}$. Then by the proof of Theorem \ref{thmopimalregularity}, since $u\in W^{2,p}_\loc(\Omega)$, then $\nabla w_i^{\pm}=0$ on $\{w_i=0\}\cap K_\delta$. It follows that
	\begin{eqnarray*}
		\int_{K_\delta }\nabla w_i^\pm\nabla \phi dx\leq \lambda\|\nabla v\|_{L^\infty(K_{2\delta})}\int_{K_\delta}\phi dx
	\end{eqnarray*}
	for any nonnegative $\phi\in C_0^\infty(K_\delta)$. And the density of $C_0^\infty(K_\delta)$ in $W^{1,2}_0(K_\delta)$ implies that
	\begin{eqnarray*}
		\int_{K_\delta }\nabla w_i^\pm\nabla \eta dx\leq \lambda\|\nabla v\|_{L^\infty(K_{2\delta})}\int_{K_\delta}\eta dx
	\end{eqnarray*}
	for any nonnegative $\eta\in W^{1,2}_0(K_\delta)$. We now choose $\eta=\psi_\epsilon(w_i^\pm)\phi$ with
	\begin{eqnarray*}
		\psi_\epsilon(t)=\left\{
		\begin{array}{ll}
			0,\ &t\leq0,\\
			\frac{t}{\epsilon},\ & 0\leq t\leq \epsilon,\\
			1,\ & t\geq\epsilon
		\end{array}
		\right.
	\end{eqnarray*}
	and nonnegative $\phi\in C_0^\infty(K_\delta)$ satisfying $\phi|_K=1$. We obtain
	\begin{eqnarray*}
		\lambda\|\nabla v\|_{L^\infty(K_{2\delta})}\int_{K_\delta}\psi_\epsilon(w_i^\pm)\phi dx
		&\geq&\int_{K_\delta}\nabla w_i^\pm\nabla(\psi_\epsilon(w_i^\pm)\phi)dx\\
		&=&\epsilon^{-1}\int_{K_\delta\cap\{0<w_i^\pm<\epsilon\}}|\nabla w_i^\pm|^2\phi dx+\int_{K_\delta}\nabla w_i^\pm\nabla \phi\psi_\epsilon(w_i^\pm)dx,
	\end{eqnarray*}
	which leads to
	\begin{eqnarray*}
		\epsilon^{-1}\int_{K \cap\{0<w_i^\pm<\epsilon\}\cap\Omega(u)}|\nabla w_i|^2 dx&\leq& \epsilon^{-1}\int_{K_\delta\cap\{0<w_i^\pm<\epsilon\}}|\nabla w_i|^2\phi dx\\
		&\leq& \int_{K_\delta}|\nabla w_i||\nabla \phi| dx+\lambda\|\nabla v\|_{L^\infty(K_{2\delta})}\int_{K_\delta}|\phi| dx\\
		&\leq& \|D^2 u\|_{L^\infty(K_\delta)}\int_{K_\delta}|\nabla \phi| dx+\lambda\|\nabla v\|_{L^\infty(K_{2\delta})}\int_{K_\delta}|\phi| dx.
	\end{eqnarray*}
	Thus summing over $i=1,2,\cdots,n$, we arrive at an estimate
	\begin{eqnarray}\label{5.1}
		|K\cap E_\epsilon|\leq \epsilon n2^{n-1} \left(\|D^2 u\|_{L^\infty(K_\delta)}\int_{K_\delta}|\nabla \phi| dx+\lambda\|\nabla v\|_{L^\infty(K_{2\delta})}\int_{K_\delta}|\phi| dx\right).
	\end{eqnarray}
	Using the Besicovitch covering lemma, there exists a covering of $\Gamma(u)\cap K$ by a finite family $\{B^i_\epsilon\}_{i\in I}$ centered on $\Gamma \cap K$ such that no more than $N=N_n$ balls from this family overlap.  For $\epsilon>0$ sufficiently small, we may assume that $B_\epsilon^i\subset K'$ for a slightly larger compact $K'$ such that $K\subset\Int(K')\Subset K_\delta\subset\Omega$. Notice that in every $B_\epsilon^i, |\nabla u|\leq \|D^2 u\|_{L^\infty(K_\delta)}\epsilon$, which shows that $B^i_\epsilon\cap\Omega(u)\subset E_{\kappa}$ with $\kappa=\|D^2 u\|_{L^\infty(K_\delta)}\epsilon$. The facts Lemma \ref{lemdensity} and \eqref{5.1} imply that
	\begin{eqnarray*}
		\displaystyle \sum_{i\in I}|B_\epsilon^i|&\leq& \frac{1}{\beta_1}\sum_{i\in I}|B_\epsilon^i\cap \Omega(u)|\\
		&\leq&\frac{1}{\beta_1}\sum_{i\in I}|B_\epsilon^i\cap E_\kappa|\\
		&\leq& \frac{N}{\beta_1}|K'\cap E_\kappa|\\
		&\leq& \frac{Nn 2^{n-1}\epsilon \|D^2 u\|_{L^\infty(K_\delta)}}{\beta_1}\left(\|D^2 u\|_{L^\infty(K_\delta)}\int_{K_\delta}|\nabla \phi| dx+\lambda\|\nabla v\|_{L^\infty(K_{2\delta})}\int_{K_\delta}|\phi| dx\right).
	\end{eqnarray*}
	Thus
	\begin{eqnarray*}
		\displaystyle \sum_{i\in I}{\rm diam}(B^i_\epsilon)^{n-1}\leq \frac{Nn\|D^2 u\|_{L^\infty(K_\delta)} }{\beta_1 \omega_n}\left(\|D^2 u\|_{L^\infty(K_\delta)}\int_{K_\delta}|\nabla \phi| dx+\lambda\|\nabla v\|_{L^\infty(K_{2\delta})}\int_{K_\delta}|\phi| dx\right),
	\end{eqnarray*}
	letting $\epsilon\to 0$, we conclude that
	\begin{eqnarray*}
		\mathcal{H}^{n-1}(\Gamma(u)\cap K)\leq \frac{Nn\|D^2 u\|_{L^\infty(K_\delta)}}{\beta_1 \omega_n}\left(\|D^2 u\|_{L^\infty(K_\delta)}\int_{K_\delta}|\nabla \phi| dx+\lambda\|\nabla v\|_{L^\infty(K_{2\delta})}\int_{K_\delta}|\phi| dx\right).
	\end{eqnarray*}
Fixing a $\phi\in C_0^\infty(K_\delta)$ and we can find a positive constant $C(n,K,\dist(K,\partial\Omega), \lambda, v, \Omega)$ satisfying
	\begin{eqnarray*}
		\mathcal{H}^{n-1}(\Gamma(u)\cap K)\leq C(n,K,\dist(K,\partial\Omega), \lambda, v, \Omega)
	\end{eqnarray*}
	and complete the proof that $\Gamma(u)$ is a  set of finite $(n-1)$-dimensional Hausdorff measure locally in $\Omega$. Similarly, as $\lambda\geq0$, we can obtain that $\Gamma(v)$ is a  set of finite $(n-1)$-dimensional Hausdorff measure locally in $\Omega$.
	
	(ii)\ The case $\lambda<0, v(z)<-\frac{1}{\lambda}$ for any $z\in \Gamma(u)$, for any $K\Subset \Omega$, there exists  $\varsigma>0$ such that $\sup_{(K\cap \Gamma(u))_\varsigma}v<-\frac{1}{\lambda}$, where $(K\cap \Gamma(u))_\varsigma:=\{x|\dist(x,K\cap \Gamma(u)\}<\varsigma$. Thus there exists $\sigma_0>0$ satisfies
	\begin{eqnarray*}
		\displaystyle \sigma_0\leq|\Delta u|^2\ {\rm in}\ \Omega(u)\cap (K\cap \Gamma(u))_\varsigma.
	\end{eqnarray*}
	Furthermore, using the demonstration above and $\epsilon$ small enough,  we can find a positive constant  \\ $C(n,K,\dist(K,\partial\Omega), \lambda, v,  \Omega)$ satisfying
	\begin{eqnarray*}
		\mathcal{H}^{n-1}(\Gamma(u)\cap K)\leq C(n,K,\dist(K,\partial\Omega), \lambda, v,  \Omega)
	\end{eqnarray*}
	and complete the proof that $\Gamma(u)$ is a  set of finite $(n-1)$-dimensional Hausdorff measure locally in $\Omega$. Similarly, as $\lambda<0$, we can obtain that $\Gamma(v)$ is a  set of finite $(n-1)$-dimensional Hausdorff measure locally in $\Omega$.
	
	Therefore, we have the desired result.\hbx
	\begin{rem}\label{rem1}
		The estimate \eqref{5.1} essentially means that
		\begin{eqnarray*}
			|K\cap \{|\nabla u|<\epsilon\}\cap\Omega(u)|\leq \epsilon n2^{n-1} \left(\|D^2 u\|_{L^\infty(K_\delta)}\int_{K_\delta}|\nabla \phi| dx+\lambda\|\nabla v\|_{L^\infty(K_{2\delta})}\int_{K_\delta}|\phi| dx\right)
		\end{eqnarray*}
		for any $\epsilon>0$. Furthermore, we obtain that $|K\cap \{|\nabla u|=0\}\cap\Omega(u)|=0$.
	\end{rem}

\subsection{Proof of Proposition \ref{prolimit}.}

\setcounter{equation}{0}

\renewcommand\theequation{B.\arabic{equation}}		
In this subsection, we give the proof of  Proposition \ref{prolimit} about sequence of rescaling of solutions.
{\bf Proof of Proposition \ref{prolimit}.}(i) The implications in this part are immediate corollaries of the $C^1_\loc$ convergence.\\
	(ii) Suppose the results of (ii) fail, then there exists $i_k\to\infty$ such that $|u_{r_{i_k}}(y^k)|>0$ for some $y^k\in B_{\delta/2}(x^0)$ and $|u_0|=0$ in $B_\delta(x^0)$.
	
	Now we first prove the case of $u$ with $\lambda\geq0$. Since $y^k\in \Omega(u_{r_{i_k}})$, by the nondegeneracy   at $y^k$ and using the fact $B_{\frac\delta4}(y^k)\subset B_{\frac{3\delta}{4}}(x^0)$, we obtain
	\begin{eqnarray*}
		\displaystyle \sup_{B_{\frac{3\delta}{4}}(x^0)}|u_{r_{i_k}}|\geq\frac{\delta^2}{32n}.
	\end{eqnarray*}
	Then passing to the limit, we arrive that
	\begin{eqnarray*}
		\displaystyle \sup_{B_{\frac{3\delta}{4}}(x^0)}|u_0|\geq\frac{\delta^2}{32n}.
	\end{eqnarray*}
	Contradicting the assumption. Similarly, the results of $v$ hold. For the case $\lambda<0$ with the additional assumption, the nondegeneracy of $u_{r_{i_k}}$ holds when $\delta/t$ small satisfying $B_{\frac{\delta}{t}}(y^k)\subset B_{\frac{3\delta}{4}}(x^0)$.\\
	(iii) Let $x^0\in \Gamma(u)$, without loss of generality, we may assume $u_{\lambda_i}$ is bounded in $W^{2,p}(K), 1<p\leq \infty$ for any $K\Subset \Omega$ and hence $u_0\in W^{2,p}_\loc(\Omega)$. Therefore it is enough to verify that $u_0$ satisfies
	\begin{eqnarray*}
		\Delta u_0=(1+\lambda v(x^0))\chi_{\{u_0>0\}}.
	\end{eqnarray*}
	Since $\nabla u_0=0$ on $\Omega^c(u_0)$ and $\Delta u_0=0$ a.e. there. If $x^0\in \Omega(u_0)$, then from (i) above, we obtain that $B_\delta(x^0)\subset \Omega(u_{r_i})$ for some $\delta>0$ and $i\geq i^0$. Consequently,
	\begin{eqnarray*}
		\Delta u_{r_i}=1+\lambda v(x^0+r_ix)\ \ {\rm in}\ B_\delta(x^0), i\geq i^0.
	\end{eqnarray*}
	For any $\eta\in C_0^\infty(\Omega(u_0))$, then
	\begin{eqnarray*}
		-\int_{\Omega(u_0)}\nabla u_{r_i}\nabla \eta dx=\int_{\Omega(u_0)} (1+\lambda v(x^0+r_ix))\eta dx.
	\end{eqnarray*}
	Passing to the limit we obtain
	\begin{eqnarray*}
		\Delta u_0=(1+\lambda v(x^0))\chi_{\{u_0>0\}}\ {\rm a.e.\ in}\ \Omega.
	\end{eqnarray*}
	For the case $x^0\in\Gamma(v)$, we can get that similarly,
	\begin{eqnarray*}
		\Delta v_0=(1+\lambda u(x^0))\chi_{\{v_0>0\}}\ {\rm a.e.\ in}\ \Omega.
	\end{eqnarray*}
	(iv) If $x^{i_k}\in \Gamma(u_{r_{i_k}})\subset \Omega^c(u_{r_{i_k}})$, the implications in (i) show that $x^0\in \Omega^c(u_0)$, so the question is whether $x^0\in \Gamma(u_0)$ or there exists a ball $B_\delta(x^0)\subset \Omega^c(u_0)$ with $\delta$ small. Now note that in such a ball $|u_0|=0$. Thus, by (ii), we will have that $B_{\delta/2}(x^{i_k})\subset \Omega^c(u_{r_{i_k}})$, contradicting the fact that $x^{i_k}\in \Gamma(u_{r_{i_k}})$. The result for $v$ can be established similarly.\\
	(v)\ Let $x^0\in \Gamma(u)$. In order to prove that $u_{r_{i_k}}\to u_0$ strongly in $W^{2,p}_\loc(\Omega)$ for any $1<p<\infty$. It will suffice to show that
	\begin{eqnarray}
		D^2 u_{r_{i_k}}\to D^2u_0\ {\rm a.e.\ in}\ \Omega.\label{5.3}
	\end{eqnarray}
	Since $D^2 u_{r_{i_k}}$ are uniformly bounded in any $K\Subset\Omega$. Now  the subsequence pointwise convergence in $\Omega(u_0)$ follows from the subsequent proof. For any $x^*\in \Omega(u_0)$  and $B_\delta(x^*)\subset \Omega(u_0)$. The part (i) implies that
	\begin{eqnarray*}
		\left\{
		\begin{array}{ll}
			\Delta u_{r_{i_k}}=1+\lambda v(x^0+r_{i_k}x),\ {\rm in}\ B_\delta(x^*),\\
			\Delta u_0=1+\lambda v(x^0),\ {\rm in}\ B_\delta(x^*).
		\end{array}
		\right.
	\end{eqnarray*}
	It follows that $\Delta(u_{r_{i_k}}-u_0)=\lambda v(x^0+r_{i_k}x)-\lambda v(x^0)$ in $B_\delta(x^*)$. Applying  \cite[Theorem 1.1]{2012PetrosyanShahgholianUraltseva}, then
	\begin{eqnarray*}
		\|u_{r_{i_k}}-u_0\|_{W^{2,p}(\overline{B_{\delta/2}(x^*)})} \leq C\|u_{r_{i_k}}-u_0\|_{L^1(B_\delta(x^*))}+C\|v(x^0+r_{i_k}x)-v(x^0)\|_{L^p(B_\delta(x^*))}.
	\end{eqnarray*}
	Thus there exists a subsequence satisfying $D^2 u_{r_{i_k}}\to D^2u_0$ a.e. in $\Omega(u_0)$. The pointwise convergence on $\Int(\Lambda(u_0))$ follows from part (ii). Using  the fact from Corollary \ref{corGamma}, we can verify the result \eqref{5.3}. Combining with \eqref{5.2}, the part (v) holds for $u$ and the conclusion for the case $x^0\in \Gamma(v)$ is similar.\\
	(vi)\ It is obvious that
	\begin{eqnarray*}
		&& \int_\Omega|\chi_{\{u_{r_i}>0\}}-\chi_{\{u_0>0\}}|dx\\
		&&=\int_{\Omega\cap\{u_{r_i}>0\}\cap\{u_0>0\} }|\chi_{\{u_{r_i}>0\}}-\chi_{\{u_0>0\}}|dx+\int_{\Omega\cap\{u_{r_i}=0\}\cap\{u_0>0\} }|\chi_{\{u_{r_i}>0\}}-\chi_{\{u_0>0\}}|dx\\
		&&+\int_{\Omega\cap\{u_{r_i}>0\}\cap\{u_0=0\} }|\chi_{\{u_{r_i}>0\}}-\chi_{\{u_0>0\}}|dx+\int_{\Omega\cap\{u_{r_i}=0\}\cap\{u_0=0\} }|\chi_{\{u_{r_i}>0\}}-\chi_{\{u_0>0\}}|dx\\
		&&=\int_{\Omega\cap\{u_{r_i}=0\}\cap\{u_0>0\} }dx+\int_{\Omega\cap\{u_{r_i}>0\}\cap\{u_0=0\} }dx.
	\end{eqnarray*}
	To obtain the desired results, we claim that
	\begin{eqnarray}
		\int_{\Omega\cap\{u_{r_i}=0\}\cap\{u_0>0\} }dx\to 0\  \ {\rm and}\ \
		\int_{\Omega\cap\{u_{r_i}>0\}\cap\{u_0=0\} }dx \to 0.\label{5.4}
	\end{eqnarray}
	If $x\in \{u_0>0\}$, by part(i), there exists $B_\delta(x)$ such that $u_{r_i}(x)>0$ for $x\in B_\delta(x)$ as $i\geq i^0$. Thus $|\{u_{r_i}=0\}\cap\{u_0>0\}|=0$ for large $i$. \\
	If $x\in \Int(\Lambda(u_0))$. Then the $C^1$ convergence implies that there exists $B_\delta(x)\subset \Int(\Lambda(u_0))$ such that $u_{r_i}(x)=0$ on $B_\delta(x)$ for large $i$. Thus $|\{u_0=0\}\cap\{u_{r_i}>0\}|=0$ for $i\ge i^0$ with $i^0$ large enough.
	Therefore, we can get that the claim \eqref{5.4} holds.\hbx

\subsection{Weiss-type monotonicity formula}

\setcounter{equation}{0}

\renewcommand\theequation{D.\arabic{equation}}

According to the free boundary theory of classical obstacle problem and Proposition \ref{thmclassification}, we know that the free boundary of system \eqref{system} can be divided into regular and singular points. Now, we give the Weiss type monotonicity formula, which plays an role on that showing blowup limit is homogeneous of degree two(see also Proposition \ref{thmclassification}) and  Subsection \ref{subsec5.5}.

\begin{pro}[Weiss-type monotonicity formula]\label{weissmonotonocity}
		Let $(u,v)$ be a solution of system \eqref{system}, then if  $x^0\in \Gamma(u)$,
		\begin{eqnarray*}
			\frac{d}{dr}\left(W(r,u,v,x^0)-F_1(r)\right)=\frac{1}{r}\int_{\partial B_1}\left(x\cdot\nabla u_r(x)-2u_r(x)\right)^2d\mathcal{H}^{n-1},
		\end{eqnarray*}
		where
		\begin{eqnarray*}
			W(r,u,v,x^0):=\frac{1}{r^{n+2}}\int_{B_r(x^0)}\frac{1}{2}|\nabla u|^2+u+\lambda vudx-\frac{1}{r^{n+3}}\int_{\partial B_r(x^0)}u^2d \mathcal{H}^{n-1},
		\end{eqnarray*}
		\begin{eqnarray*}
			F_1(r)=\int_0^rf_1(s)ds,\ \  f_1(s)=\lambda\int_{B_1}\frac{u(sx+x^0)}{s^2}\nabla v(sx+x^0)\cdot xdx.
		\end{eqnarray*}
		If  $x^0\in \Gamma(v)$, then
		\begin{eqnarray*}
			\frac{d}{dr}\left(W(r,v,u,x^0)-F_2(r)\right)=\frac{1}{r}\int_{\partial B_1}\left(x\cdot\nabla v_r(x)-2v_r(x)\right)^2d\mathcal{H}^{n-1},
		\end{eqnarray*}
		where
		\begin{eqnarray*}
			F_2(r)=\int_0^rf(s)ds,\ \  f_2(s)=\lambda\int_{B_1}\frac{v(sx+x^0)}{s^2}\nabla u(sx+x^0)\cdot xdx.
		\end{eqnarray*}
		Furthermore, the blowup limit $u_0$ of $u$ or the blowup limit $v_0$ of $v$ is homogeneous of degree two, that is
		\begin{eqnarray*}
			u_0(rx)=r^2u_0(x),\ x\in \Rn, r>0\ {\rm if}\ x_0\in\Gamma(u),\\
			v_0(rx)=r^2v_0(x),\ x\in \Rn, r>0\ {\rm if}\ x_0\in\Gamma(v).
		\end{eqnarray*}
	\end{pro}
	\pf  If $x^0\in\Gamma(u)$, it is obvious that $W(r,u,v(x),x^0)=W(1,u_r,v(rx+x^0),0)$, and
	\begin{eqnarray*}
		&&\frac{d}{dr}W(r,u,v,x^0)\\
		&=&\int_{B_1}-\Delta u_r(x)\frac{du_r(x)}{dr}+\frac{du_r(x)}{dr}+\lambda \frac{du_r(x)}{dr}v(rx+x^0)+\lambda u_r(x)\frac{d(v(rx+x^0))}{dr}dx\\
		&&+\int_{\partial B_1}\left(\frac{\partial u_r(x)}{\partial  \nu}\frac{du_r(x)}{dr}-2u_r(x)\frac{du_r(x)}{dr}\right)d\mathcal{H}^{n-1}\\
		&=&\int_{B_1}\lambda u_r(x)\frac{d(v(rx+x^0))}{dr}dx+\int_{\partial B_1}\left(\frac{\partial u_r(x)}{\partial  \nu}\frac{du_r(x)}{dr}-2u_r(x)\frac{du_r(x)}{dr}\right)d\mathcal{H}^{n-1}\\
		&=&\lambda \int_{B_1}\frac{u(rx+x^0)}{r^2}x\cdot \nabla v(rx+x^0)dx+\int_{\partial B_1}\left(x\cdot\nabla u_r(x)-2u_r(x)\right)^2d\mathcal{H}^{n-1},
	\end{eqnarray*}
	where using the facts
	\begin{eqnarray*}
		\Delta u_r=(1+\lambda v(rx+x^0)\ \ \  {\rm in}\ \{u_r(x)>0\},\\
		\frac{d u_r(x)}{dr}=0\ \ \   {\rm in}\ \{ u_r(x)=0\},\\
		\frac{du_r(x)}{dr}=\frac{1}{r}\left(x\cdot\nabla u_r(x)-2u_r(x)\right).
	\end{eqnarray*}
	Now we define functions $F_1(r)$ and $f_1(s)$ as in Theorem \ref{weissmonotonocity}.	Then
	\begin{eqnarray*}
		\frac{d}{dr}\left(W(r,u,v,x^0)-F_1(r)\right)=\frac{1}{r}\int_{\partial B_1}\left(x\cdot\nabla u_r(x)-2u_r(x)\right)^2d\mathcal{H}^{n-1}.
	\end{eqnarray*}
	Since $\displaystyle \lim_{r\to0^+}F_1(r)=0$, we get that $\displaystyle \lim_{r\to0^+}W(r,u,v,x^0)$ exists. By the regularity of $u$ and $x^0\in\Gamma(u)$, we get that $W(0^+,u,v,x^0)>-\infty$. Since
	\begin{eqnarray*}
		W(r_k\rho,u,v,x^0)&=&W(\rho,u_{r_k},v(r_kx+x^0),0)\\
		&=&\frac{1}{\rho^{n+2}}\int_{B_\rho}\frac{1}{2}|\nabla u_{r_k}|^2+u_{r_k}+\lambda v(r_k x+x^0)u_{r_k}dx-\int_{\partial B_\rho}u_{r_k}^2d\mathcal{H}^{n-1}.
	\end{eqnarray*}
	Then
	\begin{eqnarray*}
		\displaystyle W(0^+,u,v,x^0)&=&\lim_{k\to\infty}W(r_k\rho,u,v,x^0)\\
		&=&\lim_{k\to\infty}W(\rho,u_{r_k},v(r_kx+x^0),0)\\
		&=&W(\rho,u_0,v(x^0),0).
	\end{eqnarray*}
	We obtain that $W(\rho,u_0,v(x^0),0)$ is constant in $\rho$. Since $\Delta u_0=(1+\lambda v(x^0))\chi_{\{u_0>0\}}$ in $B_\rho$, using the proof above replacing $v(x)$ by $v(x^0)$, then $x\cdot\nabla u_0-2u_0\equiv0$ in $\Rn$, and therefore $u_0$ is homogeneous of degree $2$. The case of $v$ is similar.\hbx

	We define
	\begin{eqnarray*}
		\left.
		\begin{array}{ll}
			\displaystyle \omega(x^0,u_0):=W(0^+,u,v,x^0)=\lim_{r\to0^+}W(r,u,v,x^0),\ {\rm if}\ x^0\in\Gamma(u),\\
			\displaystyle \omega(x^0,v_0):=W(0^+,v,u,x^0)=\lim_{r\to0^+}W(r,v,u,x^0),\ {\rm if}\ x^0\in\Gamma(v),
		\end{array}
		\right.
	\end{eqnarray*}
	which exist by Theorem \ref{weissmonotonocity}, are called the balanced energy of $u$ and $v$ at $x^0$. More generality, if $(u,v)$ is a solution of system \eqref{system}, $\omega(x^0,u_0)$ defines a function $\omega:\Gamma(u)\to \mathbb{R}$ and $\omega(x^0,v_0)$ defines a function $\omega:\Gamma(v)\to \mathbb{R}$.
	
	By the proof above, if $x^0\in\Gamma(u)$,
	\begin{eqnarray*}
		\displaystyle \omega(x^0,u_0)=\lim_{\rho\to0^+}W(r\rho,u,v,x^0)&=&\lim_{\rho\to0^+}W(r,u_\rho,v(\rho x+x^0),0)\\
		&=&W(r,u_0,v(x^0),0)\\
		&=&W(1,u_0,v(x^0),0).
	\end{eqnarray*}
	
	If $x^0\in\Gamma(v)$,
	\begin{eqnarray*}
		\displaystyle \omega(x^0,v_0)=\lim_{\rho\to0^+}W(r\rho,v,u,x^0)&=&\lim_{\rho\to0^+}W(r,v_\rho,u(\rho x+x^0),0)\\
		&=&W(r,v_0,u(x^0),0)\\
		&=&W(1,v_0,u(x^0),0).
	\end{eqnarray*}

The proof is over. \hbx

\subsection*{Monneau-type monotonicity formulas}

\setcounter{equation}{0}

\renewcommand\theequation{E.\arabic{equation}}
Next, we establish   Monneau type monotonicity formulas.
	\begin{pro}\label{thm9.4}
		Let $(u,v)$ be a solution of system \eqref{system}.
		\begin{itemize}
			\item If $0\in\Sigma(u)$. As $\lambda\geq0$ or \eqref{gammaulambda} holds.  Then for any $q\in \mathcal{Q}^+(u)$, the functional
			\begin{eqnarray*}
				r\to M(r,u,q)+\tilde{F}_1(r)-\bar{F}_1(r)
			\end{eqnarray*}
			is monotone nondecreasing for $r\in(0,1)$, where
			\begin{eqnarray*}
				&\displaystyle \tilde{F}_1(r)=\int_0^t\tilde{f}_1(s)ds,\ \  \tilde{f}_1(s):=\frac{2\lambda}{s^{n+3}}\int_{B_s}\left(v-v(0))(u+q(x)\right)\chi_{\{u>0\}}dx,\\
				&\displaystyle \bar{F}_1(r)=4\int_0^r\frac{F_1(t)}{t}dt
			\end{eqnarray*}
			and $F_1(t)$ as in Proposition  {\rm\ref{weissmonotonocity}}.
			
			\item If $0\in\Sigma(v)$. As $\lambda\geq0$ or \eqref{gammavlambda} holds. Then for any $q\in \mathcal{Q}^+(v)$, the functional
			\begin{eqnarray*}
				r\to M(r,v,q)+\tilde{F}_2(r)-\bar{F}_2(r)
			\end{eqnarray*}
			is monotone nondecreasing for $r\in(0,1)$, where
			\begin{eqnarray*}
				&\displaystyle \tilde{F}_2(r)=\int_0^t\tilde{f}_2(s)ds,\ \  \tilde{f}_2(s):=\frac{2\lambda}{s^{n+3}}\int_{B_s}\left(v-v(0))(v+q(x)\right)\chi_{\{v>0\}}dx,\\
				&\displaystyle \bar{F}_2(r)=4\int_0^r\frac{F_2(t)}{t}dt
			\end{eqnarray*}
			and $F_2(t)$ as in Proposition  {\rm \ref{weissmonotonocity}}.
		\end{itemize}
		
	\end{pro}

	\begin{rem}\label{rem3}
For the scaling version, we can introduce Monneau's functional centered at any $x^0\in\Sigma(u)$ or $x^0\in\Sigma(u)$ by
		\begin{eqnarray*}
			M(r,u,q,x^0)=\frac{1}{r^{n+3}}\int_{\partial B_r(x^0)}[u(x)-q(x-x^0)]^2d\mathcal{H}^{n-1},\\
			M(r,v,q,x^0)=\frac{1}{r^{n+3}}\int_{\partial B_r(x^0)}[v(x)-q(x-x^0)]^2d\mathcal{H}^{n-1}.
		\end{eqnarray*}
		The results of Proposition  {\rm\ref{thm9.4}} becomes
		\begin{eqnarray*}
			r\to M(r,u,q,x^0)+\tilde{F}_1(r)-\bar{F}_1(r)
		\end{eqnarray*}
		is monotone nondecreasing for $r\in(0,1)$  with
		\begin{eqnarray*}
			&\displaystyle \tilde{F}_1(r)=\int_0^t\tilde{f}_1(s)ds,\ \  \tilde{f}_1(s):=\frac{2\lambda}{s^{n+3}}\int_{B_s}\left(v-v(x^0))(u+q(x)\right)\chi_{\{u>0\}}dx.
		\end{eqnarray*}
		And
		\begin{eqnarray*}
			r\to M(r,v,q)+\tilde{F}_2(r)-\bar{F}_2(r)
		\end{eqnarray*}
		is monotone nondecreasing for $r\in(0,1)$, where
		\begin{eqnarray*}
			&\displaystyle \tilde{F}_2(r)=\int_0^t\tilde{f}_2(s)ds,\ \  \tilde{f}_2(s):=\frac{2\lambda}{s^{n+3}}\int_{B_s}\left(v-v(x^0))(v+q(x)\right)\chi_{\{v>0\}}dx.
		\end{eqnarray*}
	\end{rem}
{\bf Proof of Proposition  \ref{thm9.4}.} Since
	\begin{eqnarray*}
		W(r,u,v,0):&=&\frac{1}{r^{n+2}}\int_{B_r}\frac{1}{2}|\nabla u|^2+u+\lambda vudx-\frac{1}{r^{n+3}}\int_{\partial B_r}u^2d \mathcal{H}^{n-1}\\
		&=&\int_{B_1}(\frac{1}{2}|\nabla u_{r}|^2+u_{r}+\lambda v(ry)u_{r,x^0}dx-\int_{\partial B_1}u_{r}^2d \mathcal{H}^{n-1},
	\end{eqnarray*}
	thus
	\begin{eqnarray*}
		W(0^+,u,v,0)&=&\int_{B_1}\frac{1}{2}|\nabla q|^2+q+\lambda v(0)qdx-\int_{\partial B_1}q^2d \mathcal{H}^{n-1}\\
		&=&\frac{1}{r^{n+2}}\int_{B_r}\frac{1}{2}|\nabla q|^2+q+\lambda v(0)qdx-\frac{1}{r^{n+3}}\int_{\partial B_r}q^2d \mathcal{H}^{n-1}\\
		&=& W(r,q,v(0),0).
	\end{eqnarray*}
	Furthermore, let $w=u-q(x)$, since
	\begin{eqnarray*}
		&&\frac{1}{r^{n+2}}\int_{B_r}\nabla w\nabla q(x)dx\\
		&=&\frac{1}{r^{n+2}}\int_{B_r}-\Delta q(x)\cdot wdx+\frac{1}{r^{n+2}}\int_{\partial B_r}\nabla q(x)\cdot\frac{x}{r}wd \mathcal{H}^{n-1}\\
		&=&\frac{1}{r^{n+2}}\int_{B_r}-(1+v(0))\cdot wdx+\frac{1}{r^{n+2}}\int_{\partial B_r}\nabla q(x)\cdot\frac{x}{r}wd \mathcal{H}^{n-1},
	\end{eqnarray*}
	we obtain combining with $q\in \mathcal{Q}^+$ and $W(0^+,u,v,0)=W(r,q,v(0),0)$ that
	\begin{eqnarray*}
		&& W(r,u,v,0)-W(0^+,u,v,0)\\
		&=& W(r,u,v,0)-W(r,q,v(0),0)\\
		&=&\frac{1}{r^{n+2}}\int_{B_r}\frac{1}{2}|\nabla w|^2+\nabla w \nabla q(x)+w+\lambda vu-\lambda v(0)q(x)dx\\
		&&-\frac{1}{r^{n+3}}\int_{\partial B_r}(w^2+2wq(x))d \mathcal{H}^{n-1}\\
		&=&\frac{1}{r^{n+2}}\int_{B_r}\frac{1}{2}|\nabla w|^2+w+\lambda vu-\lambda v(0)q(x)dx\\
		&&-\frac{1}{r^{n+2}}\int_{B_r}(1+\lambda v(0))wdx+\frac{1}{r^{n+3}}\int_{\partial B_r}x\cdot\nabla q(x)wd \mathcal{H}^{n-1}\\
		&&-\frac{1}{r^{n+3}}\int_{\partial B_r}(w^2+2wq(x))d \mathcal{H}^{n-1}\\
		&=&\frac{1}{r^{n+2}}\int_{B_r}\frac{1}{2}|\nabla w|^2-\lambda v(0)w+\lambda vu-\lambda v(0)q(x)dx-\frac{1}{r^{n+3}}\int_{\partial B_r}w^2d \mathcal{H}^{n-1}\\
		&&+\frac{1}{r^{n+3}}\int_{\partial B_r}\left(x\cdot\nabla q(x)-2q(x)\right)wd \mathcal{H}^{n-1}\\
		&=&\frac{1}{r^{n+2}}\int_{B_r}\frac{1}{2}|\nabla w|^2-\lambda v(0)w+\lambda vu-\lambda v(0)q(x)dx-\frac{1}{r^{n+3}}\int_{\partial B_r}w^2d \mathcal{H}^{n-1}.
	\end{eqnarray*}
	Since
	\begin{eqnarray*}
		\int_{B_r}\frac{1}{2}|\nabla w|^2dx=\int_{B_r}\frac{1}{2}(-\Delta w)wdx+\frac{1}{2}\int_{\partial B_r}\nabla w\cdot\frac{x}{r}wd\mathcal{H}^{n-1},
	\end{eqnarray*}
	we get that
	\begin{eqnarray*}
		W(r,u,v,0)-W(0^+,u,v,0)&=&\frac{1}{r^{n+2}}\int_{B_r}\frac{1}{2}(-\Delta w)w+\lambda vu-\lambda v(0)q(x)-\lambda v(0)wdx\\
		&&+\frac{1}{2r^{n+3}}\int_{\partial B_r}[x\cdot\nabla w-2w]wd\mathcal{H}^{n-1}.
	\end{eqnarray*}
	In fact,
	\begin{eqnarray*}
		-\frac{1}{2}\Delta w w+\lambda vu-\lambda v(0)q(x)-\lambda v(0)(u-q(x))=-\frac{1}{2}\Delta(u-q(x))(u-q(x))+\lambda(v-v(0))u.
	\end{eqnarray*}
	On the set $\{u(x)=0\}$, one has
	\begin{eqnarray}\label{9.1}
		\left.
		\begin{array}{ll}
			-\frac{1}{2}\Delta(u-q(x))(u-q(x))+\lambda(v-v(0))u&=-\frac{1}{2}\Delta q(x) q(x)\\[0.5em]
			&=-\frac{1}{2}(1+\lambda v(0))q(x)\\[0.5em]
			&=-\frac{1}{2}(1+\lambda v(0))q(x)\\[0.5em]
			&\leq0.
		\end{array}
		\right.
	\end{eqnarray}
	On the set $\{u(x)>0\}$, we have
	\begin{eqnarray}\label{9.2}
		\left.
		\begin{array}{ll}
			&-\frac{1}{2}\Delta(u-q(x))(u-q(x))+\lambda(v-v(0))u\\ [0.5em]
			&=-\frac{\lambda}{2}(v-v(0))[u-q(x)]+\lambda(v-v(0))u\\[0.5em]
			&=\frac{\lambda}{2}(v-v(0))(u+q(x)).
		\end{array}
		\right.
	\end{eqnarray}
	It is obvious that
	\begin{eqnarray*}
		\frac{dM(r,u,q)}{dr}=\frac{2}{r^{n+4}}\int_{\partial B_r}(w\nabla w\cdot x-2w)d \mathcal{H}^{n-1},
	\end{eqnarray*}
	based on \eqref{9.1} and \eqref{9.2},
	\begin{eqnarray*}
		&&\frac{dM(r,u,q,0)}{dr}\\
		&=&\frac{4}{r}\Big\{W(r,u,v,0)-W(0^+,u,v,0) \\
		&&\ \left.-\frac{1}{r^{n+2}}\int_{B_r}\frac{1}{2}(-\Delta w)w+\lambda vw-\lambda v(0)q(x)-\lambda v(0)wdx\right\}\\
		&\geq&\frac{4}{r}\Big\{W(r,u,v,0)-W(0^+,u,v,0) \\
		&&\left.-\frac{\lambda}{2r^{n+2}}\int_{B_r}\left((v-v(0))(u+q(x))
		\right)\chi_{\{u>0\}}dx\right\}.
	\end{eqnarray*}
	We denote that
	\begin{eqnarray*}
		\tilde{f}_1(s)&=&\frac{2\lambda}{s^{n+3}}\int_{B_s}\left(vu+vq(x)-v(0)u-v(0)q(x)
		\right)\chi_{\{u>0\}}dx\\
		&=&\frac{2\lambda}{s}\int_{B_1}\left(v(sx)\frac{u(sx)}{s^2}+v(sx)q(x)-v(0)\frac{u(sx)}{s^2}-v(x^0)q(x)
		\right)\chi_{\{u(sx)>0\}}dx.
	\end{eqnarray*}
	The limit $\displaystyle f(0^+)$ exists. Since $v\in C^{1,1}_\loc (B_1)$, taking $s$ small such that
	\begin{eqnarray*}
		&&\left|\frac{1}{s}\int_{B_1}\left((v(sx)-v(0))q(x)+[v(sx)-v(0)]\frac{u(sx)}{s^2}\right)\chi_{\{u(sx)>0\}}dx\right|\\
		&\leq&C(v)\int_{B_1}q(x)dx+\frac{u(sx)}{s^2}dx\\
		&<&+\infty.
	\end{eqnarray*}
	Thus $\tilde{f}_1(s)\in L^1((0,r))$ for $t$ small. Set $\tilde{F}_1(r)=\int_0^t\tilde{f}_1(s)ds$, then $\displaystyle \lim_{r\to0^+}\tilde{F}_1(r)=0$. Recalling Proposition \ref{weissmonotonocity},
	\begin{eqnarray*}
		F_1(r)=\int_0^rf_1(s)ds
	\end{eqnarray*}
	and
	\begin{eqnarray*}
		f_1(s)=\lambda\int_{B_1}u_s(x)\nabla v(sx)\cdot xdx.
	\end{eqnarray*}
	We denote
	\begin{eqnarray*}
		\bar{F}_1(r)=4\int_0^r\frac{F_1(t)}{t}dt.
	\end{eqnarray*}
	Then $\displaystyle \lim_{r\to0^+}\bar{F}(r)=0$.
	Thus
	\begin{eqnarray*}
		&&\frac{d(M(r,u,q,x^0)+\tilde{F}_1(r)-\bar{F}_1(r))}{dr}\geq \frac{4}{r}\left\{W(r,u,v,x^0)-F_1(r)-W(0^+,u,v,x^0)\right\}\geq0.
	\end{eqnarray*}
	The results of $v$ are similar.\hbx
	
	\section*{Acknowledgment} \ This work  is supported by the National Natural Science Foundation of China grant 11971331, 12125102, 12301258, Sichuan Youth Science and Technology Foundation 2021JDTD0024, and Sichuan Science and Technology Program(2023NSFSC1298).


\vspace{0.3cm}
	\begin{thebibliography}{99}
		\footnotesize
		
		\bibitem{1992Adams} D. R. Adams, Weakly elliptic systems with obstacle constraints: part I.-A $2\times2$ model problem. In: \emph{IMA Vol. Math. Appl.,} vol. 42, \emph{Springer, New York,} 1992, pp. 1-14.
		
		\bibitem{2000Adams} D. R. Adams, Weakly elliptic systems with obstacle constraints:part. II.-An $N\times N$ model problem. \emph{J. Geom. Anal.,} \textbf{10}(2000) 375-412.
		
		\bibitem{2022AFSW} G. Aleksanyan, M. Fotouhi, H. Shahgholian, G. S.  Weiss,  Regularity of the free boundary for a parabolic cooperative system. \emph{Calc. Var. Partial Differential Equations,}  \textbf{61}(2022),  Paper No. 124, 38 pp.
		
		\bibitem{2019AllenShaholian} M. Allen, H. Shahgholian, A new boundary Harnack principle(equations with right hand side). \emph{Arch. Ration. Mech. Anal.,} \textbf{234}(2019) 1413--1444.
		
		\bibitem{1984AltCaffarelliFriedman} H. W. Alt, L. A. Caffarelli,  A. Friedman, Variational problems with two phases and their free boundaries. \emph{Trans. Amer. Math. Soc.,} \textbf{282}(1984) 431--461.
		
		\bibitem{2015AnderssonShahgholianUraltsevaWeiss} J. Andersson, H. Shahgholian, N. N. Uraltseva, G. S. Weiss,  Equilibrium points of a singular cooperative system with free boundary. \emph{Adv. Math.,} \textbf{280}(2015), 743-771.
		
\bibitem{1977Caffarelli} L. A. Caffarelli, The regularity of free boundaries in higher dimensions. \emph{Acta Math.,} \textbf{139}(1977), 155-184.
		
\bibitem{1980Caffarelli} L. A. Caffarelli, Compactness methods in free boundary problems. \emph{Comm. Partial Differential Equations,} \textbf{5}(1980), 427-448.
		
\bibitem{1998Caffarelli} L. A. Caffarelli, The obstacle problem revisited. \emph{J. Fourier Anal. Appl.,} \textbf{4}(1998), 383-402.
		
		\bibitem{2018CaffarelliShahgholianYeressian} L. A. Caffarelli, H. Shahgholian, K. Yeressian, A minimization problem with free boundary related to a cooperative system. \emph{Duke Math. J.,} \textbf{167}(2018), 1825-1882.
		
		
		\bibitem{1985ChipotVergara-Caffarelli}  M. Chipot, G. Vergara-Caffarelli, The N-membranes problem. \emph{Appl. Math. Optim.,} \textbf{13}(1985) 231-249.
		
		\bibitem{2004DKNF} B. Deconinck, P. G. Kevrekidis, H. E. Nistazakis, D. J. Frantzeskakis, Linearly coupled Bose-Einstein condensates: From rabi oscillations and quasiperiodic solutions to oscillating domain walls and spiral waves, \emph{Phys. Rev. A,} \textbf{70}(2004) 063605.
		
		
		\bibitem{1986DuzaarFuchs} F. Duzaar, M. Fuchs, Optimal regularity theorems for variational problems with obstacles. \emph{Manuscripta Math.,} \textbf{56}(1986) 209-234.
		
		
		\bibitem{21ElHajjShahgholian} L. El Hajj, H. Shahgholian, Remarks on the convexity of free boundaries(scalar and system cases). \emph{Algebra i Analiz,} \textbf{32}(2020), 146-165; \emph{St. Petersburg Math. J.,} \textbf{32}(2021), 713-727.
		
		\bibitem{1979EvansFriedman}  L. C. Evans, A. Friedman, Optimal stochastic switching and the Dirichlet problem for the Bellman equation. \emph{Trans. Amer. Math. Soc.,} \textbf{253}(1979) 365-389.
		
		\bibitem{2022FR} X. Fern\'{a}ndez-Real, X. Ros-Oton, Regularity theory for elliptic PDE. \emph{Zurich Lectures in Advanced Mathematics, 28. EMS Press, Berlin,} 2022.
		
		\bibitem{2018Figalli} A. Figalli, Free boundary regularity in obstacle problems, \emph{Journ\'{e}es \'{e}quations aux d\'{e}riv\'{e}es partielles,} (2018),  1-26.
		
		\bibitem{2019FigalliSerra} A. Figalli, J. Serra, On the fine structure of the free boundary for the classical obstacle problem. \emph{Invent. Math.,} \textbf{215}(2019), 311-366.
		
		\bibitem{2015FigalliShahgholian} A. Figalli, H. Shahgholian, An overview of unconstrained free boundary problems. \emph{Philos. Trans. Roy. Soc. A,} \textbf{373}(2015), 20140281, 11 pp.
		
		\bibitem{1987Fuchs} M. Fuchs, An elementary partial regularity proof for vector-valued obstacle problems. \emph{Math. Ann.,} \textbf{279}(1987) 217-226.
		
\bibitem{1972Frehse} J. Frehse, On the regularity of the solution of a second order variational inequality. \emph{Boll. Un. Mat. Ital.,} \textbf{6}(1972), 312-315.
		
		
		\bibitem{1984FriedmanPhillips} A. Friedman, D. Phillips, The free boundary of a semilinear elliptic equation. \emph{Trans. Amer. Math. Soc.,} \textbf{282}(1984), 153-182.
		
		
		\bibitem{2001GilbargTrudinger} D. Gilbarg, N. S. Trudinger, Elliptic partial differential equations of second order. Reprint of the 1998 edition.  \emph{Springer-Verlag, Berlin,} 2001.
		
		\bibitem{1971Lions} J. L. Lions, Optimal Control of Systems Governed by Partial Differential Equations.  \emph{Springer-Verlag, New York, Berlin,} 1971, translated from French by S.K. Mitter.
		
		\bibitem{1985Kawohl} B. Kawohl, When are solutions to nonlinear elliptic boundary value problems convex? \emph{Comm. Partial Differential Equations,} \textbf{10}(1985),  1213-1225.
		
		\bibitem{1977KinderlehrerNirenberg} D. Kinderlehrer, L. Nirenberg, Regularity in free boundary problems. \emph{Ann. Scuola Norm. Sup. Pisa Cl. Sci.,} \textbf{4}(1977), 373-391.
		
		\bibitem{2020MazzoleniTerraciniVelichkov} D. Mazzoleni, S. Terracini, B. Velichkov, Regularity of the free boundary for the vectorial Bernoulli problem. \emph{Anal. PDE,} \textbf{13}(2020), 741-764.
		
		\bibitem{2017MazzoleniTerraciniVelichkov} D. Mazzoleni, S. Terracini, B. Velichkov, Regularity of the optimal sets for some spectral functionals. \emph{Geom. Funct. Anal.,} \textbf{27}(2017), 373-426.
		
		
		\bibitem{1996MCSS} M. Mitchell, Z. Chen, M. Shih, M. Segev, Self-trapping of partially spatially incoherent light, \emph{Phys. Rev. Lett.,} \textbf{77}(1996) 490-493.
		
		
		
		
		\bibitem{2003Monneau} R. Monneau, On the number of singularities for the obstacle problem in two dimensions. \emph{J. Geom. Anal.,} \textbf{13}(2003), 359-389.
		
		\bibitem{2012PetrosyanShahgholianUraltseva}  A. Petrosyan, H. Shahgholian, N. N. Uraltseva, Regularity of free boundaries in obstacle-type problems. \emph{American Mathematical Society, Providence, RI,} 2012.
		
		\bibitem{2003RCFGKMWHV} Ch. R\"{u}egg, N. Cavadini, A. Furrer, H.-U. G\"{u}del, K. Kr\"{a}mer, H. Mutka, A. Wildes, K. Habicht, P. Vorderwischu, Bose-Einstein condensation of the triplet states in the magnetic insulator TlCuCl$_3$, \emph{Nature,} \textbf{423}(2003) 62-65.
		
		\bibitem{2019SpolaorVelichkov} L. Spolaor,  B. Velichkov,  An epiperimetric inequality for the regularity of some free boundary problems: the 2-dimensional case.  \emph{Comm. Pure Appl. Math.,} \textbf{72}(2019), 375-421.
		
		\bibitem{1964Stampacchia} G. Stampacchia, Formes bilin\'{e}aires coercitives sur les ensembles convexes. \emph{C. R. Acad. Sci. Paris,} \textbf{258}(1964), 4413-4416.
		
		\bibitem{2008struwe} M. Struwe, Variational methods. Applications to nonlinear partial differential equations and Hamiltonian systems. Fourth edition. \emph{Springer-Verlag, Berlin,} 2008.
		
		\bibitem{1999Weiss} G. S. Weiss, A homogeneity improvement approach to the obstacle problem. \emph{Invent. Math.,} \textbf{138}(1999), 23-50.
		
		\bibitem{1934Whitney} H. Whitney, Analytic extensions of differentiable functions defined in closed sets. \emph{Trans. Amer. Math. Soc.,} \textbf{36}(1934) 63--89.
		
		
		
	\end{thebibliography}
\end{document}